\documentclass[a4paper,fleqn]{cas-sc}

\usepackage[numbers,sort&compress]{natbib}

\usepackage{caption}
\usepackage{amssymb}
\usepackage{pifont}
\newcommand{\cmark}{\ding{51}}
\newcommand{\xmark}{\ding{55}}

\begin{document}
\let\WriteBookmarks\relax
\def\floatpagepagefraction{1}
\def\textpagefraction{.001}
\shorttitle{}
\shortauthors{}

\title [mode = title]{Developing solution algorithm for $LR$-type fully interval-valued intuitionistic fuzzy linear programming problems using lexicographic-ranking method}                      



\author[1]{Manisha Malik}
\ead{mmalik@ma.iitr.ac.in}

\credit{Conceptualization of this study, Methodology, Software, Validation, Formal analysis, Writing - original draft}

\address[1]{Department of Mathematics, Indian Institute of Technology Roorkee, 247 667, India}

\author[1]{S. K. Gupta}
\cormark[1]

\ead{s.gupta@ma.iitr.ac.in}

\credit{Visualization, Methodology, Investigation, Supervision, Writing - review $\&$ editing}

\author[2]{Manuel Arana-Jiménez}

\ead{manuel.arana@uca.es}
\address[2]{Department of Statistics and Operational Research, University of Cádiz, 11406, Spain}

\credit{Methodology, Investigation, Supervision}



\cortext[cor1]{Corresponding author.}


\begin{abstract}
In this article, a new concept of $LR$-type interval-valued intuitionistic fuzzy numbers ($LR$-type IVIFN) has been introduced. The theory has also been enriched by demonstrating diagrammatic representations of $LR$-type IVIFNs and establishing arithmetic operations among these fuzzy numbers. The total order properties of lexicographic criteria have been used for ranking $LR$-type IVIFNs. Further, a linear programming problem having both equality as well as inequality type constraints with all the parameters as $LR$-type IVIFNs and unrestricted decision variables has been formulated. An algorithm to find a unique optimal solution to the problem using the lexicographic ranking method has been developed. In the proposed methodology, the given linear programming problem is converted to an equivalent mixed $0$ \textendash $1$ lexicographic non-linear programming problem.  Various theorems have been proved to show the equivalence of the proposed problem and its different constructions. The model formulation, algorithm and discussed results  have not only developed a new idea but also generalized various well-known related works existing in the literature. A numerical problem has also been exemplified to show the steps involved in the approach. Finally, a practical application in production planning is framed, solved and analyzed to establish the applicability of the study.
\end{abstract}


%

%
%


\maketitle

\section{Introduction}

Linear programming problems (LPPs) are the simplest kind of optimization problem that are widely used to solve many real-life problems. In conventional LPPs, all the parameters and decision variables are taken to be precise real numbers. However, in practical situations due to various uncontrollable factors, the data may not be available as crisp values. It may involve some vagueness/ambiguity in all or some of the parameters and/or decision variables of the problem. A general fuzzy LPP can be modeled as follows:
\begin{center}
$\hspace{-4.39cm}\textbf{(P)}~\max~ (\mbox{or} \min)~ \displaystyle\sum_{j=1}^{n} c_j {x}_j$\\
$~~~~~~~~~~\mbox{subject to}~~\displaystyle\sum_{j=1}^{n} a_{ij}  {x}_j ~\{ \preceq,=,\succeq \}~ b_i,~~i=1,2, \dots ,m,$\\
$~~~~~~~~x_j \succeq 0,~~~~~~j=1,2, \dots ,n.$
\end{center}

Zadeh \cite{ref62} developed the concept of fuzzy sets which incorporated imprecision in the data successfully. Motivated by Zadeh's concept of fuzzy sets, Zimmermann \cite{ref64} initiated and developed the theory to solve fuzzy linear programming problems (FLPPs). Tanaka and Asai \cite{ref56} had first introduced the FLPPs in which both the parameters and decision variables were represented by fuzzy numbers. Initially, the researchers had extended the classical methods which were used to solve crisp LPPs to deal with FLPPs such as simplex algorithm, two-phase approach, etc. But, later on, these were proven to be incompatible with fuzzy theory. After that, linear ranking functions have been widely employed to convert FLPPs into crisp optimization problems. However, the ranking functions fail to order two such fuzzy numbers (FNs) which seems to be distinguished to a decision-maker. To overcome such limitations for the ordering of FNs, the idea of lexicographic ranking criteria came which uses multiple parameters at a time associated with an FN and hence is a more effective and powerful ordering criterion. The fuzzy theory was later on extended to the intuitionistic fuzzy set theory, which is more general than the former. The detailed literature can be seen in Section 2 of this article. It was observed that while handling the uncertain and hesitant data, the intuitionistic fuzzy sets utilize exact or crisp real numbers to assign a membership and non-membership degrees for each element of the set. However, in practical situations, a decision-maker may fails to give these degrees with full confidence. Consequently, the notion of intuitionistic fuzzy sets are extended to the interval-valued intuitionistic fuzzy (IVIF) sets, which is the key motivation for our present study. The IVIF sets use an interval to define the acceptance and rejection degrees of each element in the set. Moreover, to represent a realistic situation mathematically, $LR$-type fuzzy numbers play a crucial role in optimization theory since any type of variation in the input data can be reflected in the mathematical model using different $L$ and $R$ functions. Thus, we have firstly defined the concept of $LR$-type IVIF numbers and then considered an LPP having all the parameters and decision variables as $LR$-type IVIF numbers. Further, an approach for solving such LPPs using a lexicographic ranking methodology has been proposed. 

To the best of our knowledge, there is no study in the literature, describing the arithmetic operations on $LR$-type IVIF numbers (IVIFNs) and to find the unique optimal IVIF solution for an interval-valued intuitionistic fuzzy linear programming problem (IVIFLPP) having both linear equalities and inequalities with all the parameters represented by $LR$-type IVIFNs and decision variables as unrestricted $LR$-type IVIFNs. In many practical problems like selling / purchase of some units or profit / loss etc., unrestricted decision variables are required which can be handled through this formulation. To sum up, the key features of  the present work are listed as follows:
\begin{enumerate}
\item On the basis of $(\alpha, \beta)$-cut of $LR$-type IVIFNs, we define the score and accuracy indices of these numbers.
\item The basic arithmetic operations on unrestricted $LR$-type IVIFNs are developed using the $(\alpha, \beta)$-cut. 
\item Using the total order properties of the lexicographic criterion, a ranking of $LR$-type IVIFNs has been proposed.
\item Based on the introduced lexicographic ranking criterion, the $LR$-type IVIFLPP is converted to an equivalent mixed $0$\textendash$1$ lexicographic non-linear programming problem for finding a unique optimal IVIF solution of the $LR$-type fully IVIFLPP.
\item Various theorems are established to show the equivalence between the various problems obtained in the proposed algorithm.
\item A practical application in production planning is constructed, solved and examined using the proposed technique.
\end{enumerate}
        
The rest of the paper is summarized as follows: In Section 2, a detailed literature review on fuzzy, intuitionistic fuzzy and IVIF theory is given. Section 3 includes some basic definitions and arithmetic operations on $LR$-type IVIFNs. A lexicographic ranking criterion is also proposed to rank two $LR$-type IVIFNs. The mathematical formulation of an IVIFLPP is described in Section 4 and a lexicographic method has been proposed to find the unique IVIF optimal solution of $LR$-type IVIFLPPs. In Section 5, the advantages of the proposed method are listed. Section 6 illustrates a numerical example to describe the proposed algorithm. A production planning problem along with its managerial insights is discussed in Section 7. The last section sums up the conclusions and some interesting future directions.

\section{Literature review}

\subsection{Fuzzy LPP}

In the model (P), if all the parameters are taken to be fuzzy numbers (FNs), then the problem is described as a fuzzy linear programming problem (FLPP). In the literature, several methods have been proposed to solve these models depending on which parameters and/or decision variables are taken to be fuzzy. A comprehensive survey on FLPPs can be found in  Ebrahimnejad and Verdegay \cite[ch. 2-4]{ref17}. Hashemi et al. \cite{ref25} considered a fully FLPP with inequality constraints having all the parameters and decision variables to be given by symmetric $LR$-type FNs and proposed a two-phase solution approach by using the lexicographic comparison of the mean and standard deviation of FNs. Allahviranloo et al. \cite{ref3} used a ranking function to develop a solution algorithm to deal with the fully FLPPs having inequality constraints. Later on, Kumar et al. \cite{ref36} proposed a methodology to solve the fully FLPPs with equality constraints where parameters were taken to be unrestricted and decision variables as non-negative triangular FNs. After that, Najafi and Edalatpanah \cite{ref43} proposed some corrections to the methodology of Kumar et al. \cite{ref36}. Khan et al. \cite{ref34} studied a fully FLPP where parameters and decision variables were taken to be triangular FNs and proposed a method by making use of some ranking function. Ozkok et al. \cite{ref46} extended the method of Kumar and Kaur \cite{ref35} to solve fully FLPP with all types of constraints having parameters as unrestricted and decision variables as non-negative triangular FNs.

Najafi et al. \cite{ref44} examined a fully FLPP  having equality constraints with parameters as well as decision variables to be expressed by unrestricted triangular FNs and developed a solution technique by converting the original model to a non-linear model. Later, Gong and Zhao \cite{ref22} considered a fully FLPP with equality constraints and proposed a method in which the problem is first transformed into a crisp multi-objective LPP and then solved by using various approaches. Arana-Jiménez \cite{refarana} presented a new method to find fuzzy optimal (nondominated) solutions of fully FLPPs having inequality constraints with triangular fuzzy numbers and not necessarily symmetric, via solving a multiobjective linear problem with crisp numbers. Kaur and Kumar \cite{ref31} analyzed that by employing the existing approaches for solving fully FLPP, the obtained optimal solution is not necessarily unique. To overcome this limitation, they have defined a lexicographic criterion for ranking trapezoidal FNs and introduced an approach to find the unique optimal solution of fully FLPP having equality constraints with unrestricted parameters and non-negative decision variables. On  similar lines, Ezzati et al. \cite{ref18} introduced a lexicographic method to solve a fully FLPP with equality constraints having parameters as unrestricted triangular fuzzy and decision variables to be non-negative triangular fuzzy. Further, Mottaghi et al. \cite{ref41} solved a fully FLPP with inequality constraints by introducing non-negative fuzzy slack and surplus variables for converting the inequalities into fuzzy equality constraints.

 Based on a lexicographic criterion for ranking of $LR$-type FNs, Hosseinzadeh and Edalatpanah \cite{ref28} devised a method to solve a fully FLPP having only equality constraints where the parameters and decision variables were taken to be non-positive or non-negative $LR$-type FNs. After that, Kaur and Kumar \cite{ref33} introduced a lexicographic technique for obtaining the unique optimal solution of a fully FLPP with equality constraints having parameters and decision variables as unrestricted $LR$-type FNs. They pointed out that no method exists to obtain the unique optimal solution of a fully FLPP having inequalities in the set of constraints. But, some researchers \cite{ref21,ref22,ref41} solved fully FLPPs with inequality constraints by transforming them into equality constraints using fuzzy slack and surplus variables. However, in the case of FNs, such transformations are not correct mathematically and may lead to infeasible solutions for the considered FLPP. Later, Das et al. \cite{ref16} proposed a lexicographic method to solve a fully FLPP with all types of constraints keeping parameters as unrestricted and decision variables as non-negative trapezoidal FNs. But Ebrahimnejad and Verdegay \cite[p. 298]{ref17} demonstrated that this method is not suitable to deal with the fully FLPP having inequality constraints as the authors utilized different order relation for inequality constraints than that was used for the objective function, which is clearly false. Consequently, Ebrahimnejad and Verdegay \cite[p. 299]{ref17} suggested a correction by replacing the inequalities with a set of crisp linear inequalities. Pérez-Cañedo and Concepción-Morales \cite{ref47} introduced a method to solve a fully FLPP having equality and inequality constraints with parameters and decision variables as unrestricted $LR$-type FNs, using the lexicographic ranking criterion for the objective function and the set of inequality constraints. Recently, Tadesse et al. \cite{ref55} described a geometrical approach to handle the fully FLLP having non-negative decision variables. Further, some other significant applications of fuzzy theory can be found in studies of \cite{reffs1, reffs2, reffs3}.

\subsection{Intuitionistic fuzzy LPP}

Atanassov \cite{ref7} generalized Zadeh's concept of fuzzy sets by introducing intuitionistic fuzzy sets (IFSs) in order to include uncertainty as well as hesitation in the involved parameters. Angelov \cite{ref4} was the first to apply the IFS theory to optimization problems. Mahapatra and Roy \cite{ref40} developed the arithmetic operations on triangular intuitionistic fuzzy numbers (IFNs) and did reliability evaluation using these numbers. A linear programming problem (P) having equality and inequality constraints with all the parameters and decision variables expressed by IFNs is classified as a fully intuitionistic fuzzy linear programming problem (IFLPP). Nagoorgani and Ponnalagu \cite{ref42} proposed a method to solve an IFLPP with inequality constraints only. Using a ranking function for IFNs, Suresh et al. \cite{ref54} introduced a method to solve IFLPPs. Singh and Yadav \cite{ref52} suggested the modelling and optimization of the multi-objective non-linear programming problem in an intuitionistic fuzzy environment.

Later, Arefi and Taheri \cite{ref5} proposed the product of $LR$-type IFNs when both the numbers are either non-negative or non-positive or one is non-negative, and the other is non-positive. However, the remaining cases are not discussed. Then, Singh and Yadav \cite{ref53} introduced the product of unrestricted $LR$-type IFNs using $(\alpha, \beta)$-cut and proposed a method for solving fully IFLPPs using score and accuracy indices of $LR$-type IFNs. More review of IFS theory and its application to fully IFLPP can be seen in \cite{ref8, ref37, ref45, ref49, ref51, ref57}. Later on, Pérez-Cañedo and Concepción-Morales \cite{ref48} proposed a method using the total order properties of the lexicographic ranking criterion for finding the unique optimal intuitionistic fuzzy solution of a fully IFLPP having equality as well as inequality constraints with all the parameters and/or decision variables represented by unrestricted $LR$-type IFNs. Recently, Akram et al. \cite{ref2} introduced a class of fully Pythagorean fuzzy linear programming problems with equality constraints and suggested a linear ranking function based approach to handle such problems.

\subsection{Interval-valued intuitionistic fuzzy LPP}

 In view of real-life situations, it is more flexible and viable to represent the membership and non-membership degrees of an element by intervals rather than crisp real numbers. Hence, Atanassov and Gargov \cite{ref6} proposed the concept of interval-valued intuitionistic fuzzy (IVIF) sets. Optimizing a linear objective function over a set of linear constraints (P) where all the parameters and decision variables expressed by interval-valued intuitionistic fuzzy numbers (IVIFNs) is termed as a fully interval-valued intuitionistic fuzzy linear programming problem (IVIFLPP). Several researchers had  used the idea of IVIF theory for dealing with realistic decision-making problems. Ishibuchi and Tanaka \cite{ref30} were the first to solve a multi-objective programming problem in which coefficients of the objective function are intervals instead of crisp numbers. Şahin \cite{ref50} suggested a ranking of IVIFNs. The basic theory and various rankings of interval-valued fuzzy numbers can be reviewed in works of \cite{ref13, ref14, ref15, ref27, ref29, ref59}. Yang et al. \cite{ref61} had studied the combination of interval-valued fuzzy sets and soft sets.

Zhang et al. \cite{ref63} defined the $LR$-type interval-valued triangular FNs and proposed a method for solving multi-criteria decision-making problems with $LR$-type interval \textendash valued triangular fuzzy assessments and unknown weights. Garg et al. \cite{ref20} gave an intuitionistic fuzzy optimization approach using an interval environment to solve multi-objective reliability optimization problems.
 Later on, Akbari and Hesamian \cite{ref1} introduced signed-distance measures to rank $LR$-type interval-valued FNs and applied it to solve a multi-criteria group-decision making problem. Bharati and Singh \cite{ref11}  proposed a method to solve a multi-objective LPP in IVIF situations. Recently, Bharati and Singh \cite{ref12} introduced an approach for solving an IVIFLPP having unrestricted parameters while decision variables are taken to be non-negative.

A brief description of the various approaches to deal with FLPPs, IFLPPs and our proposed methodology is presented in Table \ref{table1}.

\begin{table*}[width=1\linewidth,cols=5,pos=h]
\caption{Existing approaches to solve FLPPs, IFLPPs and contribution of our present study} \label{table1}
\begin{scriptsize}
\begin{tabular*}{\tblwidth}{@{} LLCLL@{} }
\toprule
Existing methods & Type of FN/IFN/ & Unrestricted  & Criterion & Type of constraints \\
& IVIFN & variables & & (equality and/or inequality)\\
& & &  & \\
\midrule
Hashemi et al. \cite{ref25} & $LR$-type FN & \xmark & lexicographic & inequality only\\
Lotfi et al. \cite{ref39} & Triangular FN & \xmark & lexicographic & equality only \\
Kaur and Kumar \cite{ref31} & Trapezoidal FN & \xmark & lexicographic & equality only\\
Kaur and Kumar  \cite{ref32} & $LR$-type FN & \cmark & ranking function & equality and inequality both\\
Hosseinzadeh and Edalatpanah \cite{ref28} & $LR$-type FN & \xmark & lexicographic & equality only\\
Kaur and Kumar \cite{ref33} & $LR$-type FN & \cmark & lexicographic & equality only\\
Pérez-Cañedo and Concepción-Morales \cite{ref47} & $LR$-type FN & \cmark & lexicographic & equality and inequality both\\
Tadesse et al. \cite{ref55} & Triangular FN & \xmark & geometric approach & inequality only\\
Nagoorgani and Ponnalagu \cite{ref42} & Triangular IFN & \xmark & score function & inequality only\\
Singh and Yadav \cite{ref53} & $LR$-type IFN & \cmark & weighted sum of score  & equality and inequality both\\ 
 & & & and accuracy indices & \\
Pérez-Cañedo and Concepción-Morales \cite{ref48} & $LR$-type IFN & \cmark & lexicographic & equality and inequality both\\
Akram et al. \cite{ref2} & Pythagorean FN & \cmark & ranking function & equality only\\
Bharati and Singh \cite{ref12} & Triangular IVIFN & \xmark & expected value function & equality and inequality both\\
\textbf{Present study} & $LR$-type IVIFN & \cmark & lexicographic & equality and inequality both\\
\bottomrule
\end{tabular*}
\end{scriptsize}
\end{table*}

\section{Preliminaries}
In this section, we have introduced the basic concepts related to $LR$-type IVIFNs followed by the arithmetic operations on them.\\ 
 
\noindent{\bf{Definition 3.1}} \cite{ref6}. Let $X$ be the universal set and $Int[0,1]$ denote the set of all subintervals of the interval $[0,1]$. An interval-valued intuitionistic fuzzy set (IVIFS) is defined as a set $\tilde{A}=\{ (x, \mu_{\tilde{A}}(x), \nu_{\tilde{A}}(x)): x \in X \}$, where $\mu_{\tilde{A}} : X \rightarrow Int[0,1]$  and $\nu_{\tilde{A}} : X \rightarrow Int[0,1]$ represent the interval-valued membership and non-membership functions respectively, provided $~~0 \leq \text{Sup} (\mu_{\tilde{A}}(x)) +\text{Sup}( \nu_{\tilde{A}}(x)) \leq 1,$ $\forall~ x \in X$.\\

\noindent{\bf{Definition 3.2}} A set $\tilde{A}=\{ (x, \mu_{\tilde{A}}(x), \nu_{\tilde{A}}(x)): x \in X \}$ where $\mu_{\tilde{A}}=[\mu_{\tilde{A}}^{L},\mu_{\tilde{A}}^{U}]$ and $\nu_{\tilde{A}}=[\nu_{\tilde{A}}^{L},\nu_{\tilde{A}}^{U}]$ is called a convex IVIFS if $\forall~ x_1, x_2 \in X,~0 \leq \lambda \leq 1,$ the following conditions
are satisfied:
\begin{itemize}
\item $\mu_{\tilde{A}}^L(\lambda x_1+(1-\lambda) x_2) \geq \min \{ \mu_{\tilde{A}}^L(x_1), \mu_{\tilde{A}}^L(x_2)\},$
\item $\mu_{\tilde{A}}^U(\lambda x_1+(1-\lambda) x_2) \geq \min \{ \mu_{\tilde{A}}^U(x_1), \mu_{\tilde{A}}^U(x_2)\},$
\item $\nu_{\tilde{A}}^L(\lambda x_1+(1-\lambda) x_2) \leq \max \{\nu_{\tilde{A}}^L(x_1), \nu_{\tilde{A}}^L(x_2)\}$ and
\item $\nu_{\tilde{A}}^U(\lambda x_1+(1-\lambda) x_2) \leq \max \{\nu_{\tilde{A}}^U(x_1), \nu_{\tilde{A}}^U(x_2)\}.$
\end{itemize}

\noindent{\bf{Definition 3.3}} An IVIF set $\tilde{A}$ in $X$ is called normal IVIFS if there exist $x_1,x_2 \in X$ such that $\mu_{\tilde{A}}(x_1)=1$ and $\nu_{\tilde{A}}(x_2)=1$.\\  

\noindent{\bf{Definition 3.4}} An IVIF set $\tilde{A}=\{ (x, \mu_{\tilde{A}}(x), \nu_{\tilde{A}}(x)): x \in \mathbb{R} \}$ is called an IVIFN if the following conditions hold:
\begin{itemize}
\item $\tilde{A}$ is a convex IVIFS in $\mathbb{R}$,
\item $\tilde{A}$ is a normal IVIFS and
\item $\mu_{\tilde{A}}^L$, $\mu_{\tilde{A}}^U$, $\nu_{\tilde{A}}^L$ and $\nu_{\tilde{A}}^U$ are piecewise continuous functions from $\mathbb{R}$ to [0,1].
\end{itemize}
\noindent Mathematically, lower - upper membership and non-membership functions of an IVIFN $\tilde{A}$ can be represented as:
\begin{alignat*}{2}
    & \begin{aligned} 
     &\mu_{\tilde{A}}^{L}(x)=\begin{cases}
1,             & \mbox{if}  ~~  x=a,\\
g_1(x),   & \mbox{if}   ~~  a-l^{\mu}_L <x<a,\\
g_2(x),  & \mbox{if}    ~~ a<x<a+r^{\mu}_L,\\
0,            & \mbox{otherwise},
\end{cases}\\
\end{aligned}
    &  \hskip 2em  \hskip 3em & \mu_{\tilde{A}}^{U}(x)=
    \begin{aligned} 
    \begin{cases}
1,             & \mbox{if}  ~~  x=a,\\
h_1(x),   & \mbox{if}   ~~  a-l'^{\mu}_U<x<a,\\
h_2(x),  & \mbox{if}    ~~ a<x<a+r'^{\mu}_U,\\
0,            & \mbox{otherwise,}
\end{cases}\\[5.1ex]
  \end{aligned}
\end{alignat*}
\vspace{-1cm}
\begin{alignat*}{2}
    & \begin{aligned} 
     &\nu_{\tilde{A}}^{L}(x)=\begin{cases}
0,             & \mbox{if}  ~~  x=a,\\
l_1(x),   & \mbox{if}   ~~  a-l^{\nu}_L <x<a,\\
l_2(x),  & \mbox{if}    ~~ a<x<a+r^{\nu}_L,\\
1,            & \mbox{otherwise}
\end{cases}\\
\end{aligned}
    &  \hskip 2em \mbox{and} \hskip 2em & \nu_{\tilde{A}}^{U}(x)=  
     \begin{aligned} 
      \begin{cases}
0,             & \mbox{if}  ~~  x=a,\\
m_1(x),   & \mbox{if}   ~~  a-l'^{\nu}_U<x<a,\\
m_2(x),  & \mbox{if}    ~~ a<x<a+r'^{\nu}_U,\\
1,            & \mbox{otherwise}
\end{cases}\\[5.1ex]
  \end{aligned}
\end{alignat*}

where
\begin{enumerate}
\item [$(i).$] $g_1$, $h_1$, $l_2$ and $m_2$ are piecewise continuous and strictly increasing functions,
\item [$(ii).$]$g_2$, $h_2$, $l_1$ and $m_1$ are piecewise continuous and strictly decreasing functions,
\item [$(iii).$] $g_1(x) \leq h_1(x),~g_2(x) \leq h_2(x),~l_1(x) \leq m_1(x),~ l_2(x) \leq m_2(x),~~\forall~x \in \mathbb{R},$
\item [$(iv).$] $a$ is called the mean value of $\tilde{A}$,  
\item [$(v).$] $l^{\mu}_{L},~l'^{\mu}_{U},~l^{\nu}_{L}~\mbox{and}~l'^{\nu}_{U}$ are respectively the left spreads of $\mu_{\tilde{A}}^{L},~ \mu_{\tilde{A}}^{U}, ~\nu_{\tilde{A}}^{L}~\mbox{and}~\nu_{\tilde{A}}^{U}$ and 
\item [$(vi).$] $r^{\mu}_{L},~r'^{\mu}_{U},~ r^{\nu}_{L}~\mbox{and}~r'^{\nu}_{U}$ are respectively the right spreads of $\mu_{\tilde{A}}^{L},~ \mu_{\tilde{A}}^{U},~ \nu_{\tilde{A}}^{L}~\mbox{and}~\nu_{\tilde{A}}^{U}$.
\end{enumerate} 

\begin{figure*}
\centering
 \includegraphics[scale=0.5]{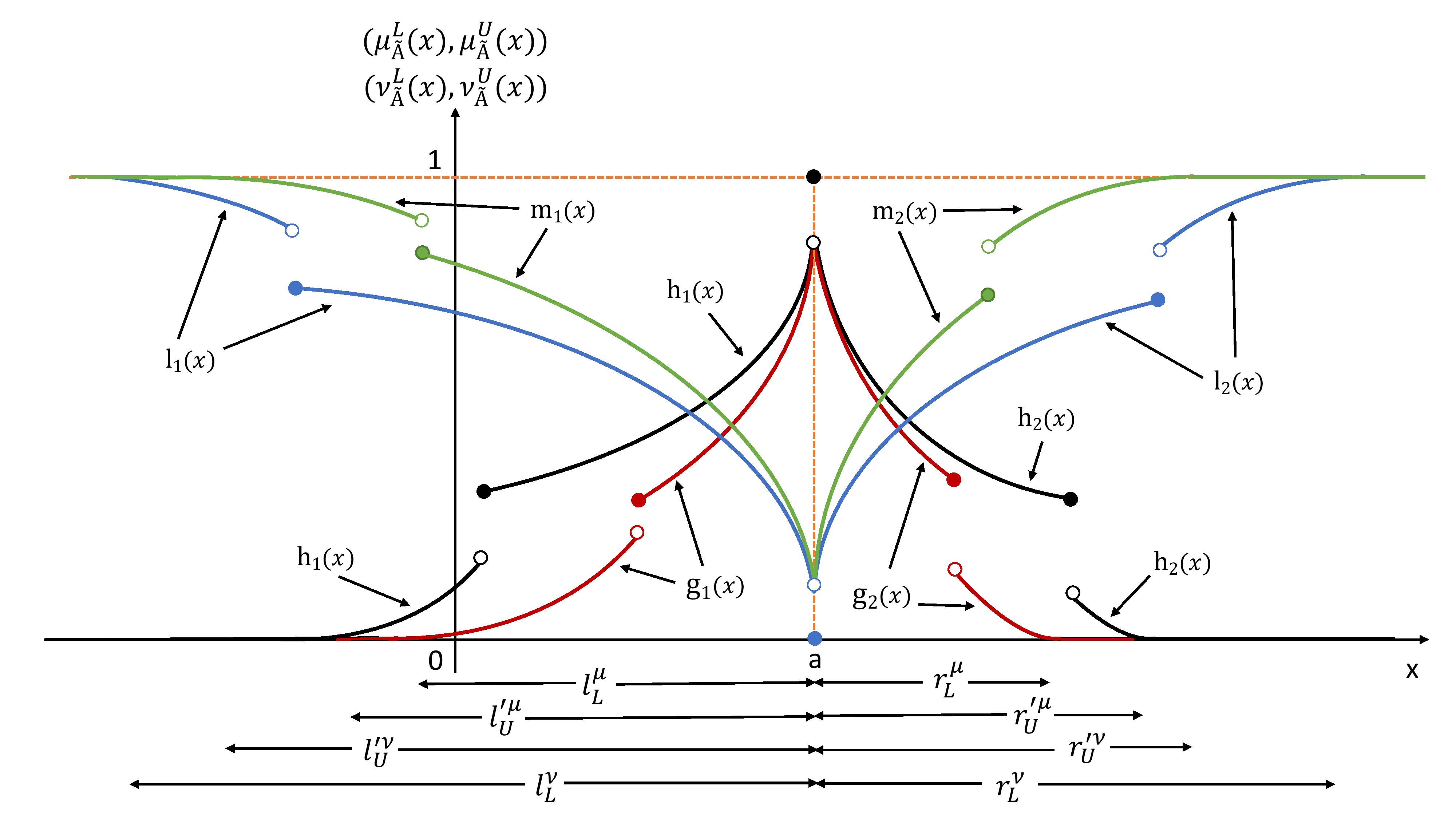}
\caption{Graphical representation of an IVIFN}
\label{FIG:1}
\end{figure*}

 It can be represented as $\tilde{A}=(a; l_{L}^{\mu},r_{L}^{\mu},l'^{\mu}_{U},r'^{\mu}_{U};l^{\nu}_{L},r^{\nu}_{L},l'^{\nu}_{U},r'^{\nu}_{U})$. The graphical representation of an IVIFN $\tilde{A}$ is given in Fig.~\ref{FIG:1}.\\

\begin{figure*}
\centering
\includegraphics[scale=0.5]{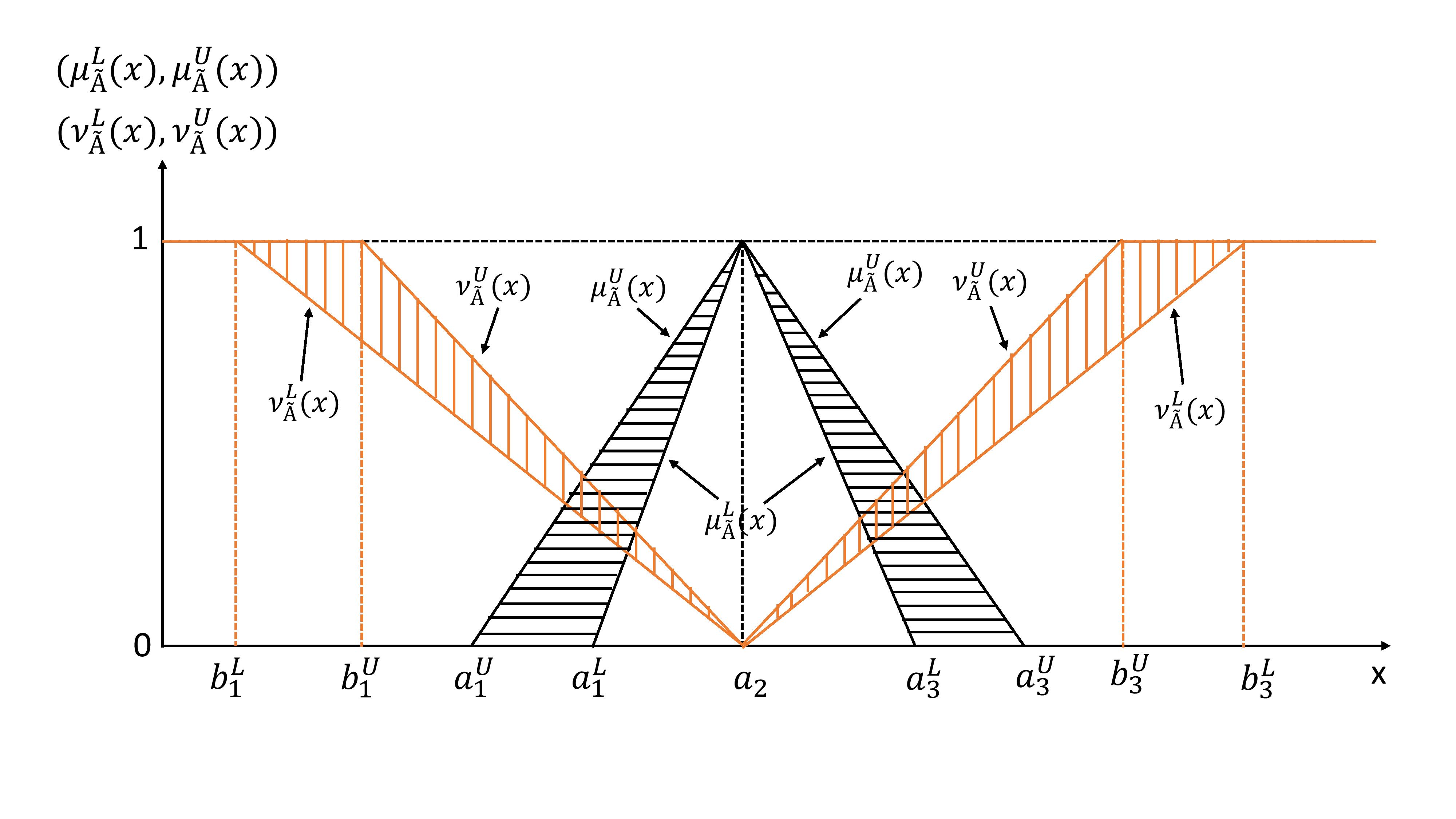}
\caption{Triangular interval-valued intuitionistic fuzzy number}\label{fig2}
\end{figure*}

\noindent{\bf{Definition 3.5}} \cite{ref10}. A triangular IVIFN (TIVIFN) is denoted by $\tilde{A}=\{ (a_1^U, a_1^L, a_2, a_3^L, a_3^U), (b_1^{L}, b_1^U, a_2, b_3^U, b_3^L) \}$, and its membership and non-membership degrees are defined as follows:
\begin{itemize}
\item Lower and upper membership functions are respectively given by:
\begin{alignat*}{2}
    & \begin{aligned} 
     &\mu_{\tilde{A}}^{L}(x)=\begin{cases}
1,             & \mbox{if}  ~~  x=a_2,\\
\displaystyle\frac{x-a_1^L}{a_2-a_1^L},   & \mbox{if}   ~~  a_1^L <x<a_2,\\
\displaystyle\frac{a_3^L-x}{a_3^L-a_2},  & \mbox{if}    ~~ a_2<x<a_3^L,\\
0,            & \mbox{otherwise}
\end{cases}\\
\end{aligned}
    &  \hskip 2em \mbox{and} \hskip 2em &\mu_{\tilde{A}}^{U}(x)=
    \begin{aligned} 
    \begin{cases}
1,             & \mbox{if}  ~~  x=a_2,\\
\displaystyle\frac{x-a_1^U}{a_2-a_1^U},   & \mbox{if}   ~~  a_1^U <x<a_2,\\
\displaystyle\frac{a_3^U-x}{a_3^U-a_2},  & \mbox{if}    ~~ a_2<x<a_3^U,\\
0,            & \mbox{otherwise}
\end{cases}\\[5.1ex]
  \end{aligned}
\end{alignat*}

\item Lower and upper non-membership functions are respectively defined as:
\begin{alignat*}{2}
    & \begin{aligned} 
     &\nu_{\tilde{A}}^{L}(x)=\begin{cases}
\medskip
0,             & \mbox{if}  ~~  x=a_2,\\
\medskip
\displaystyle\frac{a_2-x}{a_2-b_1^L},   & \mbox{if}   ~~  b_1^L <x<a_2,\\
\displaystyle\frac{a_2-x}{a_2-b_3^L},  & \mbox{if}    ~~ a_2<x<b_3^L,\\
1,            & \mbox{otherwise}
\end{cases}\\
\end{aligned}
    &  \hskip 2em \mbox{and} \hskip 2em &\nu_{\tilde{A}}^{U}(x)=
       \begin{aligned}  \begin{cases}
 \medskip
0,             & \mbox{if}  ~~  x=a_2,\\
\medskip
\displaystyle\frac{x-a_2}{b_1^U-a_2},   & \mbox{if}   ~~  b_1^U <x<a_2,\\
\displaystyle\frac{x-a_2}{b_3^U-a_2},  & \mbox{if}    ~~ a_2<x<b_3^U,\\
1,            & \mbox{otherwise}
\end{cases}\\[5.1ex]
  \end{aligned}
\end{alignat*}
\end{itemize}
\vspace{-0.4cm}
\noindent where ~$b_{1}^{L}\leq b_{1}^{U}\leq a_{1}^{U}\leq a_{1}^{L}\leq a_{2}\leq a_{3}^{L}\leq a_{3}^{U}\leq b_{3}^{U}\leq b_{3}^{L}$. The diagrammatic representation of a TIVIFN is shown in Fig. \ref{fig2}.\\

\noindent{\bf{Definition 3.6}} \cite{ref53}. A function $f: [0, \infty) \to [0, 1]$ is said to be shape function or reference function if it satisfies the following conditions:
\begin{enumerate}
\item [$(i)$] $f(0) = 1$,
\item [$(ii)$] $f$ is invertible on $[0, \infty)$,
\item [$(iii)$] $f$ is continuous function on $[0, \infty)$,
\item [$(iv)$] $f$ is strictly decreasing on $[0, \infty)$  and
\item [$(v)$] $\displaystyle \lim_{x \to \infty} f(x)=0.$
\end{enumerate}

\begin{figure*}
\centering
\includegraphics[scale=0.5]{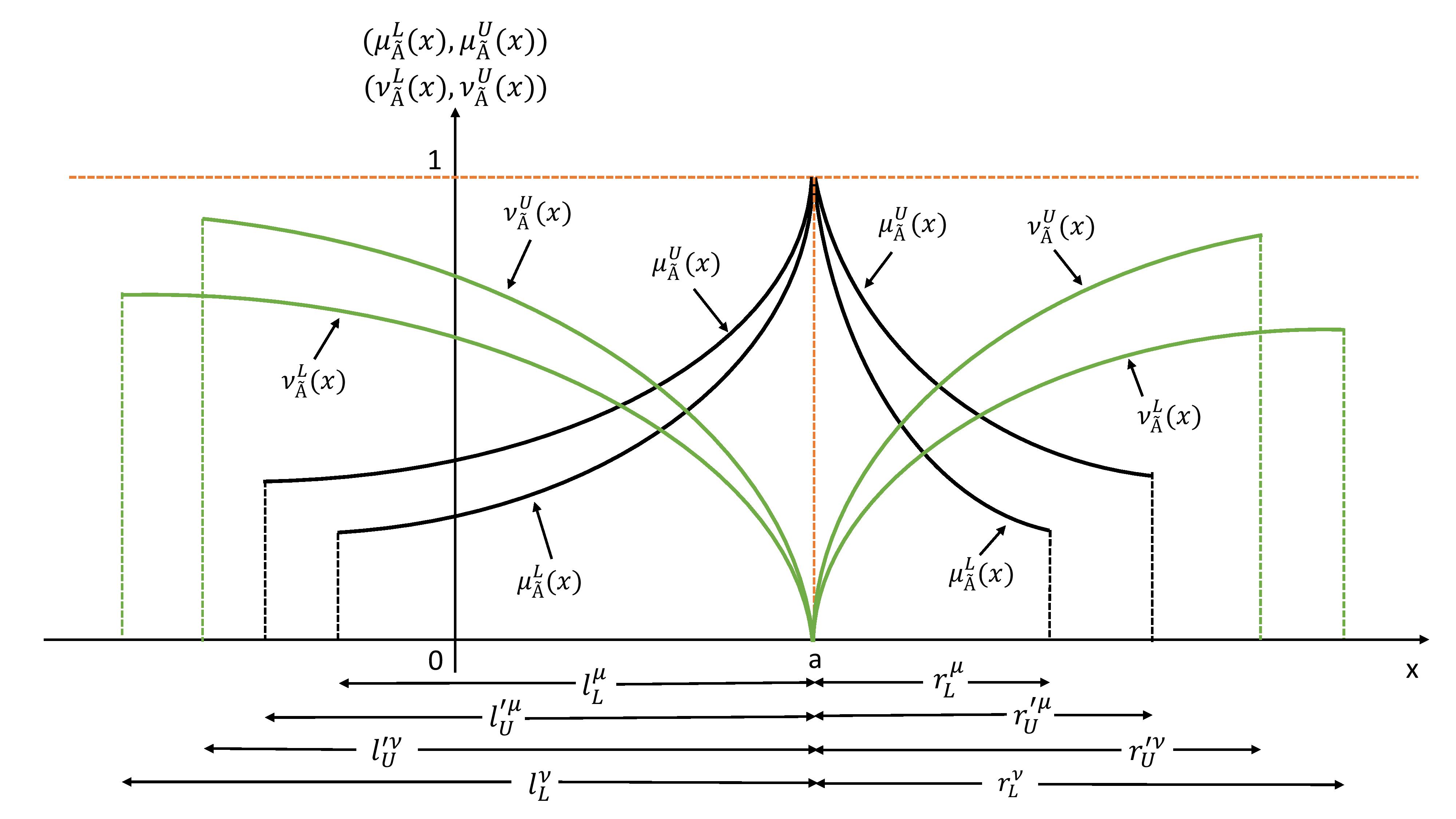}
\caption{Graphical representation of an $LR$-type IVIFN}
\label{fig3}
\end{figure*}

\noindent{\bf{Definition 3.7}} An IVIFN $\tilde{A}$ is said to be $LR$-type IVIFN if there exist shape functions $L, R, L'~ \mbox{and}~ R'$, and positive real constants $l^{\mu}_{L},~r^{\mu}_{L},~l'^{\mu}_{U} ,~r'^{\mu}_{U} ,~l^{\nu}_{L} ,~r^{\nu}_{L} ,~l'^{\nu}_{U}~\mbox{and}~r'^{\nu}_{U}$, such that its
\begin{itemize}
\item Lower and upper membership functions, respectively are defined as:
\begin{alignat*}{2}
    & \begin{aligned} 
     &\mu_{\tilde{A}}^{L}(x)=\begin{cases}
L\Bigg(\displaystyle\frac{a-x}{l^{\mu}_{L}}\Bigg),             & a-l^{\mu}_{L} \leq x \leq a,\\\\
R\Bigg(\displaystyle\frac{x-a}{r^{\mu}_{L}}\Bigg),             & a \leq x \leq a+r^{\mu}_{L},\\
~0, & \mbox{otherwise}
\end{cases}\\
\end{aligned}
    &  \hskip 2em \mbox{and} \hskip 2em & \mu_{\tilde{A}}^{U}(x)=
    \begin{aligned} 
    \begin{cases}
L'\Bigg(\displaystyle\frac{a-x}{l'^{\mu}_{U}}\Bigg),             & a-l'^{\mu}_{U} \leq x \leq a,\\\\
R'\Bigg(\displaystyle\frac{x-a}{r'^{\mu}_{U}}\Bigg),             & a \leq x \leq a+r'^{\mu}_{U},\\
~0, & \mbox{otherwise}
\end{cases}\\[5.1ex]
  \end{aligned}
\end{alignat*}

\item Lower and upper non-membership functions, respectively are given by:
\begin{alignat*}{2}
    & \begin{aligned} 
     &\nu_{\tilde{A}}^{L}(x)=\begin{cases}
1-L\Bigg(\displaystyle\frac{a-x}{l^{\nu}_{L}}\Bigg),             & a-l^{\nu}_{L} \leq x \leq a,\\\\
1-R\Bigg(\displaystyle\frac{x-a}{r^{\nu}_{L}}\Bigg),             & a \leq x \leq a+r^{\nu}_{L},\\
~1, & \mbox{otherwise}
\end{cases}\\
\end{aligned}
    &  \hskip 2em \mbox{and} \hskip 2em & \nu_{\tilde{A}}^{U}(x)=
       \begin{aligned} 
       \begin{cases}
1-L'\Bigg(\displaystyle\frac{a-x}{l'^{\nu}_{U}}\Bigg),             & a-l'^{\nu}_{U} \leq x \leq a,\\\\
1-R'\Bigg(\displaystyle\frac{x-a}{r'^{\nu}_{U}}\Bigg),             & a \leq x \leq a+r'^{\nu}_{U},\\
~1, & \mbox{otherwise}\end{cases}\\[5.1ex]
  \end{aligned}
\end{alignat*}
\end{itemize}

 where $l'^{\mu}_{U} \geq l^{\mu}_{L} ,~~r'^{\mu}_{U} \geq r^{\mu}_{L} ,~~l^{\nu}_{L} \geq l'^{\nu}_{U} ,~~r^{\nu}_{L} \geq r'^{\nu}_{U} ,~~l^{\nu}_{L} \geq l^{\mu}_{L} ,~r^{\nu}_{L} \geq r^{\mu}_{L} ,~~l'^{\nu}_{U} \geq l'^{\mu}_{U} ,~~r'^{\nu}_{U} \geq r'^{\mu}_{U}$ and\\ $ 0 \leq \mbox{Sup}\{\mu_{\tilde{A}}(x)\}+\mbox{Sup}\{\nu_{\tilde{A}}(x)\} \leq 1,~\forall ~x \in \mathbb{R}$. $a$ is called the mean value of $\tilde{A}$; $l^{\mu}_{L},~l'^{\mu}_{U},~l^{\nu}_{L}~\mbox{and}~l'^{\nu}_{U}$ are respectively the left spreads of $\mu_{\tilde{A}}^{L},~ \mu_{\tilde{A}}^{U}, ~\nu_{\tilde{A}}^{L}~\mbox{and}~\nu_{\tilde{A}}^{U}$, and $r^{\mu}_{L},~r'^{\mu}_{U},~ r^{\nu}_{L}~\mbox{and}~r'^{\nu}_{U}$ are the respective right spreads of $\mu_{\tilde{A}}^{L},~ \mu_{\tilde{A}}^{U},~ \nu_{\tilde{A}}^{L}~\mbox{and}~\nu_{\tilde{A}}^{U}$. An $LR$-type IVIFN is denoted by $\tilde{A}=(a; l_{L}^{\mu},r_{L}^{\mu},l'^{\mu}_{U},r'^{\mu}_{U};l^{\nu}_{L},r^{\nu}_{L},l'^{\nu}_{U},r'^{\nu}_{U})_{LR}$ and its possible general graphical representation is shown in Fig. \ref{fig3}. Let $IV(\mathbb{R})$ represents the set of all $LR$-type IVIFNs.\\

\noindent{\bf{Remark 3.1}} Taking $L(x)=R(x)=L'(x)=R'(x)=\max \{0, 1-x \},~\forall~ x \in \mathbb{R}$, the Definition $3.7$ reduces to Definition $3.5$.\\

%
%

\noindent{\bf{Definition 3.8}} An $LR$-type IVIFN $\tilde{A}=(a; l_{L}^{\mu},r_{L}^{\mu},l'^{\mu}_{U},r'^{\mu}_{U};l^{\nu}_{L},r^{\nu}_{L},l'^{\nu}_{U},r'^{\nu}_{U})_{LR}$ is called an unrestricted $LR$-type IVIFN if $a$ is any real number.\\

\noindent{\bf{Definition 3.9}} An $LR$-type IVIFN $\tilde{A}=(a; l_{L}^{\mu},r_{L}^{\mu},l'^{\mu}_{U},r'^{\mu}_{U};l^{\nu}_{L},r^{\nu}_{L},l'^{\nu}_{U},r'^{\nu}_{U})_{LR}$ is called non-negative (positive) if\\ $a-l^{\nu}_{L} \geq (>)~0$ and non-positive (negative) if $a+r^{\nu}_{L} \leq(<)~ 0$.\\


\noindent{\bf{Theorem 3.1.}} \textit{Let $\tilde{A}=(a; l_{L}^{\mu},r_{L}^{\mu},l'^{\mu}_{U},r'^{\mu}_{U};l^{\nu}_{L},r^{\nu}_{L},l'^{\nu}_{U},r'^{\nu}_{U})_{LR}$ be an $LR$-type IVIFN. Then, $\forall~\alpha,\beta \in (0,1]~ \mbox{and} ~\alpha+\beta \leq 1,$} 
\begin{enumerate}
\item[$(i).$] \textit{its lower $\alpha$-cut for membership and lower $\beta$-cut for non-membership are respectively, given by\\\\
\noindent$A_{\alpha}^{L}= [a-l^{\mu}_{L}L^{-1}(\alpha), a+r^{\mu}_{L}R^{-1}(\alpha)]~~~\mbox{and}~~~A_{\beta}^{L}= [a-l^{\nu}_{L}L^{-1}(1-\beta), a+r^{\nu}_{L}R^{-1}(1-\beta)].$} \\
\item[$(ii).$] \textit{its upper $\alpha$-cut for membership and upper $\beta$-cut for non-membership, respectively are\\\\
$A_{\alpha}^{U}= [a-l'^{\mu}_{U}(L')^{-1}(\alpha), a+r'^{\mu}_{U}(R')^{-1}(\alpha)]~~~\mbox{and}~~~A_{\beta}^{U}= [a-l'^{\nu}_{U}(L')^{-1}(1-\beta), a+r'^{\nu}_{U}(R')^{-1}(1-\beta)].$} \\
\item[$(iii).$] \textit{its lower and upper $(\alpha, \beta)$-cut respectively, are\\\\
$A_{\alpha, \beta}^{L}= [a-l^{\mu}_{L}L^{-1}(\alpha), a+r^{\mu}_{L}R^{-1}(\alpha)]~ \cap ~[a-l^{\nu}_{L}L^{-1}(1-\beta), a+r^{\nu}_{L}R^{-1}(1-\beta)]~~\mbox{and} $\\\\
\noindent  $A_{\alpha, \beta}^{U}=[a-l'^{\mu}_{U}(L')^{-1}(\alpha), a+r'^{\mu}_{U}(R')^{-1}(\alpha)]~ \cap ~[a-l'^{\nu}_{U}(L')^{-1}(1-\beta), a+r'^{\nu}_{U}(R')^{-1}(1-\beta)].$}
\end{enumerate}
\noindent{\bf{Proof.}} \begin{enumerate}
\item[$(i).$] For $\alpha \in (0,1]$, $~\mu_{\tilde{A}}^{L}(x) \geq  \alpha ~~~\mbox{implies}$
$$L\Bigg(\displaystyle\frac{a-x}{l^{\mu}_{L}}\Bigg) \geq \alpha~\mbox{and}~R\Bigg(\displaystyle\frac{x-a}{r^{\mu}_{L}}\Bigg) \geq \alpha.$$

\noindent $\mbox{Since $L$ and $R$ are decreasing functions, therefore}$
$$\displaystyle\frac{a-x}{l^{\mu}_{L}} \leq L^{-1}(\alpha),~\displaystyle\frac{x-a}{r^{\mu}_{L}} \leq R^{-1}(\alpha).$$
\noindent$\mbox{It further yields}$
$$a-l^{\mu}_{L}L^{-1}(\alpha) \leq x \leq a+r^{\mu}_{L}R^{-1}(\alpha).$$

\noindent Hence, 
$$A_{\alpha}^{L}= [a-l^{\mu}_{L}L^{-1}(\alpha), a+r^{\mu}_{L}R^{-1}(\alpha)].$$

\noindent Now, for $\beta \in (0,1]~ \mbox{such that}~ \alpha+\beta \leq 1,$\\

$\nu_{\tilde{A}}^{L}(x) \leq  \beta~~\mbox{gives}$
$$1-L\Bigg(\displaystyle\frac{a-x}{l^{\nu}_{L}}\Bigg) \leq \beta~\mbox{and}~1-R\Bigg(\displaystyle\frac{x-a}{r^{\nu}_{L}}\Bigg) \leq \beta$$

\noindent$\mbox{which implies}$
$$\displaystyle\frac{a-x}{l^{\nu}_{L}} \leq L^{-1}(1-\beta),~\displaystyle\frac{x-a}{r^{\nu}_{L}} \leq R^{-1}(1-\beta).$$

\noindent Thus, 
$$A_{\beta}^{L}= [a-l^{\nu}_{L}L^{-1}(1-\beta), a+r^{\nu}_{L}R^{-1}(1-\beta)].$$
\medspace
\noindent This proves $(i)$.
\item[$(ii).$] Applying $\alpha$-cut on the upper membership function, that is, $\mu_{\tilde{A}}^{U}(x) \geq  \alpha$, $\alpha \in (0,1]$, we get
$$ L'\Bigg(\displaystyle\frac{a-x}{l'^{\mu}_{U}}\Bigg) \geq \alpha~\mbox{and}~R'\Bigg(\displaystyle\frac{x-a}{r'^{\mu}_{U}}\Bigg) \geq \alpha.$$

\noindent Using the fact that $L^{\prime}$ and $R^{\prime}$ are decreasing functions, it follows that
$$\displaystyle\frac{a-x}{l'^{\mu}_{U}} \leq (L')^{-1}(\alpha)~\mbox{and}~\displaystyle\frac{x-a}{r'^{\mu}_{U}} \leq (R')^{-1}(\alpha).$$

\noindent This after simplification gives
$$a-l'^{\mu}_{U}(L')^{-1}(\alpha) \leq x \leq a+r'^{\mu}_{U}(R')^{-1}(\alpha).$$

\noindent Therefore,
$$A_{\alpha}^{U}= [a-l'^{\mu}_{U}(L')^{-1}(\alpha), a+r'^{\mu}_{U}(R')^{-1}(\alpha)].$$

\noindent Similarly, for $\beta \in (0,1]~\mbox{such that}~ \alpha+\beta \leq 1,$ the expression $\nu_{\tilde{A}}^{U}(x) \leq  \beta$ yields
$$1-L'\Bigg(\displaystyle\frac{a-x}{l'^{\nu}_{U}}\Bigg) \leq \beta~\mbox{and}~1-R'\Bigg(\displaystyle\frac{x-a}{r'^{\nu}_{U}}\Bigg) \leq \beta.$$

\noindent This finally gives
$$A_{\beta}^{U}= [a-l'^{\nu}_{U}(L')^{-1}(1-\beta), a+r'^{\nu}_{U}(R')^{-1}(1-\beta)].$$
\medspace
\noindent Hence proved part $(ii)$.
\item[$(iii).$] From $(i)$, the lower $\alpha$-cut for membership and the lower $\beta$-cut for non-membership of $\tilde{A}$ are respectively, given by 
$$A_{\alpha}^{L}= [a-l^{\mu}_{L}L^{-1}(\alpha), a+r^{\mu}_{L}R^{-1}(\alpha)]~~~~\mbox{and}~~~A_{\beta}^{L}= [a-l^{\nu}_{L}L^{-1}(1-\beta), a+r^{\nu}_{L}R^{-1}(1-\beta)].$$

\noindent It yields\\
$$\hspace{-9.65cm}A_{\alpha, \beta}^{L}=A_{\alpha}^{L} \cap A_{\beta}^{L}$$
$$= [a-l^{\mu}_{L}L^{-1}(\alpha), a+r^{\mu}_{L}R^{-1}(\alpha)] ~\cap~ [a-l^{\nu}_{L}L^{-1}(1-\beta), a+r^{\nu}_{L}R^{-1}(1-\beta)].$$

\noindent On the same lines, the proof of $A_{\alpha, \beta}^{U}$ can also be obtained. Hence, the result. 
\end{enumerate}


\noindent{\bf{Definition 3.10}} Let $\tilde{A}=(a; l_{L}^{\mu},r_{L}^{\mu},l'^{\mu}_{U},r'^{\mu}_{U};l^{\nu}_{L},r^{\nu}_{L},l'^{\nu}_{U},r'^{\nu}_{U})_{LR}$ be an $LR$-type IVIFN. Then, the score and accuracy indices of $\tilde{A}$ are denoted by $S(\tilde{A})$ and $A(\tilde{A})$, respectively and are defined by:\\
$S(\tilde{A}):=\displaystyle\frac{1}{4}\displaystyle\int_{0}^{1} \big(a-l^{\mu}_{L}L^{-1}(\alpha)+a+r^{\mu}_{L}R^{-1}(\alpha)+a-l'^{\mu}_{U}(L')^{-1}(\alpha)+a+r'^{\mu}_{U}(R')^{-1}(\alpha)\big) d\alpha -\displaystyle\frac{1}{4}\displaystyle\int_{0}^{1} \big(a-l^{\nu}_{L}L^{-1}(1-\beta)$\\

$~~~~~~~~~+a+r^{\nu}_{L}R^{-1}(1-\beta)+a-l'^{\nu}_{U}(L')^{-1}(1-\beta)+a+r'^{\nu}_{U}(R')^{-1}(1-\beta)\big) d\beta $\\

\noindent $A(\tilde{A}):=\displaystyle\frac{1}{4}\displaystyle\int_{0}^{1} \big(a-l^{\mu}_{L}L^{-1}(\alpha)+a+r^{\mu}_{L}R^{-1}(\alpha)+a-l'^{\mu}_{U}(L')^{-1}(\alpha)+a+r'^{\mu}_{U}(R')^{-1}(\alpha)\big) d\alpha +\displaystyle\frac{1}{4}\displaystyle\int_{0}^{1} \big(a-l^{\nu}_{L}L^{-1}(1-\beta)$\\

$~~~~~~~~~+a+r^{\nu}_{L}R^{-1}(1-\beta)+a-l'^{\nu}_{U}(L')^{-1}(1-\beta)+a+r'^{\nu}_{U}(R')^{-1}(1-\beta)\big) d\beta $\\

\noindent{\bf{Remark 3.2}} If $l_{L}^{\mu}= l'^{\mu}_{U},~r_{L}^{\mu}= r'^{\mu}_{U},~l^{\nu}_{L}=l'^{\nu}_{U}$ and $r^{\nu}_{L}=r'^{\nu}_{U}$, then the Definition 3.10 reduces to the corresponding definition for $LR$-type IFNs given in Singh and Yadav \cite{ref53}. \\

\noindent{\bf{Theorem 3.2.}} \textit{Let $\tilde{A}=(a; a-a_1^L, a_3^L-a, a-a_1^U, a_3^U-a; a-b_1^{L}, b_3^L-a, a-b_1^U, b_3^U-a)_{LR}$ be an $LR$-type TIVIFN. Then, the score and accuracy indices of $LR$-type TIVIFN $\tilde{A}$ are respectively, given by:
$$S(\tilde{A})=\displaystyle\frac{a_1^L+a_3^L+a_1^U+a_3^U-b_1^L-b_3^L-b_1^U-b_3^U}{8},$$ 
$$A(\tilde{A})=\displaystyle\frac{a_1^L+a_3^L+a_1^U+a_3^U+8a+b_1^L+b_3^L+b_1^U+b_3^U}{8}.$$} 
\noindent{\bf{Proof.}} From Definition 3.10, we have\\
$S(\tilde{A}):=\displaystyle\frac{1}{4}\displaystyle\int_{0}^{1} \big(a-l^{\mu}_{L}L^{-1}(\alpha)+a+r^{\mu}_{L}R^{-1}(\alpha)+a-l'^{\mu}_{U}(L')^{-1}(\alpha)+a+r'^{\mu}_{U}(R')^{-1}(\alpha)\big) d\alpha -\displaystyle\frac{1}{4}\displaystyle\int_{0}^{1} \big(a-l^{\nu}_{L}L^{-1}(1-\beta)$
\begin{equation}
\hspace{0.5cm}+a+r^{\nu}_{L}R^{-1}(1-\beta)+a-l'^{\nu}_{U}(L')^{-1}(1-\beta)+a+r'^{\nu}_{U}(R')^{-1}(1-\beta)\big) d\beta 
\end{equation}

\noindent $A(\tilde{A}):=\displaystyle\frac{1}{4}\displaystyle\int_{0}^{1} \big(a-l^{\mu}_{L}L^{-1}(\alpha)+a+r^{\mu}_{L}R^{-1}(\alpha)+a-l'^{\mu}_{U}(L')^{-1}(\alpha)+a+r'^{\mu}_{U}(R')^{-1}(\alpha)\big) d\alpha +\displaystyle\frac{1}{4}\displaystyle\int_{0}^{1} \big(a-l^{\nu}_{L}L^{-1}(1-\beta)$
\begin{equation}
\hspace{0.5cm}+a+r^{\nu}_{L}R^{-1}(1-\beta)+a-l'^{\nu}_{U}(L')^{-1}(1-\beta)+a+r'^{\nu}_{U}(R')^{-1}(1-\beta)\big) d\beta \end{equation}

\noindent Now, since $\tilde{A}$ is a TIVIFN, therefore
$$L(x)=R(x)=L'(x)=R'(x)=\max \{0, 1-x \},~\forall~x \in \mathbb{R}.$$
Hence, for $\alpha \in (0,1]$, we have
\begin{equation} \hspace{4.5cm} L(\alpha)=L^{\prime}(\alpha)=R(\alpha)=R^{\prime}(\alpha)=1-\alpha.\end{equation}
This further implies
\begin{equation} \hspace{4cm} L^{-1}(\alpha)=R^{-1}(\alpha)=(L')^{-1}(\alpha)=(R')^{-1}(\alpha)=1-\alpha. \end{equation}

\noindent Substituting the expressions from the equations $(3)$ and $(4)$ in $(1)$ and $(2)$, we obtain\\ 
$$S(\tilde{A})=\displaystyle\frac{a_1^L+a_3^L+a_1^U+a_3^U-b_1^L-b_3^L-b_1^U-b_3^U}{8}~~~~\mbox{and}$$ 
$$A(\tilde{A})=\displaystyle\frac{a_1^L+a_3^L+a_1^U+a_3^U+8a+b_1^L+b_3^L+b_1^U+b_3^U}{8}.$$  
Hence the result.\\

\subsection{Arithmetic operations on LR-type IVIFNs}
In this subsection, the basic arithmetic operations on $LR$-type IVIFNs are discussed. Here, we have introduced the addition operator $(\oplus)$, subtraction operator $(\ominus)$ and product operator $(\odot)$ for $LR$-type IVIFNs. The following propositions discuss the detailed expressions for the addition, subtraction, scalar multiplication and product operations on these numbers.\\

\noindent{\bf{Proposition 3.1.1.}} \textit{Let $\tilde{A}_1=(a_1; l_{1L}^{\mu},r_{1L}^{\mu},l'^{\mu}_{1U},r'^{\mu}_{1U};l^{\nu}_{1L},r^{\nu}_{1L},l'^{\nu}_{1U},r'^{\nu}_{1U})_{LR}$ and $\tilde{A}_2=(a_2; l_{2L}^{\mu},r_{2L}^{\mu},l'^{\mu}_{2U},r'^{\mu}_{2U};l^{\nu}_{2L},r^{\nu}_{2L},\\l'^{\nu}_{2U},r'^{\nu}_{2U})_{LR}$ be two $LR$-type IVIFNs. Then,}
\begin{enumerate}
\item [$(i).$] $\tilde{A}_1\oplus\tilde{A}_2=(a_1+a_2; l_{1L}^{\mu}+l_{2L}^{\mu},r_{1L}^{\mu}+r_{2L}^{\mu},l'^{\mu}_{1U}+l'^{\mu}_{2U},r'^{\mu}_{1U}+r'^{\mu}_{2U};l^{\nu}_{1L}+l^{\nu}_{2L},r^{\nu}_{1L}+r^{\nu}_{2L},l'^{\nu}_{1U}+l'^{\nu}_{2U},r'^{\nu}_{1U}+r'^{\nu}_{2U})_{LR},$\\\\
\noindent \textit{where the conditions for $LR-$type representation of $\tilde{A}_1\oplus\tilde{A}_2$ are satisfied.}\\
\item [$(ii).$] $\tilde{A}_1\ominus\tilde{A}_2=(a_1-a_2; l_{1L}^{\mu}+r_{2L}^{\mu},r_{1L}^{\mu}+l_{2L}^{\mu},l'^{\mu}_{1U}+r'^{\mu}_{2U},r'^{\mu}_{1U}+l'^{\mu}_{2U};l^{\nu}_{1L}+r^{\nu}_{2L},r^{\nu}_{1L}+l^{\nu}_{2L},l'^{\nu}_{1U}+r'^{\nu}_{2U},r'^{\nu}_{1U}+l'^{\nu}_{2U})_{LR},$\\

\noindent \textit{where the conditions for $LR$-type representation of $\tilde{A}_1\ominus \tilde{A}_2$ are fulfilled.}
\end{enumerate}

\noindent{\bf{Proof.}} In view of Theorem 3.1, the $\alpha$ and $\beta$-cuts of $\tilde{A}_1 ~\mbox{and} ~\tilde{A}_2$ are respectively, given by
\begin{equation}
 \renewcommand{\arraystretch}{.60}
 \hspace{2.5cm} \left.\begin{array}{r@{\;}l} 
A_{1\alpha}^{L}= [a_1-l^{\mu}_{1L}L^{-1}(\alpha), a_1+r^{\mu}_{1L}R^{-1}(\alpha)],\\\\

A_{1\alpha}^{U}= [a_1-l'^{\mu}_{1U}(L')^{-1}(\alpha), a_1+r'^{\mu}_{1U}(R')^{-1}(\alpha)],\\\\

 A_{1\beta}^{L}= [a_1-l^{\nu}_{1L}L^{-1}(1-\beta), a_1+r^{\nu}_{1L}R^{-1}(1-\beta)],\\\\
 
A_{1\beta}^{U}= [a_1-l'^{\nu}_{1U}(L')^{-1}(1-\beta), a_1+r'^{\nu}_{1U}(R')^{-1}(1-\beta)].
\end{array} \right\} \label{2}
\end{equation}
\vspace{-0.5cm} 
\begin{equation}
 \renewcommand{\arraystretch}{.60}
 \hspace{2.5cm}  \left.\begin{array}{r@{\;}l}
A_{2\alpha}^{L}= [a_2-l^{\mu}_{2L}L^{-1}(\alpha), a_2+r^{\mu}_{2L}R^{-1}(\alpha)],\\\\

A_{2\alpha}^{U}= [a_2-l'^{\mu}_{2U}(L')^{-1}(\alpha), a_2+r'^{\mu}_{2U}(R')^{-1}(\alpha)],\\\\

A_{2\beta}^{L}= [a_2-l^{\nu}_{2L}L^{-1}(1-\beta), a_2+r^{\nu}_{2L}R^{-1}(1-\beta)],\\\\
 
A_{2\beta}^{U}= [a_2-l'^{\nu}_{2U}(L')^{-1}(1-\beta), a_2+r'^{\nu}_{2U}(R')^{-1}(1-\beta)]. \\\\
\end{array} \right\} \label{2}
\end{equation} 
\begin{enumerate}
\item[$(i).$] From the Eqs. $(5)~\mbox{and}~ (6)$, we get\\

\noindent$\big(\tilde{A}_1\oplus\tilde{A}_2\big)_{\alpha}^{L}=A_{1\alpha}^{L}+A_{2\alpha}^{L}$\\
$~~~~~~~~~~~~~~~~~~~~~=[a_1+a_2-(l^{\mu}_{1L}+l^{\mu}_{2L})L^{-1}(\alpha),a_1+a_2+(r^{\mu}_{1L}+r^{\mu}_{2L})R^{-1}(\alpha)].$\\

\noindent $\big(\tilde{A}_1\oplus\tilde{A}_2\big)_{\alpha}^{U}=A_{1\alpha}^{U}+A_{2\alpha}^{U}$\\
$~~~~~~~~~~~~~~~~~~~~~~= [a_1+a_2-(l'^{\mu}_{1U}+l'^{\mu}_{2U})(L')^{-1}(\alpha),a_1+a_2+(r'^{\mu}_{1U}+r'^{\mu}_{2U})(R')^{-1}(\alpha)].$\\

\noindent $\big(\tilde{A}_1\oplus\tilde{A}_2\big)_{\beta}^{L}=A_{1\beta}^{L}+A_{2\beta}^{L}$\\
$~~~~~~~~~~~~~~~~~~~~~=[a_1+a_2-(l^{\nu}_{1L}+l^{\nu}_{2L})L^{-1}(1-\beta),a_1+a_2+(r^{\nu}_{1L}+r^{\nu}_{2L})R^{-1}(1-\beta)].$\\

\noindent $\big(\tilde{A}_1\oplus\tilde{A}_2\big)_{\beta}^{U}=A_{1\beta}^{U}+A_{2\beta}^{U}$\\
$~~~~~~~~~~~~~~~~~~~~~~=[a_1+a_2-(l'^{\nu}_{1U}+l'^{\nu}_{2U})(L')^{-1}(1-\beta),a_1+a_2+(r'^{\nu}_{1U}+r'^{\nu}_{2U})(R')^{-1}(1-\beta)].$\\

\noindent Since $L, R, L'~\mbox{and}~R'$ are decreasing functions on $[0, \infty)$ with $L(0)=R(0)=L'(0)=R'(0)=1$, there exists $\alpha_{0} \in (0,1]$, such that $L^{-1}(\alpha_0)=R^{-1}(\alpha_0)= (L')^{-1}(\alpha_0)=(R')^{-1}(\alpha_0)=1$. Hence,
\begin{equation} \hspace{2cm} \big(\tilde{A}_1\oplus\tilde{A}_2\big)_{\alpha_0}^{L}=[a_1+a_2-(l^{\mu}_{1L}+l^{\mu}_{2L}),a_1+a_2+(r^{\mu}_{1L}+r^{\mu}_{2L})]. \end{equation}
\vspace{-0.5cm}
\begin{equation} \hspace{2cm}  \big(\tilde{A}_1\oplus\tilde{A}_2\big)_{\alpha_0}^{U}=[a_1+a_2-(l'^{\mu}_{1U}+l'^{\mu}_{2U}), a_1+ a_2+(r'^{\mu}_{1U}+r'^{\mu}_{2U})]. \end{equation}
\medspace
\noindent Also, choosing $\beta_{0}=1-\alpha_{0} \in (0,1]$, we get
\begin{equation}\hspace{2cm} \big(\tilde{A}_1\oplus\tilde{A}_2\big)_{\beta_0}^{L}=[a_1+a_2-(l^{\nu}_{1L}+l^{\nu}_{2L}), a_1+ a_2+(r^{\nu}_{1L}+r^{\nu}_{2L})]. \end{equation}
\vspace{-0.5cm}
\begin{equation}\hspace{2cm}\big(\tilde{A}_1\oplus\tilde{A}_2\big)_{\beta_0}^{U}=[a_1+a_2-(l'^{\nu}_{1U}+l'^{\nu}_{2U}), a_1+ a_2+(r'^{\nu}_{1U}+r'^{\nu}_{2U})]. \end{equation}

\noindent Further,
\begin{equation}\hspace{0.5cm} \big(\tilde{A}_1\oplus\tilde{A}_2\big)_{\alpha=1}^{L}=\big(\tilde{A}_1\oplus\tilde{A}_2\big)_{\alpha=1}^{U}=\big(\tilde{A}_1\oplus\tilde{A}_2\big)_{\beta=0}^{L}
=\big(\tilde{A}_1\oplus\tilde{A}_2\big)_{\beta=0}^{U}=[a_1+a_2, a_1+a_2].\end{equation}
 
\noindent Now, since $\tilde{A}_1~\mbox{and}~\tilde{A}_2$ are $LR$-type IVIFNs, therefore\\\\
\medskip
\noindent$l'^{\mu}_{1U} \geq l^{\mu}_{1L} >0,~~r'^{\mu}_{1U} \geq r^{\mu}_{1L} >0,~~l^{\nu}_{1L} \geq l'^{\nu}_{1U} >0,~~r^{\nu}_{1L} \geq r'^{\nu}_{1U} >0,~~l^{\nu}_{1L} \geq l^{\mu}_{1L} >0,~~r^{\nu}_{1L} \geq r^{\mu}_{1L} >0,$\\
\medskip
$l'^{\nu}_{1U} \geq l'^{\mu}_{1U} >0,~~r'^{\nu}_{1U} \geq r'^{\mu}_{1U} >0,~~$ and\\
\medskip
$l'^{\mu}_{2U} \geq l^{\mu}_{2L} >0,~~r'^{\mu}_{2U} \geq r^{\mu}_{2L} >0,~~l^{\nu}_{2L} \geq l'^{\nu}_{2U} >0,~~r^{\nu}_{2L} \geq r'^{\nu}_{2U} >0,~~l^{\nu}_{2L} \geq l^{\mu}_{2L} >0,~~r^{\nu}_{2L} \geq r^{\mu}_{2L} >0,$\\
\medskip
$l'^{\nu}_{2U} \geq l'^{\mu}_{2U} >0,~~r'^{\nu}_{2U} \geq r'^{\mu}_{2U} >0.$\\\\
\noindent Thus,
\begin{equation}
 \renewcommand{\arraystretch}{.60}
  \left.\begin{array}{r@{\;}l}
 l'^{\mu}_{1U}+l'^{\mu}_{2U} \geq l^{\mu}_{1L}+l^{\mu}_{2L} >0,~
r'^{\mu}_{1U}+r'^{\mu}_{2U} \geq r^{\mu}_{1L}+r^{\mu}_{2L} >0,~ l^{\nu}_{1L}+l^{\nu}_{2L} \geq l'^{\nu}_{1U}+l'^{\nu}_{2U} >0,\\\\

r^{\nu}_{1L}+r^{\nu}_{2L} \geq r'^{\nu}_{1U}+r'^{\nu}_{2U} >0,~ l^{\nu}_{1L}+l^{\nu}_{2L} \geq l^{\mu}_{1L}+l^{\mu}_{2L} >0,~r^{\nu}_{1L}+r^{\nu}_{2L} \geq r^{\mu}_{1L}+r^{\mu}_{2L} >0,\\\\

\noindent l'^{\nu}_{1U}+l'^{\nu}_{2U} \geq l'^{\mu}_{1U}+l'^{\mu}_{2U} >0,~r'^{\nu}_{1U}+r'^{\nu}_{2U} \geq r'^{\mu}_{1U}+r'^{\mu}_{2U} >0.\\
\end{array} \right\} 
\end{equation} 
\noindent Combining Eqs. $(7)$ \textendash $~(11)$, we have\\\\
\noindent $\tilde{A}_1\oplus\tilde{A}_2=(a_1+a_2; l_{1L}^{\mu}+l_{2L}^{\mu},r_{1L}^{\mu}+r_{2L}^{\mu},l'^{\mu}_{1U}+l'^{\mu}_{2U},r'^{\mu}_{1U}+r'^{\mu}_{2U};l^{\nu}_{1L}+l^{\nu}_{2L},r^{\nu}_{1L}+r^{\nu}_{2L},l'^{\nu}_{1U}+l'^{\nu}_{2U},r'^{\nu}_{1U}+r'^{\nu}_{2U})_{LR},$\\

\noindent where the conditions for $LR$-type form of $\tilde{A}_1\oplus\tilde{A}_2$ holds from $(12)$.\\
\noindent Hence, $(i)$ is proved.\\
\item[$(ii).$] Using equations $(5)$ and $(6)$, we can write\\\\
\noindent$\big(\tilde{A}_1\ominus\tilde{A}_2\big)_{\alpha}^{L}=A_{1\alpha}^{L}-A_{2\alpha}^{L}=[a_1-a_2-l^{\mu}_{1L}L^{-1}(\alpha)-r^{\mu}_{2L}R^{-1}(\alpha),~a_1-a_2+r^{\mu}_{1L}R^{-1}(\alpha)+l^{\mu}_{2L}L^{-1}(\alpha)].$\\

\noindent $\big(\tilde{A}_1\ominus\tilde{A}_2\big)_{\alpha}^{U}=A_{1\alpha}^{U}-A_{2\alpha}^{U}= [a_1-a_2-l'^{\mu}_{1U}(L')^{-1}(\alpha)-r'^{\mu}_{2U}(R')^{-1}(\alpha),~a_1-a_2+r'^{\mu}_{1U}(R')^{-1}(\alpha)+l'^{\mu}_{2U}(L')^{-1}(\alpha)].$\\

\noindent $\big(\tilde{A}_1\ominus \tilde{A}_2\big)_{\beta}^{L}=A_{1\beta}^{L}-A_{2\beta}^{L}=[a_1-a_2-l^{\nu}_{1L}L^{-1}(1-\beta)-r^{\nu}_{2L}R^{-1}(1-\beta), a_1-a_2+r^{\nu}_{1L}R^{-1}(1-\beta)+l^{\nu}_{2L}L^{-1}(1-\beta)].$\\

\noindent $\big(\tilde{A}_1\ominus \tilde{A}_2\big)_{\beta}^{U}=A_{1\beta}^{U}-A_{2\beta}^{U}=[a_1-a_2-l'^{\nu}_{1U}(L')^{-1}(1-\beta)-r'^{\nu}_{2U}(R')^{-1}(1-\beta),a_1-a_2+r'^{\nu}_{1U}(R')^{-1}(1-\beta)+$\\

$~~~~~~~~~~~~~~~~~~~~~~~~~~~~~~~~~~~~~~~~~~~~~~~~~l'^{\nu}_{2U}(L')^{-1}(1-\beta)].$\\

\noindent Now, taking $\alpha=\alpha_{0}~\mbox{and}~\beta=\beta_{0} \in (0,1],$ we get
\begin{equation}
 \renewcommand{\arraystretch}{.60}
 \hspace{2cm} \left.\begin{array}{r@{\;}l}
\big(\tilde{A}_1\ominus\tilde{A}_2\big)_{\alpha_0}^{L}=[a_1-a_2-l^{\mu}_{1L}-r^{\mu}_{2L},a_1-a_2+r^{\mu}_{1L}+l^{\mu}_{2L}],\\\\
\big(\tilde{A}_1\ominus\tilde{A}_2\big)_{\alpha_0}^{U}=[a_1-a_2-l'^{\mu}_{1U}-r'^{\mu}_{2U},a_1-a_2+r'^{\mu}_{1U}+l'^{\mu}_{2U}], \\\\
\big(\tilde{A}_1\ominus \tilde{A}_2\big)_{\beta_0}^{L}=[a_1-a_2-l^{\nu}_{1L}-r^{\nu}_{2L},a_1-a_2+r^{\nu}_{1L}+l^{\nu}_{2L}], \\\\
\big(\tilde{A}_1\ominus \tilde{A}_2\big)_{\beta_0}^{U}=[a_1-a_2-l'^{\nu}_{1U}-r'^{\nu}_{2U},a_1-a_2+r'^{\nu}_{1U}+l'^{\nu}_{2U}].\\
\end{array} \right\} \label{13}
\end{equation} 

\noindent Further, on substituting $\alpha=1~\mbox{and}~ \beta=0$, we obtain
\begin{equation}\hspace{0.5cm}\big(\tilde{A}_1\ominus\tilde{A}_2\big)_{\alpha=1}^{L}=\big(\tilde{A}_1\ominus\tilde{A}_2\big)_{\alpha=1}^{U}=\big(\tilde{A}_1\ominus\tilde{A}_2\big)_{\beta=0}^{L}
=\big(\tilde{A}_1\ominus\tilde{A}_2\big)_{\beta=0}^{U}=[a_1-a_2, a_1-a_2].\end{equation}
\noindent Also, using the fact $\tilde{A}_1~\mbox{and}~\tilde{A}_2$ are $LR$-type IVIFNs, we have\\\\
$l'^{\mu}_{1U}+r'^{\mu}_{2U} \geq l^{\mu}_{1L}+r^{\mu}_{2L} >0,~~r'^{\mu}_{1U}+l'^{\mu}_{2U} \geq r^{\mu}_{1L}+l^{\mu}_{2L} >0,~~l^{\nu}_{1L}+r^{\nu}_{2L} \geq l'^{\nu}_{1U}+r'^{\nu}_{2U} >0,~~r^{\nu}_{1L}+l^{\nu}_{2L} \geq r'^{\nu}_{1U}+l'^{\nu}_{2U} >0,$\\\\
$l^{\nu}_{1L}+r^{\nu}_{2L} \geq l^{\mu}_{1L}+r^{\mu}_{2L} >0,~~r^{\nu}_{1L}+l^{\nu}_{2L} \geq r^{\mu}_{1L}+l^{\mu}_{2L} >0,~~l'^{\nu}_{1U}+r'^{\nu}_{2U} \geq l'^{\mu}_{1U}+r'^{\mu}_{2U} >0,~~r'^{\nu}_{1U}+l'^{\nu}_{2U} \geq r'^{\mu}_{1U}+l'^{\mu}_{2U} >0.$\\\\  
\noindent Finally, from the Eqs. $(13)~\mbox{and}~(14)$, we have\\\\
\noindent $\tilde{A}_1\ominus\tilde{A}_2=(a_1-a_2; l_{1L}^{\mu}+r_{2L}^{\mu},r_{1L}^{\mu}+l_{2L}^{\mu},l'^{\mu}_{1U}+r'^{\mu}_{2U},r'^{\mu}_{1U}+l'^{\mu}_{2U};l^{\nu}_{1L}+r^{\nu}_{2L},r^{\nu}_{1L}+l^{\nu}_{2L},l'^{\nu}_{1U}+r'^{\nu}_{2U},r'^{\nu}_{1U}+l'^{\nu}_{2U})_{LR}$,\\

\noindent along-with $\tilde{A}_1\ominus\tilde{A}_2$ retains the form of a $LR$-type IVIFN. This proves $(ii)$.\\
\end{enumerate}

\noindent{\bf{Proposition 3.1.2.}} \textit{Let $\tilde{A}=(a; l_{L}^{\mu},r_{L}^{\mu},l'^{\mu}_{U},r'^{\mu}_{U};l^{\nu}_{L},r^{\nu}_{L},l'^{\nu}_{U},r'^{\nu}_{U})_{LR}$ be an $LR$-type IVIFN and $\lambda$ be any real number. Then}\\
$$\lambda \tilde{A}=\begin{cases}
(\lambda a;\lambda l_{L}^{\mu},\lambda r_{L}^{\mu},\lambda l'^{\mu}_{U},\lambda r'^{\mu}_{U};\lambda l^{\nu}_{L},\lambda r^{\nu}_{L},\lambda l'^{\nu}_{U},\lambda r'^{\nu}_{U})_{LR} & \mbox{if}   ~~\lambda \geq 0,\\\\
(\lambda a;-\lambda r_{L}^{\mu},-\lambda l_{L}^{\mu},-\lambda r'^{\mu}_{U},-\lambda l'^{\mu}_{U};-\lambda r^{\nu}_{L},-\lambda l^{\nu}_{L},-\lambda r'^{\nu}_{U},-\lambda l'^{\nu}_{U})_{LR} & \mbox{if}   ~~\lambda < 0.\\
\end{cases}$$

\noindent{\bf{Proof.}} From Theorem $3.1$, the $\alpha$ and $\beta$-cuts of $\tilde{A}$ are given by
\begin{equation}
 \renewcommand{\arraystretch}{.60}
\hspace{3cm} \left.\begin{array}{r@{\;}l}
A_{\alpha}^{L}= [a-l^{\mu}_{L}L^{-1}(\alpha), a+r^{\mu}_{L}R^{-1}(\alpha)],\\\\
A_{\alpha}^{U}= [a-l'^{\mu}_{U}(L')^{-1}(\alpha), a+r'^{\mu}_{U}(R')^{-1}(\alpha)],\\\\
 A_{\beta}^{L}= [a-l^{\nu}_{L}L^{-1}(1-\beta), a+r^{\nu}_{L}R^{-1}(1-\beta)],\\\\
A_{\beta}^{U}= [a-l'^{\nu}_{U}(L')^{-1}(1-\beta), a+r'^{\nu}_{U}(R')^{-1}(1-\beta)].
\end{array} \right\} \label{2}
\end{equation} 

\noindent Using the expression $(15)$, we have\\\\
\noindent$~~~~~~~~~~~~~~~~~~~~~~~~~~~~~~~~~~~~~~~(\lambda \tilde{A})_{\alpha}^{L}=\lambda A_{\alpha}^{L}= [\lambda,\lambda]A_{\alpha}^{L}=[\lambda,\lambda][a-l^{\mu}_{L}L^{-1}(\alpha), a+r^{\mu}_{L}R^{-1}(\alpha)],$\\\\
\noindent$~~~~~~~~~~~~~~~~~~~~~~~~~~~~~~~~~~~~~~~(\lambda \tilde{A})_{\alpha}^{U}=\lambda A_{\alpha}^{U}= [\lambda,\lambda]A_{\alpha}^{U}=[\lambda,\lambda][a-l'^{\mu}_{U}(L')^{-1}(\alpha), a+r'^{\mu}_{U}(R')^{-1}(\alpha)],$\\\\
\noindent$~~~~~~~~~~~~~~~~~~~~~~~~~~~~~~~~~~~~~~~(\lambda \tilde{A})_{\beta}^{L}=\lambda A_{\beta}^{L}= [\lambda,\lambda]A_{\beta}^{L}=[\lambda,\lambda][a-l^{\nu}_{L}L^{-1}(1-\beta), a+r^{\nu}_{L}R^{-1}(1-\beta)],$\\\\
\noindent$~~~~~~~~~~~~~~~~~~~~~~~~~~~~~~~~~~~~~~~(\lambda \tilde{A})_{\beta}^{U}=\lambda A_{\beta}^{U}= [\lambda,\lambda]A_{\beta}^{U}=[\lambda,\lambda][a-l'^{\nu}_{U}(L')^{-1}(1-\beta), a+r'^{\nu}_{U}(R')^{-1}(1-\beta)]$.\\

\noindent \textbf{Case 1.} $\tilde{A}$ is a non-negative $LR$-type IVIFN.\\

\noindent \underline{Sub-case 1.} $\lambda \geq 0$. \\

\noindent Since, $\tilde{A}$ is a non-negative $LR$-type IVIFN, i.e., $a-l^{\nu}_L \geq 0$. Thus, $a-l^{\mu}_L \geq 0,~a-l'^{\mu}_{U} \geq 0,~ a-l^{\nu}_{L} \geq 0,~a-l'^{\nu}_{U} \geq 0.$\\
$ \mbox{This further implies}~a-l^{\mu}_L L^{-1}(\alpha) \geq 0,~a-l'^{\mu}_{U}(L')^{-1}(\alpha) \geq 0,~a-l^{\nu}_{L}L^{-1}(1-\beta) \geq 0,~a-l'^{\nu}_{U}(L')^{-1}(1-\beta) \geq 0,$\\
$\forall~\alpha,\beta \in [0,1]$. Hence, we get\\\\
\noindent $(\lambda \tilde{A})_{\alpha}^{L}=[\lambda,\lambda][a-l^{\mu}_{L}L^{-1}(\alpha), a+r^{\mu}_{L}R^{-1}(\alpha)]=[\lambda(a-l^{\mu}_{L}L^{-1}(\alpha)),\lambda (a+r^{\mu}_{L}R^{-1}(\alpha))],$\\\\
\noindent $(\lambda \tilde{A})_{\alpha}^{U}=[\lambda(a-l'^{\mu}_{U}(L')^{-1}(\alpha)),\lambda (a+r'^{\mu}_{U}(R')^{-1}(\alpha))],$\\\\
\noindent $(\lambda \tilde{A})_{\beta}^{L}=[\lambda(a-l^{\nu}_{L}L^{-1}(1-\beta)),\lambda (a+r^{\nu}_{L}R^{-1}(1-\beta))],$\\\\
\noindent $(\lambda \tilde{A})_{\beta}^{U}=[\lambda (a-l'^{\nu}_{U}(L')^{-1}(1-\beta)),\lambda (a+r'^{\nu}_{U}(R')^{-1}(1-\beta))]$.\\\\
\noindent Further, as $L, R, L'~\mbox{and}~R'$ are decreasing functions on $[0, \infty)$ with $L(0)=R(0)=L'(0)=R'(0)=1$, there exists $\alpha_{0} \in (0,1]$, such that $L^{-1}(\alpha_0)=R^{-1}(\alpha_0)= (L')^{-1}(\alpha_0)=(R')^{-1}(\alpha_0)=1$. Therefore,
\begin{equation}
 \renewcommand{\arraystretch}{.60}
\hspace{2cm} \left.\begin{array}{r@{\;}l}
(\lambda \tilde{A})_{\alpha_0}^{L}=[\lambda(a-l^{\mu}_{L}),\lambda (a+r^{\mu}_{L})]=[\lambda a-\lambda l^{\mu}_{L},\lambda a+\lambda r^{\mu}_{L}],\\\\
(\lambda \tilde{A})_{\alpha_0}^{U}=[\lambda(a-l'^{\mu}_{U}),\lambda (a+r'^{\mu}_{U})]=[\lambda a-\lambda l'^{\mu}_{U},\lambda a+\lambda r'^{\mu}_{U}].
\end{array} \right\} \label{2}
\end{equation} 

\noindent Choosing $\beta_{0}=1-\alpha_0 \in (0,1]$, we have
\begin{equation}
 \renewcommand{\arraystretch}{.60}
\hspace{2cm}  \left.\begin{array}{r@{\;}l}
(\lambda \tilde{A})_{\beta_0}^{L}=[\lambda(a-l^{\nu}_{L}),\lambda (a+r^{\nu}_{L})]=[\lambda a-\lambda l^{\nu}_{L},\lambda a+\lambda r^{\nu}_{L}],\\\\
 (\lambda \tilde{A})_{\beta_0}^{U}=[\lambda (a-l'^{\nu}_{U}),\lambda (a+r'^{\nu}_{U})] =[\lambda a-\lambda l'^{\nu}_{U},\lambda a+\lambda r'^{\nu}_{U}].
\end{array}\right\} \label{2}
\end{equation} 
\noindent Putting $\alpha=1~\mbox{and}~ \beta=0$, we get
\begin{equation} \hspace{3cm}(\lambda \tilde{A})_{\alpha=1}^{L}=(\lambda \tilde{A})_{\alpha=1}^{U}=(\lambda \tilde{A})_{\beta=0}^{L}
=(\lambda \tilde{A})_{\beta=0}^{U}=[\lambda a,\lambda a]. \end{equation} 
\noindent Since $\tilde{A}$ is an $LR$-type IVIFN and $\lambda \geq 0$, we obtain\\

\noindent$\lambda l'^{\mu}_{U} \geq \lambda l^{\mu}_{L} >0,~~\lambda r'^{\mu}_{U} \geq \lambda r^{\mu}_{L} >0,~~\lambda l^{\nu}_{L} \geq \lambda l'^{\nu}_{U} >0,~~\lambda r^{\nu}_{L} \geq \lambda r'^{\nu}_{U} >0,~~\lambda l^{\nu}_{L} \geq \lambda l^{\mu}_{L} >0,~~\lambda r^{\nu}_{L} \geq \lambda r^{\mu}_{L} >0,$\\\\
$\lambda l'^{\nu}_{U} \geq \lambda l'^{\mu}_{U} >0,~~\lambda r'^{\nu}_{U} \geq \lambda r'^{\mu}_{U} >0.$ \\

\noindent Hence, combining $(16)$ \textendash $~(18)$, we have
$$\lambda \tilde{A}=(\lambda a;\lambda l_{L}^{\mu},\lambda r_{L}^{\mu},\lambda l'^{\mu}_{U},\lambda r'^{\mu}_{U};\lambda l^{\nu}_{L},\lambda r^{\nu}_{L},\lambda l'^{\nu}_{U},\lambda r'^{\nu}_{U})_{LR}.$$

\noindent \underline{Sub-case 2.} $~~\lambda < 0$. \\

\noindent Since, $\tilde{A}$ is a non-negative $LR$-type IVIFN, i.e., $a-l^{\nu}_L \geq 0$. Thus, $a-l^{\mu}_L \geq 0,~a-l'^{\mu}_{U} \geq 0,~ a-l^{\nu}_{L} \geq 0,~a-l'^{\nu}_{U} \geq 0 $\\
$\implies~a-l^{\mu}_L L^{-1}(\alpha) \geq 0,~a-l'^{\mu}_{U}(L')^{-1}(\alpha) \geq 0,~~ a-l^{\nu}_{L}L^{-1}(1-\beta) \geq 0,~a-l'^{\nu}_{U}(L')^{-1}(1-\beta) \geq 0,~~\forall~\alpha,\beta \in [0,1]$.\\
It follows that\\
\noindent $~~~~~~~~~~~~~~~~~~~(\lambda \tilde{A})_{\alpha}^{L}=[\lambda,\lambda][a-l^{\mu}_{L}L^{-1}(\alpha), a+r^{\mu}_{L}R^{-1}(\alpha)]=[\lambda(a+r^{\mu}_{L}R^{-1}(\alpha)),\lambda (a-l^{\mu}_{L}L^{-1}(\alpha))],$\\

\noindent $~~~~~~~~~~~~~~~~~~(\lambda \tilde{A})_{\alpha}^{U}=[\lambda(a+r'^{\mu}_{U}(R')^{-1}(\alpha)),\lambda (a-l'^{\mu}_{U}(L')^{-1}(\alpha))],$\\

\noindent $~~~~~~~~~~~~~~~~~~(\lambda \tilde{A})_{\beta}^{L}=[\lambda(a+r^{\nu}_{L}R^{-1}(1-\beta)),\lambda (a-l^{\nu}_{L}L^{-1}(1-\beta))],$\\

\noindent $~~~~~~~~~~~~~~~~~~(\lambda \tilde{A})_{\beta}^{U}=[\lambda(a+r'^{\nu}_{U}(R')^{-1}(1-\beta)) ,\lambda (a-l'^{\nu}_{U}(L')^{-1}(1-\beta))]$.\\

\noindent Taking $\alpha=\alpha_0$, $\beta=\beta_0=1-\alpha_0$ and using the fact that $L^{-1}(\alpha_0)=R^{-1}(\alpha_0)= (L')^{-1}(\alpha_0)=(R')^{-1}(\alpha_0)=1$, we have
\begin{equation}
 \renewcommand{\arraystretch}{.60}
 \hspace{2cm} \left.\begin{array}{r@{\;}l}
(\lambda \tilde{A})_{\alpha_0}^{L}=[\lambda(a+r^{\mu}_{L}),\lambda(a-l^{\mu}_{L}) ]=[\lambda a+\lambda r^{\mu}_{L},\lambda a-\lambda l^{\mu}_{L}],\\\\\

(\lambda \tilde{A})_{\alpha_0}^{U}=[\lambda(a+r'^{\mu}_{U}),\lambda (a-l'^{\mu}_{U})]=[\lambda a+\lambda r'^{\mu}_{U} ,\lambda a-\lambda l'^{\mu}_{U}],\\\\

(\lambda \tilde{A})_{\beta_0}^{L}=[\lambda (a+r^{\nu}_{L}),\lambda (a-l^{\nu}_{L})]=[\lambda a+\lambda r^{\nu}_{L},\lambda a-\lambda l^{\nu}_{L}],\\\\

(\lambda \tilde{A})_{\beta_0}^{U}=[\lambda (a+r'^{\nu}_{U}),\lambda (a-l'^{\nu}_{U})]=[\lambda a+\lambda r'^{\nu}_{U},\lambda -\lambda l'^{\nu}_{U}].\\
\end{array} \right\} \label{2}
\end{equation} 
\noindent Putting $\alpha=1~\mbox{and}~ \beta=0$, we obtain
\begin{equation} \hspace{3cm} (\lambda \tilde{A})_{\alpha=1}^{L}=(\lambda \tilde{A})_{\alpha=1}^{U}=(\lambda \tilde{A})_{\beta=0}^{L}
=(\lambda \tilde{A})_{\beta=0}^{U}=[\lambda a,\lambda a]. \end{equation}

\noindent Further, from the fact that $\tilde{A}$ is an $LR$-type IVIFN and $\lambda < 0$, we have\\\\
\noindent$-\lambda l'^{\mu}_{U} \geq -\lambda l^{\mu}_{L} >0,~~-\lambda r'^{\mu}_{U} \geq -\lambda r^{\mu}_{L} >0,~~-\lambda l^{\nu}_{L} \geq -\lambda l'^{\nu}_{U} >0,~~-\lambda r^{\nu}_{L} \geq -\lambda r'^{\nu}_{U} >0,~~-\lambda l^{\nu}_{L} \geq -\lambda l^{\mu}_{L} >0,$\\\\
\noindent$-\lambda r^{\nu}_{L} \geq -\lambda r^{\mu}_{L} >0,~~-\lambda l'^{\nu}_{U} \geq -\lambda l'^{\mu}_{U} >0,~~-\lambda r'^{\nu}_{U} \geq -\lambda r'^{\mu}_{U} >0.$ \\

\noindent Hence, the expressions $(19)~\mbox{and}~(20)$ finally yield
$$\lambda \tilde{A}=(\lambda a;-\lambda r_{L}^{\mu},-\lambda l_{L}^{\mu},-\lambda r'^{\mu}_{U},-\lambda l'^{\mu}_{U};-\lambda r^{\nu}_{L},-\lambda l^{\nu}_{L},-\lambda r'^{\nu}_{U},-\lambda l'^{\nu}_{U})_{LR}.$$

\noindent Therefore, the result is proved for a non-negative $LR$-type IVIFN.\\

\noindent \textbf{Case 2.} $\tilde{A}$ is an $LR$-type IVIFN such that $a-l^{\nu}_{L} <0$ and $a-l'^{\nu}_{U} \geq 0$.

\noindent Since, $~~~~~~a-l^{\mu}_L \geq 0,~a-l'^{\mu}_{U} \geq 0,~ a-l'^{\nu}_{U} \geq 0$\\
 $~~~~~~\implies~a-l^{\mu}_L L^{-1}(\alpha) \geq 0,~a-l'^{\mu}_{U}(L')^{-1}(\alpha) \geq 0,~~ a-l'^{\nu}_{U}(L')^{-1}(1-\beta) \geq 0,~\forall~\alpha,\beta \in [0,1]$.\\

\noindent So, expressions for $(\lambda \tilde{A})_{\alpha}^{L}$, $(\lambda \tilde{A})_{\alpha}^{U}$ and $(\lambda \tilde{A})_{\beta}^{U}$ are same as derived in Case 1, and
$$(\lambda \tilde{A})_{\beta}^{L}=\lambda A_{\beta}^{L}= [\lambda,\lambda]A_{\beta}^{L}=[\lambda,\lambda] [a-l^{\nu}_{L}L^{-1}(1-\beta), a+r^{\nu}_{L}R^{-1}(1-\beta)].$$

\noindent Now, if $a-l^{\nu}_{L} <0$, then either
$$ a-l^{\nu}_{L}L^{-1}(1-\beta) <0 \iff \displaystyle\frac{a}{l^{\nu}_{L}}<L^{-1}(1-\beta) \iff \beta >1-L\Bigg(\displaystyle\frac{a}{l^{\nu}_{L}}\Bigg)$$

 \noindent or $$a-l^{\nu}_{L}L^{-1}(1-\beta) \geq 0 \iff \displaystyle\frac{a}{l^{\nu}_{L}} \geq L^{-1}(1-\beta) \iff \beta \leq 1-L\Bigg(\displaystyle\frac{a}{l^{\nu}_{L}}\Bigg).$$

\noindent$\mbox{Hence,}~~ a-l^{\nu}_{L}L^{-1}(1-\beta) <0 ~\mbox{for}~\beta >1-L\Bigg(\displaystyle\frac{a}{l^{\nu}_{L}}\Bigg)~~\mbox{and}~~a-l^{\nu}_{L}L^{-1}(1-\beta) \geq 0 ~\mbox{for}~\beta \leq 1-L\Bigg(\displaystyle\frac{a}{l^{\nu}_{L}}\Bigg)$.\\

\noindent \textbf{(a).} Let $1-L\Bigg(\displaystyle\frac{a}{l^{\nu}_{L}}\Bigg) < \beta \leq 1~$ or $~a-l^{\nu}_{L}L^{-1}(1-\beta) <0 .$\\

\noindent \underline{Sub-case 1.} $~~\lambda \geq 0$.\\

\noindent$(\lambda \tilde{A})_{\beta}^{L}=[\lambda,\lambda][a-l^{\nu}_{L}L^{-1}(1-\beta), a+r^{\nu}_{L}R^{-1}(1-\beta)]=[\lambda(a-l^{\nu}_{L}L^{-1}(1-\beta)),\lambda (a+r^{\nu}_{L}R^{-1}(1-\beta))]$.\\

\noindent On the similar lines as in Case 1, we get\\

\noindent$(\lambda \tilde{A})_{\beta_0}^{L}=[\lambda (a-l^{\nu}_{L}),\lambda (a+r^{\nu}_{L})]=[\lambda a-\lambda l^{\nu}_{L},\lambda a+\lambda r^{\nu}_{L}]$,\\

\noindent This yields
$$\lambda \tilde{A}=(\lambda a;\lambda l_{L}^{\mu},\lambda r_{L}^{\mu},\lambda l'^{\mu}_{U},\lambda r'^{\mu}_{U};\lambda l^{\nu}_{L},\lambda r^{\nu}_{L},\lambda l'^{\nu}_{U},\lambda r'^{\nu}_{U})_{LR}.$$

\noindent \underline{Sub-case 2.} $~~\lambda <0$.
$$(\lambda \tilde{A})_{\beta}^{L}=[\lambda(a+r^{\nu}_{L}R^{-1}(1-\beta)),\lambda (a-l^{\nu}_{L}L^{-1}(1-\beta))].$$
\noindent Further, we have
$$(\lambda \tilde{A})_{\beta_0}^{L}=[\lambda (a+r^{\nu}_{L}),\lambda (a-l^{\nu}_{L})]=[\lambda a+\lambda r^{\nu}_{L},\lambda a-\lambda l^{\nu}_{L}]$$
\noindent which gives
$$\lambda \tilde{A}=(\lambda a;-\lambda r_{L}^{\mu},-\lambda l_{L}^{\mu},-\lambda r'^{\mu}_{U},-\lambda l'^{\mu}_{U};-\lambda r^{\nu}_{L},-\lambda l^{\nu}_{L},-\lambda r'^{\nu}_{U},-\lambda l'^{\nu}_{U})_{LR}.$$

\noindent \textbf{(b).} Let $0 \leq \beta \leq 1-L\Bigg(\displaystyle\frac{a}{l^{\nu}_{L}}\Bigg)~$ or $~a-l^{\nu}_{L}L^{-1}(1-\beta) \geq 0. $\\

\noindent \underline{Sub-case 1.} $~~\lambda \geq 0$.\\

\noindent$(\lambda \tilde{A})_{\beta}^{L}=[\lambda,\lambda][a-l^{\nu}_{L}L^{-1}(1-\beta), a+r^{\nu}_{L}R^{-1}(1-\beta)]=[\lambda(a-l^{\nu}_{L}L^{-1}(1-\beta)),\lambda (a+r^{\nu}_{L}R^{-1}(1-\beta))]$.\\

\noindent Following the steps of Case 1, we obtain
$$(\lambda \tilde{A})_{\beta_0}^{L}=[\lambda (a-l^{\nu}_{L}),\lambda (a+r^{\nu}_{L})]=[\lambda a-\lambda l^{\nu}_{L},\lambda a+\lambda r^{\nu}_{L}]$$
\noindent which yields
$$\lambda \tilde{A}=(\lambda a;\lambda l_{L}^{\mu},\lambda r_{L}^{\mu},\lambda l'^{\mu}_{U},\lambda r'^{\mu}_{U};\lambda l^{\nu}_{L},\lambda r^{\nu}_{L},\lambda l'^{\nu}_{U},\lambda r'^{\nu}_{U})_{LR}.$$

\noindent \underline{Sub-case 2.} $~~\lambda <0$.\\
$$(\lambda \tilde{A})_{\beta}^{L}=[\lambda(a+r^{\nu}_{L}R^{-1}(1-\beta)),\lambda (a-l^{\nu}_{L}L^{-1}(1-\beta))].$$

\noindent Proceeding on the lines of Case 1 and taking $\beta=\beta_0,$ we get
$$(\lambda \tilde{A})_{\beta_0}^{L}=[\lambda (a+r^{\nu}_{L}),\lambda (a-l^{\nu}_{L})]=[\lambda a+\lambda r^{\nu}_{L},\lambda a-\lambda l^{\nu}_{L}].$$
\noindent Therefore, we have
$$\lambda \tilde{A}=(\lambda a;-\lambda r_{L}^{\mu},-\lambda l_{L}^{\mu},-\lambda r'^{\mu}_{U},-\lambda l'^{\mu}_{U};-\lambda r^{\nu}_{L},-\lambda l^{\nu}_{L},-\lambda r'^{\nu}_{U},-\lambda l'^{\nu}_{U})_{LR}.$$

\noindent Thus, the result follows. Rest all the Cases can also be proved on the similar lines. This completes the proof.\\

\noindent \textbf{Corollary 3.1.1} Proposition 3.1.2 can also be restated as: If $\tilde{A}$ be an $LR$-type IVIFN and $\lambda$ be any arbitrary real number, then\\

\noindent$\lambda \tilde{A}=\big(\lambda a;\max \{ \lambda l_{L}^{\mu},-\lambda r_{L}^{\mu}\},\max \{\lambda r_{L}^{\mu},-\lambda l_{L}^{\mu}\},\max \{\lambda l'^{\mu}_{U},-\lambda r'^{\mu}_{U}\},\max \{\lambda r'^{\mu}_{U},-\lambda l'^{\mu}_{U}\};\max \{\lambda l^{\nu}_{L},-\lambda r^{\nu}_{L}\},$\\\\
$~~~~~~~~~~~~\max \{\lambda r^{\nu}_{L},-\lambda l^{\nu}_{L}\},\max \{\lambda l'^{\nu}_{U},-\lambda r'^{\nu}_{U}\},\max \{\lambda r'^{\nu}_{U},-\lambda l'^{\nu}_{U}\}\big)_{LR}$. \\

\noindent{\bf{Proposition 3.1.3.}} \textit{Let $\tilde{A}_1=(a_1; l_{1L}^{\mu},r_{1L}^{\mu},l'^{\mu}_{1U},r'^{\mu}_{1U};l^{\nu}_{1L},r^{\nu}_{1L},l'^{\nu}_{1U},r'^{\nu}_{1U})_{LR}$ be an $LR$-type IVIFN such that $a_1-l_{1L}^{\nu} <0,$\\
$a_1-l_{1U}^{\prime\nu} \geq 0$ and $\tilde{A}_2=(a_2; l_{2L}^{\mu},r_{2L}^{\mu},l'^{\mu}_{2U},r'^{\mu}_{2U};l^{\nu}_{2L},r^{\nu}_{2L},l'^{\nu}_{2U},r'^{\nu}_{2U})_{LR}$ be any $LR$-type IVIFN. Then\\\\
  $\tilde{A}_1 \odot \tilde{A}_2=(a; l_{L}^{\mu},r_{L}^{\mu},l'^{\mu}_{U},r'^{\mu}_{U};l^{\nu}_{L},r^{\nu}_{L},l'^{\nu}_{U},r'^{\nu}_{U})_{LR}~~$ where\\
\noindent $a=a_1a_2,$\\
\noindent $l_L^{\mu}=a_1a_2-\mbox{min} \{(a_1-l_{1L}^{\mu})(a_2-l_{2L}^{\mu}),~(a_1+r_{1L}^{\mu})(a_2-l_{2L}^{\mu}) \},$\\
\noindent $r_L^{\mu}=\mbox{max} \{(a_1-l_{1L}^{\mu})(a_2+r_{2L}^{\mu}),~(a_1+r_{1L}^{\mu})(a_2+r_{2L}^{\mu}) \}-a_1a_2,$\\ 
\noindent $l_U^{\prime\mu}=a_1a_2-\mbox{min} \{(a_1-l_{1U}^{\prime\mu})(a_2-l_{2U}^{\prime\mu}),~(a_1+r_{1U}^{\prime\mu})(a_2-l_{2U}^{\prime\mu}) \},$\\
\noindent $r_U^{\prime\mu}=\mbox{max} \{(a_1-l_{1U}^{\prime\mu})(a_2+r_{2U}^{\prime\mu}),~(a_1+r_{1U}^{\prime\mu})(a_2+r_{2U}^{\prime\mu}) \}-a_1a_2,$\\ 
\noindent $l_L^{\nu}=a_1a_2-\mbox{min} \{(a_1-l_{1L}^{\nu})(a_2+r_{2L}^{\nu}),~(a_1+r_{1L}^{\nu})(a_2-l_{2L}^{\nu}) \},$\\
\noindent $r_L^{\nu}=\mbox{max} \{(a_1-l_{1L}^{\nu})(a_2-l_{2L}^{\nu}),~(a_1+r_{1L}^{\nu})(a_2+r_{2L}^{\nu}) \}-a_1a_2,$\\ 
\noindent $l_U^{\prime\nu}=a_1a_2-\mbox{min} \{(a_1-l_{1U}^{\prime\nu})(a_2-l_{2U}^{\prime\nu}),~(a_1+r_{1U}^{\prime\nu})(a_2-l_{2U}^{\prime\nu}) \},$\\
\noindent $r_U^{\prime\nu}=\mbox{max} \{(a_1+r_{1U}^{\prime\nu})(a_2+r_{2U}^{\prime\nu}),~(a_1-l_{1U}^{\prime\nu})(a_2+r_{2U}^{\prime\nu}) \}-a_1a_2,$\\ 
\noindent where the conditions for $LR-$type representation of $\tilde{A}_1\odot \tilde{A}_2$ are satisfied.}\\

\noindent{\bf{Proof.}} Let $\tilde{A}_1=(a_1; l_{1L}^{\mu},r_{1L}^{\mu},l'^{\mu}_{1U},r'^{\mu}_{1U};l^{\nu}_{1L},r^{\nu}_{1L},l'^{\nu}_{1U},r'^{\nu}_{1U})_{LR}$ and $\tilde{A}_2=(a_2; l_{2L}^{\mu},r_{2L}^{\mu},l'^{\mu}_{2U},r'^{\mu}_{2U}; l^{\nu}_{2L},r^{\nu}_{2L},l'^{\nu}_{2U},r'^{\nu}_{2U})_{LR}$ be two $LR$-type IVIFN with $a_1-l_{1L}^{\nu} <0,~a_1-l_{1U}^{\prime\nu} \geq 0$ and $a_2-l_{2L}^{\nu},~a_2-l_{2U}^{\prime\nu},~a_2-l_{2U}^{\prime\mu},~a_2-l_{2L}^{\mu},~a_2,~a_2+r_{2L}^{\mu},~a_2+r_{2U}^{\prime\mu},a_2+r_{2U}^{\prime\nu},~a_2+r_{2L}^{\nu}$ be any real numbers. Then, in view of Theorem $3.1$, we can write\\
%
%
%
%
%
%
%

\noindent$\big(\tilde{A}_1\odot \tilde{A}_2\big)_{\alpha}^{L}=A_{1\alpha}^{L} \times A_{2\alpha}^{L}=[a_1-l^{\mu}_{1L}L^{-1}(\alpha), a_1+r^{\mu}_{1L}R^{-1}(\alpha)] ~\times~[a_2-l^{\mu}_{2L}L^{-1}(\alpha), a_2+r^{\mu}_{2L}R^{-1}(\alpha)].$\\

\noindent $\big(\tilde{A}_1\odot \tilde{A}_2\big)_{\alpha}^{U}=A_{1\alpha}^{U} \times A_{2\alpha}^{U}= [a_1-l'^{\mu}_{1U}(L')^{-1}(\alpha), a_1+r'^{\mu}_{1U}(R')^{-1}(\alpha)]~\times~[a_2-l'^{\mu}_{2U}(L')^{-1}(\alpha), a_2+r'^{\mu}_{2U}(R')^{-1}(\alpha)].$\\

\noindent $\big(\tilde{A}_1\odot \tilde{A}_2\big)_{\beta}^{L}=A_{1\beta}^{L} \times A_{2\beta}^{L}=[a_1-l^{\nu}_{1L}L^{-1}(1-\beta), a_1+r^{\nu}_{1L}R^{-1}(1-\beta)] ~\times~[a_2-l^{\nu}_{2L}L^{-1}(1-\beta), a_2+r^{\nu}_{2L}R^{-1}(1-\beta)].$\\

\noindent $\big(\tilde{A}_1\odot \tilde{A}_2\big)_{\beta}^{U}=A_{1\beta}^{U} \times A_{2\beta}^{U}=[a_1-l'^{\nu}_{1U}(L')^{-1}(1-\beta), a_1+r'^{\nu}_{1U}(R')^{-1}(1-\beta)] ~\times~[a_2-l'^{\nu}_{2U}(L')^{-1}(1-\beta), a_2+r'^{\nu}_{2U}(R')^{-1}(1-\beta)].$\\

\noindent Since $a_1-l_{1L}^{\nu} <0$ and $a_1-l_{1U}^{\prime\nu} \geq 0$, therefore\\
 $~~~~~~~~a_1-l'^{\nu}_{1U}(L')^{-1}(1-\beta) \geq 0$ for $\beta \in [0,1]~~$ and $a_1-l^{\nu}_{1L}L^{-1}(1-\beta) \leq~(\geq)~ 0 $ for $\beta \leq ~(\geq)~ 1-L\left(\dfrac{a_1}{l_{1L}^{\nu}}\right).$\\

\noindent Now, to find the product $\tilde{A}_1 \odot \tilde{A}_2$, the following nine cases will arise and in each case, based on the sign of\\ $a_1-l^{\nu}_{1L}L^{-1}(1-\beta)$, two sub-cases are there.\\

\noindent \textbf{Case 1.} $a_2-l_{2L}^{\nu} \geq 0$.\\
Then, $a_2-l^{\nu}_{2L}L^{-1}(1-\beta) \geq 0,~a_2-l'^{\nu}_{2U}(L')^{-1}(1-\beta) \geq 0, a_2-l^{\mu}_{2L}L^{-1}(\alpha) \geq 0,~~a_2-l'^{\mu}_{2U}(L')^{-1}(\alpha) \geq 0~\forall~ \alpha, \beta \in [0,1].$\\
 
\noindent \underline{Sub-case 1.1.} $~~a_1-l^{\nu}_{1L}L^{-1}(1-\beta) \leq 0 $.\\
Now, \\
\noindent$\big(\tilde{A}_1\odot \tilde{A}_2\big)_{\alpha}^{L}=[(a_1-l^{\mu}_{1L}L^{-1}(\alpha)) (a_2-l^{\mu}_{2L}L^{-1}(\alpha)),~(a_1+r^{\mu}_{1L}R^{-1}(\alpha)) (a_2+r^{\mu}_{2L}R^{-1}(\alpha) )]$\\

\noindent $\big(\tilde{A}_1\odot \tilde{A}_2\big)_{\alpha}^{U}= [(a_1-l'^{\mu}_{1U}(L')^{-1}(\alpha))(a_2-l'^{\mu}_{2U}(L')^{-1}(\alpha)),~ (a_1+r'^{\mu}_{1U}(R')^{-1}(\alpha)) (a_2+r'^{\mu}_{2U}(R')^{-1}(\alpha))]$\\

\noindent $\big(\tilde{A}_1\odot \tilde{A}_2\big)_{\beta}^{L}=[(a_1-l^{\nu}_{1L}L^{-1}(1-\beta))(a_2+r^{\nu}_{2L}R^{-1}(1-\beta)),~ (a_1+r^{\nu}_{1L}R^{-1}(1-\beta)) (a_2+r^{\nu}_{2L}R^{-1}(1-\beta))]$\\

\noindent $\big(\tilde{A}_1\odot \tilde{A}_2\big)_{\beta}^{U}=[(a_1-l'^{\nu}_{1U}(L')^{-1}(1-\beta))(a_2-l'^{\nu}_{2U}(L')^{-1}(1-\beta)),~(a_1+r'^{\nu}_{1U}(R')^{-1}(1-\beta)) (a_2+r'^{\nu}_{2U}(R')^{-1}(1-\beta))].$\\

\noindent Further, as $L, R, L'~\mbox{and}~R'$ are decreasing functions on $[0, \infty)$ with $L(0)=R(0)=L'(0)=R'(0)=1$, there exists $\alpha_{0} \in (0,1]$, such that $L^{-1}(\alpha_0)=R^{-1}(\alpha_0)= (L')^{-1}(\alpha_0)=(R')^{-1}(\alpha_0)=1$. Hence,
\begin{equation}
\hspace{-1cm} \renewcommand{\arraystretch}{.60}
\hspace{3.5cm}  \left.\begin{array}{r@{\;}l}
\big(\tilde{A}_1\odot \tilde{A}_2\big)_{\alpha_0}^{L}=[(a_1-l^{\mu}_{1L}) (a_2-l^{\mu}_{2L}),~ (a_1+r^{\mu}_{1L} ) (a_2+r^{\mu}_{2L} )],\\\\
\big(\tilde{A}_1\odot \tilde{A}_2\big)_{\alpha_0}^{U}= [(a_1-l'^{\mu}_{1U}) (a_2-l'^{\mu}_{2U}),~ (a_1+r'^{\mu}_{1U}) (a_2+r'^{\mu}_{2U})]. \\
\end{array} \right\} 
\end{equation} 
Also, choosing $\beta_{0}=1-\alpha_0 \in (0,1]$, we have
\begin{equation}
\hspace{-1cm} \renewcommand{\arraystretch}{.60}
\hspace{3.5cm}   \left.\begin{array}{r@{\;}l}
  \big(\tilde{A}_1\odot \tilde{A}_2\big)_{\beta_0}^{L}=[(a_1-l^{\nu}_{1L})(a_2+r^{\nu}_{2L}),~(a_1+r^{\nu}_{1L}) (a_2+r^{\nu}_{2L} )],\\\\
\big(\tilde{A}_1\odot \tilde{A}_2\big)_{\beta_0}^{U}=[(a_1-l'^{\nu}_{1U}) (a_2-l'^{\nu}_{2U} ),~(a_1+r'^{\nu}_{1U}) (a_2+r'^{\nu}_{2U})].\\
\end{array} \right\} 
\end{equation} 
Taking $\alpha=1$ and $\beta=0$, we obtain
\begin{equation}
\hspace{1.8cm} \big(\tilde{A}_1\odot \tilde{A}_2\big)_{\alpha=1}^{L}=\big(\tilde{A}_1\odot \tilde{A}_2\big)_{\alpha=1}^{U}=\big(\tilde{A}_1\odot \tilde{A}_2\big)_{\beta=0}^{L}
=\big(\tilde{A}_1\odot \tilde{A}_2\big)_{\beta=0}^{U}=a_1a_2.
\end{equation}
Now, using the property that $\tilde{A}_1$ and $\tilde{A}_2$ are $LR$-type IVIFNs, we obtain\\
%
%

\noindent $a_1a_2-(a_1-l_{1U}^{\prime\mu})(a_2-l_{2U}^{\prime\mu})\geq a_1a_2-(a_1-l_{1L}^{\mu})(a_2-l_{2L}^{\mu}) \implies l_{U}^{\prime\mu} \geq l_{L}^{\mu}.$\\
$(a_1+r_{1U}^{\prime\mu})(a_2+r_{2U}^{\prime\mu})-a_1a_2 \geq (a_1+r_{1L}^{\mu})(a_2+r_{2L}^{\mu})-a_1a_2 \implies r_{U}^{\prime\mu} \geq r_{L}^{\mu}.$\\
$(a_1+r_{1L}^{\nu})(a_2+r_{2L}^{\nu})-a_1a_2 \geq (a_1+r_{1U}^{\prime\nu})(a_2+r_{2U}^{\prime\nu})-a_1a_2 \implies r_{L}^{\nu} \geq r_{U}^{\prime\nu}.$\\
$(a_1+r_{1L}^{\nu})(a_2+r_{2L}^{\nu})-a_1a_2 \geq (a_1+r_{1L}^{\mu})(a_2+r_{2L}^{\mu})-a_1a_2 \implies r_{L}^{\nu} \geq r_{L}^{\mu}.$\\
$a_1a_2-(a_1-l_{1U}^{\prime\nu})(a_2-l_{2U}^{\prime\nu}) \geq a_1a_2-(a_1-l_{1U}^{\prime\mu})(a_2-l_{2U}^{\prime\mu}) \implies l_{U}^{\prime\nu} \geq l_{U}^{\prime\mu}.$\\
$(a_1+r_{1U}^{\prime\nu})(a_2+r_{2U}^{\prime\nu})-a_1a_2 \geq (a_1+r_{1U}^{\prime\mu})(a_2+r_{2U}^{\prime\mu})-a_1a_2 \implies r_{U}^{\prime\nu} \geq r_{U}^{\prime\mu}.$\\

\noindent Further, since $ a_2-l_{2L}^{\nu} \leq a_2+r_{2L}^{\nu},~a_1-l_{1L}^{\nu} \leq a_1-l_{1U}^{\prime\nu}$ and $a_2-l_{2L}^{\nu} \leq a_2-l_{2U}^{\prime\nu}$, therefore\\\\
$a_1a_2-(a_1-l_{1U}^{\prime\nu})(a_2-l_{2U}^{\prime\nu}) \leq a_1a_2-(a_1-l_{1L}^{\nu})(a_2+r_{2L}^{\nu}),$ that is, $~l_U^{\prime\nu} \leq l_L^{\nu}$\\
 and from $a_1-l_{1L}^{\nu} <0,~a_2-l_{2L}^{\nu} \leq a_2+r_{2L}^{\nu},a_1-l_{1L}^{\nu} \leq a_1-l_{1L}^{\mu},~a_2-l_{2L}^{\nu} \leq a_2-l_{2L}^{\mu}$, we have\\
$ a_1a_2-(a_1-l_{1L}^{\mu})(a_2-l_{2L}^{\mu}) \leq a_1a_2-(a_2+r_{2L}^{\nu})(a_1-l_{1L}^{\nu}),$ that is, $l_L^{\mu} \leq l_L^{\nu}.$\\
\noindent Moreover, $l_{2L}^{\mu}(a_1-l_{1L}^{\mu})+a_2l_{1L}^{\mu} \geq 0~~~$ and $~~a_1r_{2L}^{\mu}+a_2r_{1L}^{\mu}+r_{1L}^{\mu}r_{2L}^{\mu} \geq 0~~$ yield $~~l_L^{\mu} \geq 0$ and $~r_L^{\mu} \geq 0$, respectively.\\

\noindent Finally, in view of these inequalities and combining expressions $(21),(22)$ and $(23)$, we get\\

\noindent $\tilde{A}_1 \odot \tilde{A}_2=\big(a_1a_2;a_1a_2-\left(a_1-l^{\mu}_{1L}\right) \left(a_2-l^{\mu}_{2L}\right),\left(a_1+r^{\mu}_{1L} \right) \left(a_2+r^{\mu}_{2L} \right)-a_1a_2,a_1a_2-\left(a_1-l'^{\mu}_{1U}\right) \left(a_2-l'^{\mu}_{2U}\right),$\\
$~~~~~~~~~~~~~~~~~~~\left(a_1+r'^{\mu}_{1U}\right) \left(a_2+r'^{\mu}_{2U}\right)-a_1a_2;a_1a_2-\left(a_1-l^{\nu}_{1L}\right)\left(a_2+r^{\nu}_{2L}\right),\left(a_1+r^{\nu}_{1L}\right) \left(a_2+r^{\nu}_{2L} \right)-a_1a_2,$\\
$~~~~~~~~~~~~~~~~~~~~~a_1a_2-\left(a_1-l'^{\nu}_{1U}\right) \left(a_2-l'^{\nu}_{2U} \right),\left(a_1+r'^{\nu}_{1U}\right) \left(a_2+r'^{\nu}_{2U}\right)-a_1a_2\big)_{LR}.$\\

\noindent \underline{Sub-case 1.2.} $~~~a_1-l^{\nu}_{1L}L^{-1}(1-\beta) \geq 0 $.\\
Now,\\
\noindent $\big(\tilde{A}_1\odot \tilde{A}_2\big)_{\beta}^{L}=\left[\left(a_1-l^{\nu}_{1L}L^{-1}(1-\beta)\right)\left(a_2-l^{\nu}_{2L}L^{-1}(1-\beta)\right),~ \left(a_1+r^{\nu}_{1L}R^{-1}(1-\beta)\right) \left(a_2+r^{\nu}_{2L}R^{-1}(1-\beta) \right)\right].$\\

\noindent Choosing $ \beta =\beta_0$, we have
\begin{equation}
 \hspace{2cm} \big(\tilde{A}_1\odot \tilde{A}_2\big)_{\beta _0}^{L}=\left[\left(a_1-l^{\nu}_{1L}\right)\left(a_2-l^{\nu}_{2L} \right),\left(a_1+r^{\nu}_{1L} \right) \left(a_2+r^{\nu}_{2L} \right) \right].
\end{equation} 
Hence, we have $l_L^{\nu}=a_1a_2-\left(a_1-l^{\nu}_{1L}\right)\left(a_2-l^{\nu}_{2L} \right)$ and remaining spreads of $\tilde{A}_1 \odot \tilde{A}_2$ will be same as derived in Sub-case 1.1. Further,\\\\
 $~~~~~~~~~~~~~~~~~~~~~~~~~~~~~~~~(a_1-l^{\nu}_{1L})(a_2-l^{\nu}_{2L}) \leq (a_1-l^{\prime\nu}_{1U})(a_2-l^{\prime\nu}_{2U}) \implies l^{\nu}_{L} \geq l^{\prime\nu}_{U}~~~$ and\\
$~~~~~~~~~~~~~~~~~~~~~~~~~~~~~~~~(a_1-l^{\nu}_{1L})(a_2-l^{\nu}_{2L}) \leq (a_1-l^{\mu}_{1L})(a_2-l^{\mu}_{2L}) \implies l^{\nu}_{L} \geq l^{\mu}_{L}.$\\

\noindent Finally, combining Sub-cases 1.1 and 1.2, we get\\

\noindent $\tilde{A}_1 \odot \tilde{A}_2=(a; l_{L}^{\mu},r_{L}^{\mu},l'^{\mu}_{U},r'^{\mu}_{U};l^{\nu}_{L},r^{\nu}_{L},l'^{\nu}_{U},r'^{\nu}_{U})_{LR}$ where\\
\noindent $a=a_1a_2,$\\
\noindent $l_L^{\mu}=a_1a_2-(a_1-l_{1L}^{\mu})(a_2-l_{2L}^{\mu}),$\\
\noindent $r_L^{\mu}=(a_1+r_{1L}^{\mu})(a_2+r_{2L}^{\mu})-a_1a_2,$\\ 
\noindent $l_U^{\prime\mu}=a_1a_2-(a_1-l_{1U}^{\prime\mu})(a_2-l_{2U}^{\prime\mu}),$\\
\noindent $r_U^{\prime\mu}=(a_1+r_{1U}^{\prime\mu})(a_2+r_{2U}^{\prime\mu})-a_1a_2,$\\ 
\noindent $l_L^{\nu}=a_1a_2-\mbox{min} \{(a_1-l_{1L}^{\nu})(a_2+r_{2L}^{\nu}),(a_1-l^{\nu}_{1L})(a_2-l^{\nu}_{2L}) \},$\\
\noindent $r_L^{\nu}= (a_1+r_{1L}^{\nu})(a_2+r_{2L}^{\nu}) \}-a_1a_2,$\\ 
\noindent $l_U^{\prime\nu}=a_1a_2-(a_1-l_{1U}^{\prime\nu})(a_2-l_{2U}^{\prime\nu}),$\\
\noindent $r_U^{\prime\nu}=(a_1+r_{1U}^{\prime\nu})(a_2+r_{2U}^{\prime\nu})-a_1a_2.$\\

\noindent \textbf{Case 2.} $~a_2-l_{2L}^{\nu} < 0$ and $a_2-l_{2U}^{\prime\nu} \geq 0$.\\
Then, we have
$$a_2-l^{\nu}_{2L}L^{-1}(1-\beta) \leq~(\geq)~ 0~~\mbox{for}~~\beta \leq~(\geq)~ 1-L\left(\dfrac{a_2}{l_{2L}^{\nu}}\right).$$
 
\noindent \underline{Sub-case 2.1.} $~~a_1-l^{\nu}_{1L}L^{-1}(1-\beta) \leq 0 $.\\
Then, $\big(\tilde{A}_1\odot \tilde{A}_2\big)_{\alpha}^{L},~\big(\tilde{A}_1\odot \tilde{A}_2\big)_{\alpha}^{U},~\big(\tilde{A}_1\odot \tilde{A}_2\big)_{\beta}^{U}$ remain same as in Sub-case $1.1$, however, the expression for $\big(\tilde{A}_1\odot \tilde{A}_2\big)_{\beta}^{L}$ will be changed and is given by\\

\noindent $\big(\tilde{A}_1\odot \tilde{A}_2\big)_{\beta}^{L}=A_{1\beta}^{L} \times A_{2\beta}^{L}$\\
$~~~~~~~~~~~~~~~~~~~~~=\left[\left(a_1-l^{\nu}_{1L}L^{-1}(1-\beta)\right), \left(a_1+r^{\nu}_{1L}R^{-1}(1-\beta)\right)\right]~\times ~\left[\left(a_2-l^{\nu}_{2L}L^{-1}(1-\beta)\right), \left(a_2+r^{\nu}_{2L}R^{-1}(1-\beta) \right)\right]$\\
$~~~~~~~~~~~~~~~~~~~~~=\left[\mbox{min} \left\{ \left(a_1-l^{\nu}_{1L}L^{-1}(1-\beta)\right)\left(a_2+r^{\nu}_{2L}R^{-1}(1-\beta) \right),\right.\left(a_2-l^{\nu}_{2L}L^{-1}(1-\beta)\right)\left(a_1+r^{\nu}_{1L}R^{-1}(1-\beta)\right)\right\},$\\
$~~~~~~~~~~~~~~~~~~~~~~\mbox{max}\left\{\left(a_1-l^{\nu}_{1L}L^{-1}(1-\beta)\right)\left(a_2-l^{\nu}_{2L}L^{-1}(1-\beta)\right),\left.\left(a_1+r^{\nu}_{1L}R^{-1}(1-\beta)\right)\left(a_2+r^{\nu}_{2L}R^{-1}(1-\beta) \right) \right\}\right].$\\

\noindent Choosing $ \beta =\beta_0$, we get
\begin{equation}
 \renewcommand{\arraystretch}{.60}
 \hspace{2cm} \left.\begin{array}{r@{\;}l}
  \big(\tilde{A}_1\odot \tilde{A}_2\big)_{\beta_0}^{L}=\left[\mbox{min} \left\{ \left(a_1-l^{\nu}_{1L}\right)\left(a_2+r^{\nu}_{2L} \right),\right. \left(a_2-l^{\nu}_{2L}\right) \left(a_1+r^{\nu}_{1L} \right) \right\},\\\\
\mbox{max}\left\{\left(a_1-l^{\nu}_{1L}\right) \left.\left(a_2-l^{\nu}_{2L}\right),\left(a_1+r^{\nu}_{1L}\right)\left(a_2+r^{\nu}_{2L} \right) \right\}\right].
\end{array} \right\} 
\end{equation} 
Hence,\\\\
 $~~~~~~~~~~~~~~~~~~~~~~~~~~~~~~~l_L^{\nu}=a_1a_2-\mbox{min} \left\{ \left(a_1-l^{\nu}_{1L}\right)\left(a_2+r^{\nu}_{2L} \right),~ \left(a_2-l^{\nu}_{2L}\right) \left(a_1+r^{\nu}_{1L} \right) \right\},$\\

\noindent $~~~~~~~~~~~~~~~~~~~~~~~~~~~~~~r_L^{\nu}=\mbox{max}\left\{\left(a_1-l^{\nu}_{1L}\right)
\left(a_2-l^{\nu}_{2L}\right),\left(a_1+r^{\nu}_{1L}\right)\left(a_2+r^{\nu}_{2L} \right)\right\}-a_1a_2.$\\

\noindent Further, since $a_1-l_{1L}^{\nu} \leq a_1-l_{1U}^{\prime\nu},~a_2-l_{2U}^{\prime\nu} \leq a_2+r_{2U}^{\prime\nu},~a_2+r_{2U}^{\prime\nu} \leq a_2+r_{2L}^{\nu}$ and $a_1-l_{1L}^{\nu}<0$, therefore
\begin{equation}
\hspace{3.5cm} (a_1-l_{1U}^{\prime\nu})(a_2-l_{2U}^{\prime\nu}) \geq (a_1-l_{1L}^{\nu})(a_2+r_{2L}^{\nu}).
\end{equation}
Also, $~a_1-l_{1U}^{\prime\nu} \leq a_1+r_{1U}^{\prime\nu} \leq a_1+r_{1L}^{\nu}~\mbox{and}~a_2-l_{2L}^{\nu} <0 $ implies
\begin{equation}
\hspace{3.5cm} (a_1-l_{1U}^{\prime\nu})(a_2-l_{2U}^{\prime\nu}) \geq (a_1+r_{1L}^{\nu})(a_2-l_{2L}^{\nu}).
\end{equation}
It follows from the inequalities $(26)$ and $(27)$ that \\

\noindent $~~~~~~~~~~~~~~~~~~~(a_1-l_{1U}^{\prime\nu})(a_2-l_{2U}^{\prime\nu}) \geq \mbox{min} \{(a_1-l_{1L}^{\nu})(a_2+r_{2L}^{\nu}),(a_1+r_{1L}^{\nu})(a_2-l_{2L}^{\nu}) \}$, that is, $l_U^{\prime\nu} \leq l_L^{\nu}$.\\

\noindent Now, if $(a_1-l^{\nu}_{1L})(a_2-l^{\nu}_{2L}) \leq (a_1+r^{\nu}_{1L})(a_2+r^{\nu}_{2L}),$ then 
\begin{equation}
\hspace{3.5cm} (a_1+r_{1U}^{\prime\nu})(a_2+r_{2U}^{\prime\nu}) \leq (a_1+r^{\nu}_{1L})(a_2+r^{\nu}_{2L})
\end{equation}
and if $(a_1-l^{\nu}_{1L})(a_2-l^{\nu}_{2L}) \geq (a_1+r^{\nu}_{1L})(a_2+r^{\nu}_{2L}),$ then
\begin{equation}
\hspace{2.15cm} (a_1+r_{1U}^{\prime\nu})(a_2+r_{2U}^{\prime\nu}) \leq (a_1+r^{\nu}_{1L})(a_2+r^{\nu}_{2L}) \leq (a_1-l^{\nu}_{1L})(a_2-l^{\nu}_{2L}).
\end{equation}

\noindent Thus, the inequalities $(28)$ and $(29)$ yield\\

\noindent $~~~~~~~~~~~~~~~~(a_1+r_{1U}^{\prime\nu})(a_2+r_{2U}^{\prime\nu}) \leq \mbox{max} \{(a_1+r^{\nu}_{1L})(a_2+r^{\nu}_{2L}),~ (a_1-l^{\nu}_{1L})(a_2-l^{\nu}_{2L}) \}$, that is, $r_U^{\prime\nu} \leq r_L^{\nu}$.\\

\noindent Further, $a_1-l^{\nu}_{1L} \leq a_1-l^{\mu}_{1L},~a_1-l^{\nu}_{1L} <0$ give
\begin{equation}
\hspace{3.5cm}(a_1-l^{\mu}_{1L})(a_2-l^{\mu}_{2L}) \geq (a_1-l^{\nu}_{1L})(a_2+r^{\nu}_{2L})
\end{equation}
\vspace{-0.7cm}
\begin{equation} \hspace{-0.8cm}\mbox{and}  \hspace{2.2cm} (a_1-l^{\mu}_{1L})(a_2-l^{\mu}_{2L}) \geq (a_1-l^{\mu}_{1L})(a_2-l^{\nu}_{2L})
\geq (a_1+r^{\nu}_{1L})(a_2-l^{\nu}_{2L}).
\end{equation}
It follows from the expressions $(30)$ and $(31)$ that\\

\noindent $~~~~~~~~~~~~~~~~(a_1-l^{\mu}_{1L})(a_2-l^{\mu}_{2L}) \geq \mbox{min}\{(a_1-l^{\nu}_{1L})(a_2+r^{\nu}_{2L}),~(a_1+r^{\nu}_{1L})(a_2-l^{\nu}_{2L})\}$, that is, $l_L^{\mu} \leq l_L^{\nu}$.\\

\noindent Now, if $(a_1-l^{\nu}_{1L})(a_2-l^{\nu}_{2L}) \leq (a_1+r^{\nu}_{1L})(a_2+r^{\nu}_{2L})$, then 
\begin{equation}
\hspace{3cm}(a_1+r_{1L}^{\mu})(a_2+r_{2L}^{\mu}) \leq (a_1+r^{\nu}_{1L})(a_2+r^{\nu}_{2L})
\end{equation}
and if $(a_1-l^{\nu}_{1L})(a_2-l^{\nu}_{2L}) \geq (a_1+r^{\nu}_{1L})(a_2+r^{\nu}_{2L})$, then
\begin{equation}
\hspace{2.1cm} (a_1+r_{1L}^{\mu})(a_2+r_{2L}^{\mu}) \leq (a_1+r^{\nu}_{1L})(a_2+r^{\nu}_{2L}) \leq (a_1-l^{\nu}_{1L})(a_2-l^{\nu}_{2L}).
\end{equation}

\noindent Hence, the inequalities $(32)$ and $(33)$ give\\

\noindent $~~~~~~~~~~~~(a_1+r_{1L}^{\mu})(a_2+r_{2L}^{\mu}) \leq \mbox{max} \{(a_1+r^{\nu}_{1L})(a_2+r^{\nu}_{2L}),~(a_1-l^{\nu}_{1L})(a_2-l^{\nu}_{2L}) \}$, that is, $r_L^{\mu} \leq r_L^{\nu}$.\\

\noindent \underline{Sub-case 2.2.} $~~a_1-l^{\nu}_{1L}L^{-1}(1-\beta) \geq 0 $.\\
The proof of this part follows on the lines of Sub-case 1.1.\\

\noindent Hence, now combining Cases $1$ and $2$, we get\\
\noindent $\tilde{A}_1 \odot \tilde{A}_2=(a; l_{L}^{\mu},r_{L}^{\mu},l'^{\mu}_{U},r'^{\mu}_{U};l^{\nu}_{L},r^{\nu}_{L},l'^{\nu}_{U},r'^{\nu}_{U})_{LR}$ where\\
\noindent $a=a_1a_2,$\\
\noindent $l_L^{\mu}=a_1a_2-(a_1-l_{1L}^{\mu})(a_2-l_{2L}^{\mu}),$\\
\noindent $r_L^{\mu}=(a_1+r_{1L}^{\mu})(a_2+r_{2L}^{\mu})-a_1a_2,$\\
\noindent $l_U^{\prime\mu}=a_1a_2-(a_1-l_{1U}^{\prime\mu})(a_2-l_{2U}^{\prime\mu}),$\\
\noindent $r_U^{\prime\mu}=(a_1+r_{1U}^{\prime\mu})(a_2+r_{2U}^{\prime\mu})-a_1a_2,$\\ 
\noindent $l_L^{\nu}=a_1a_2-\mbox{min} \left\{ \left(a_1-l^{\nu}_{1L}\right)\left(a_2+r^{\nu}_{2L} \right),~ \left(a_2-l^{\nu}_{2L}\right) \left(a_1+r^{\nu}_{1L} \right) \right\},$\\
\noindent $r_L^{\nu}= \mbox{max}\left\{\left(a_1-l^{\nu}_{1L}\right)
\left(a_2-l^{\nu}_{2L}\right),\left(a_1+r^{\nu}_{1L}\right)\left(a_2+r^{\nu}_{2L} \right)\right\}-a_1a_2,$\\ 
\noindent $l_U^{\prime\nu}=a_1a_2-(a_1-l_{1U}^{\prime\nu})(a_2-l_{2U}^{\prime\nu}),$\\
\noindent $r_U^{\prime\nu}=(a_1+r_{1U}^{\prime\nu})(a_2+r_{2U}^{\prime\nu})-a_1a_2.$\\

\noindent \textbf{Case 3.} $a_2-l_{2U}^{\prime\nu} <0$ and $a_2-l_{2U}^{\prime\mu}  \geq 0$.\\
Then, $~~~~~~~~~~~~~~~~~~~~~~~~~a_2-l_{2U}^{\prime\nu}(L^{\prime})^{-1}(1-\beta) \leq~(\geq)~ 0~$ for $~\beta \leq~(\geq)~ 1-L^{\prime}\left(\dfrac{a_2}{l_{2U}^{\prime\nu}}\right)$. \\

\noindent \underline{Sub-case 3.1.}  $~~a_1-l^{\nu}_{1L}L^{-1}(1-\beta) \leq 0 $.\\
Then, the expressions $\big(\tilde{A}_1\odot \tilde{A}_2\big)_{\alpha}^{L},~\big(\tilde{A}_1\odot \tilde{A}_2\big)_{\alpha}^{U},\big(\tilde{A}_1\odot \tilde{A}_2\big)_{\beta}^{L}$ are same as discussed in Case $2$ but $\big(\tilde{A}_1\odot \tilde{A}_2\big)_{\beta}^{U}$ is given by:\\

\noindent $\big(\tilde{A}_1\odot \tilde{A}_2\big)_{\beta}^{U}=A_{1\beta}^{U} \times A_{2\beta}^{U}$\\
$~~~~~~~~~~~~~~~~~~~~~~=\left[\left(a_1+r'^{\nu}_{1U}(R')^{-1}(1-\beta)\right)\left(a_2-l'^{\nu}_{2U}(L')^{-1}(1-\beta)\right),\left(a_1+r'^{\nu}_{1U}(R')^{-1}(1-\beta)\right)\left(a_2+r'^{\nu}_{2U}(R')^{-1}(1-\beta)\right)\right].$\\

\noindent Taking $\beta=\beta_0$, we get
$$\big(\tilde{A}_1\odot \tilde{A}_2\big)_{\beta_0}^{U}=\left[\left(a_1+r'^{\nu}_{1U}\right)\left(a_2-l'^{\nu}_{2U}\right),\left(a_1+r'^{\nu}_{1U}\right)\left(a_2+r'^{\nu}_{2U}\right)\right]$$
which gives $~~~~~~~~~~~~~~~~~~~~~~~l_U^{\prime\nu}=a_1a_2-(a_1+r'^{\nu}_{1U})(a_2-l'^{\nu}_{2U}).$\\
Also, we have\\
$l_L^{\nu}=a_1a_2-\mbox{min} \left\{ \left(a_1-l^{\nu}_{1L}\right)\left(a_2+r^{\nu}_{2L} \right),~ \left(a_2-l^{\nu}_{2L}\right) \left(a_1+r^{\nu}_{1L} \right) \right\}~~$ and $~~l_U^{\prime\mu}=a_1a_2-(a_1-l_{1U}^{\prime\mu})(a_2-l_{2U}^{\prime\mu})$.\\

\noindent Now, we claim that $l_L^{\nu} \geq l_U^{\prime\nu}$ and $l_U^{\prime\nu} \geq l_U^{\prime\mu}.$\\\\
Since $a_2-l_{2L}^{\nu}<0,~a_2-l_{2L}^{\nu} \leq a_2-l_{2U}^{\prime\nu}$ and $a_1+r_{1U}^{\prime\nu} \leq a_1+r_{1L}^{\nu}$, therefore, we obtain
\begin{equation}
\hspace{4cm} (a_1+r_{1U}^{\prime\nu})(a_2-l_{2U}^{\prime\nu}) \geq (a_1+r_{1L}^{\nu}) (a_2-l_{2L}^{\nu}),
\end{equation}
\vspace{-0.8cm}
\begin{equation}
\hspace{4cm} (a_1+r_{1U}^{\prime\nu})(a_2-l_{2U}^{\prime\nu}) \geq (a_1-l_{1L}^{\nu})(a_2+r_{2L}^{\nu}).
\end{equation}

\noindent Thus, the inequalities $(34)$ and $(35)$ yield\\\\
$~~~~~~~~~~~~~~~~~~~(a_1+r_{1U}^{\prime\nu})(a_2-l_{2U}^{\prime\nu}) \geq \mbox{min} \{(a_1+r_{1L}^{\nu}) (a_2-l_{2L}^{\nu}),~(a_1-l_{1L}^{\nu})(a_2+r_{2L}^{\nu}) \}$ that is, $l_U^{\prime\nu} \leq l_L^{\nu}$.\\

\noindent Further, $a_2-l_{2U}^{\prime\nu}<0,~ a_1-l_{1U}^{\prime\mu} \leq a_1+r_{1U}^{\prime\nu},~a_2-l_{2U}^{\prime\nu} \leq a_2-l_{2U}^{\prime\mu}$ implies\\\\
$(a_1+r_{1U}^{\prime\nu})(a_2-l_{2U}^{\prime\nu}) \leq (a_1-l_{1U}^{\prime\mu})(a_2-l_{2U}^{\prime\nu})\leq (a_1-l_{1U}^{\prime\mu})(a_2-l_{2U}^{\prime\mu})$ that is, $l_U^{\prime\nu} \geq l_U^{\prime\mu}$.\\

\noindent \underline{Sub-case 3.2.}  $~~a_1-l^{\nu}_{1L}L^{-1}(1-\beta) \geq 0 $.\\
The proof of this part follows on the lines of Sub-cases $1.1$ and $2.2$.\\\\
Hence, it is concluded from Cases $1$~\textendash~$3$ that\\
\noindent $\tilde{A}_1 \odot \tilde{A}_2=(a; l_{L}^{\mu},r_{L}^{\mu},l'^{\mu}_{U},r'^{\mu}_{U};l^{\nu}_{L},r^{\nu}_{L},l'^{\nu}_{U},r'^{\nu}_{U})_{LR}$ where\\
\noindent $a=a_1a_2,$\\
\noindent $l_L^{\mu}=a_1a_2-(a_1-l_{1L}^{\mu})(a_2-l_{2L}^{\mu}),$\\
\noindent $r_L^{\mu}=(a_1+r_{1L}^{\mu})(a_2+r_{2L}^{\mu})-a_1a_2,$\\
\noindent $l_U^{\prime\mu}=a_1a_2-(a_1-l_{1U}^{\prime\mu})(a_2-l_{2U}^{\prime\mu}),$\\
\noindent $r_U^{\prime\mu}=(a_1+r_{1U}^{\prime\mu})(a_2+r_{2U}^{\prime\mu})-a_1a_2,$\\ 
\noindent $l_L^{\nu}=a_1a_2-\mbox{min} \left\{ \left(a_1-l^{\nu}_{1L}\right)\left(a_2+r^{\nu}_{2L} \right),~ \left(a_2-l^{\nu}_{2L}\right) \left(a_1+r^{\nu}_{1L} \right) \right\},$\\
\noindent $r_L^{\nu}= \mbox{max}\left\{\left(a_1-l^{\nu}_{1L}\right)
\left(a_2-l^{\nu}_{2L}\right),\left(a_1+r^{\nu}_{1L}\right)\left(a_2+r^{\nu}_{2L} \right)\right\}-a_1a_2,$\\ 
\noindent $l_U^{\prime\nu}=a_1a_2-\mbox{min} \left\{(a_1-l_{1U}^{\prime\nu})(a_2-l_{2U}^{\prime\nu}),~(a_1+r_{1U}^{\prime\nu})(a_2-l_{2U}^{\prime\nu})\right\},$\\
\noindent $r_U^{\prime\nu}=(a_1+r_{1U}^{\prime\nu})(a_2+r_{2U}^{\prime\nu})-a_1a_2.$\\

\noindent \textbf{Case 4.}  $a_2-l_{2U}^{\prime\mu} <0$ and $a_2-l_{2L}^{\mu} \geq 0$.\\
Then, $~~~~~~~~~~~~~~~~~~~~~~~~a_2-l_{2U}^{\prime\mu}(L^{\prime})^{-1}(\alpha) \leq~(\geq)~ 0~$ for $~\alpha \geq~(\leq)~ L^{\prime}\left(\dfrac{a_2}{l_{2U}^{\prime\mu}}\right)$. \\

\noindent \underline{Sub-case 4.1.} $~~a_1-l^{\nu}_{1L}L^{-1}(1-\beta) \leq 0 $.\\
Then, the expressions $\big(\tilde{A}_1\odot \tilde{A}_2\big)_{\alpha}^{L},~\big(\tilde{A}_1\odot \tilde{A}_2\big)_{\beta}^{L},~\big(\tilde{A}_1\odot \tilde{A}_2\big)_{\beta}^{U}$ are same as in Cases $2$ and $3$ but $\big(\tilde{A}_1\odot \tilde{A}_2\big)_{\alpha}^{U}$ is given by\\

\noindent $\big(\tilde{A}_1\odot \tilde{A}_2\big)_{\alpha}^{U}=A_{1\alpha}^{U} \times A_{2\alpha}^{U}$\\
$~~~~~~~~~~~~~~~~~~~~~~=\left[\left(a_1+r'^{\mu}_{1U}(R')^{-1}(\alpha)\right)\left(a_2-l'^{\mu}_{2U}(L')^{-1}(\alpha)\right),~\left(a_1+r'^{\mu}_{1U}(R')^{-1}(\alpha)\right)\left(a_2+r'^{\mu}_{2U}(R')^{-1}(\alpha)\right)\right].$\\

\noindent Choosing $\alpha=\alpha_0$, we get\\\\
\noindent $~~~~~~~~~~~~~~~~~~~~~~~~~~~~~~~~~~~\big(\tilde{A}_1\odot \tilde{A}_2\big)_{\alpha_0}^{U}=\left[\left(a_1+r'^{\mu}_{1U}\right)\left(a_2-l'^{\mu}_{2U}\right),~\left(a_1+r'^{\mu}_{1U}\right)\left(a_2+r'^{\mu}_{2U}\right)\right].$\\\\
This further yields $~~~~~~~~~~~~~~~~~~~~l_U^{\prime\mu}=a_1a_2-(a_1+r'^{\mu}_{1U})(a_2-l'^{\mu}_{2U})$.\\

\noindent Moreover, from the preceding Case $3$, we have\\\\
$~~~~~~~~~~~l_L^{\mu}=a_1a_2-(a_1-l_{1L}^{\mu})(a_2-l_{2L}^{\mu}),~~l_U^{\prime\nu}=a_1a_2-(a_1+r_{1U}^{\prime\nu})(a_2-l_{2U}^{\prime\nu})$.\\

\noindent Next, it remain to prove that $~~~l_U^{\prime\mu} \geq l_L^{\mu}~~$ and $~~l_U^{\prime\nu} \geq l_U^{\prime\mu}.$\\\\
\noindent Since $a_2-l_{2U}^{\prime\mu}<0,~a_2-l_{2U}^{\prime\mu} \leq a_2-l_{2L}^{\mu}~$ and $~a_1-l_{1L}^{\mu} \leq a_1+r_{1L}^{\mu} \leq a_1+r_{1U}^{\prime\mu}$, therefore\\\\
$~~~~~~~~~~~~~~~~~~~~~a_1a_2-(a_2-l_{2L}^{\mu})(a_1-l_{1L}^{\mu}) \leq a_1a_2-(a_2-l_{2U}^{\prime\mu})(a_1+r_{1U}^{\prime\mu}) \implies l_L^{\mu} \leq l_U^{\prime\mu}.$\\

\noindent Further, $a_2-l_{2U}^{\prime\nu} \leq a_2-l_{2U}^{\prime\mu},~a_2-l_{2U}^{\prime\mu} \leq 0$ yield\\\\
$~~~~~~~~~~~~~~~~~~~~~a_1a_2-(a_2-l_{2U}^{\prime\mu})(a_1+r_{1U}^{\prime\mu})\leq a_1a_2-(a_2-l_{2U}^{\prime\nu})(a_1+r_{1U}^{\prime\nu}) \implies l_U^{\prime\mu} \leq l_U^{\prime\nu}.$\\

\noindent \underline{Sub-case 4.2.} $~~a_1-l^{\nu}_{1L}L^{-1}(1-\beta) \geq 0 $.\\
It follows on the lines of Sub-case $3.2$.\\

\noindent Hence, clubbing the Cases $1$~\textendash~$4$, we obtain\\
\noindent $\tilde{A}_1 \odot \tilde{A}_2=(a; l_{L}^{\mu},r_{L}^{\mu},l'^{\mu}_{U},r'^{\mu}_{U};l^{\nu}_{L},r^{\nu}_{L},l'^{\nu}_{U},r'^{\nu}_{U})_{LR}$ where\\
\noindent $a=a_1a_2,$\\
\noindent $l_L^{\mu}=a_1a_2-(a_1-l_{1L}^{\mu})(a_2-l_{2L}^{\mu}),$\\
\noindent $r_L^{\mu}=(a_1+r_{1L}^{\mu})(a_2+r_{2L}^{\mu})-a_1a_2,$\\
\noindent $l_U^{\prime\mu}=a_1a_2-\mbox{min} \left\{(a_1-l_{1U}^{\prime\mu})(a_2-l_{2U}^{\prime\mu}),~ (a_1+r'^{\mu}_{1U})(a_2-l'^{\mu}_{2U}) \right\},$\\
\noindent $r_U^{\prime\mu}=(a_1+r_{1U}^{\prime\mu})(a_2+r_{2U}^{\prime\mu})-a_1a_2,$\\ 
\noindent $l_L^{\nu}=a_1a_2-\mbox{min} \left\{ \left(a_1-l^{\nu}_{1L}\right)\left(a_2+r^{\nu}_{2L} \right),~ \left(a_2-l^{\nu}_{2L}\right) \left(a_1+r^{\nu}_{1L} \right) \right\},$\\
\noindent $r_L^{\nu}= \mbox{max}\left\{\left(a_1-l^{\nu}_{1L}\right)
\left(a_2-l^{\nu}_{2L}\right),\left(a_1+r^{\nu}_{1L}\right)~\left(a_2+r^{\nu}_{2L} \right)\right\}-a_1a_2,$\\ 
\noindent $l_U^{\prime\nu}=a_1a_2-\mbox{min} \left\{(a_1-l_{1U}^{\prime\nu})(a_2-l_{2U}^{\prime\nu}),~(a_1+r_{1U}^{\prime\nu})(a_2-l_{2U}^{\prime\nu})\right\},$\\
\noindent $r_U^{\prime\nu}=(a_1+r_{1U}^{\prime\nu})(a_2+r_{2U}^{\prime\nu})-a_1a_2.$\\

Now, for the succeeding cases, we have only derived the expressions for newer (or changed) spreads of the product $\tilde{A}_1 \odot \tilde{A}_2$ because the proofs for the bounds on the other spreads can be carried out in the similar pattern as is discussed in the preceding four cases.\\

\noindent \textbf{Case 5.} $a_2-l_{2L}^{\mu} <0$ and $a_2 \in \mathbb{R}$.\\
Then, $~~~~~~~~~~~~~~~~~~~~~~~~~~~a_2-l_{2L}^{\mu}L^{-1}(\alpha) \leq~(\geq)~ 0~$ for $~ \alpha \geq~(\leq)~ L\left(\dfrac{a_2}{l_{2L}^{\mu}}\right)$.\\
 
\noindent \underline{Sub-case 5.1.} $~~a_1-l^{\nu}_{1L}L^{-1}(1-\beta) \leq 0$.\\\\
Then, \\
\noindent $~~~~~~~~~~\big(\tilde{A}_1\odot \tilde{A}_2\big)_{\alpha}^{L}=A_{1\alpha}^{L} \times A_{2\alpha}^{L}=[(a_1+r^{\mu}_{1L}R^{-1}(\alpha))(a_2-l^{\mu}_{2L}L^{-1}(\alpha)),~(a_1+r^{\mu}_{1L}R^{-1}(\alpha))(a_2+r^{\mu}_{2L}R^{-1}(\alpha))].$\\\\
\noindent Taking $\alpha=\alpha_0$, we get
$$\big(\tilde{A}_1\odot \tilde{A}_2\big)_{\alpha_0}^{L}=[(a_1+r^{\mu}_{1L})(a_2-l^{\mu}_{2L}),(a_1+r^{\mu}_{1L})(a_2+r^{\mu}_{2L})].$$
It further yields $~~~~~~~~~~~~~~~~~~~~~~~~~~~~~~l_L^{\mu}=a_1a_2-(a_1+r^{\mu}_{1L})(a_2-l^{\mu}_{2L}).$\\\\
Moreover, from above Cases, we have\\\\
$~~~~~~~~~~~~~~~~~~~~~~~l_U^{\prime\mu}=a_1a_2-\mbox{min} \left\{(a_1-l_{1U}^{\prime\mu})(a_2-l_{2U}^{\prime\mu}),(a_1+r'^{\mu}_{1U})(a_2-l'^{\mu}_{2U}) \right\},$\\\\
$~~~~~~~~~~~~~~~~~~~~~~~l_L^{\nu}=a_1a_2-\mbox{min} \{ \left(a_1-l^{\nu}_{1L}\right)\left(a_2+r^{\nu}_{2L} \right), \left(a_2-l^{\nu}_{2L}\right) \left(a_1+r^{\nu}_{1L} \right) \}$.\\\\
Next, the claim that $l_U^{\prime\mu} \geq l_L^{\mu}$ and $l_L^{\nu} \geq l_L^{\mu}$ can be proved following the lines of Sub-case $3.1$.\\

\noindent \underline{Sub-case 5.2.} $~~a_1-l^{\nu}_{1L}L^{-1}(1-\beta) \geq 0$.\\
The proof can be carried out following the lines of Sub-case $1.1$.\\

\noindent \textbf{Case 6.} $a_2+r_{2L}^{\mu} <0$ and $a_2+r_{2U}^{\prime\mu} \geq 0$.\\
Then, $~~~~~~~~~~~~~~~~~~~~~~~~~a_2+r_{2L}^{\mu}R^{-1}(\alpha) \leq~(\geq)~ 0~$ for $~ \alpha \leq~(\geq)~ R\left(-\dfrac{a_2}{r_{2L}^{\mu}}\right)$.\\

\noindent \underline{Sub-case 6.1.} $a_1-l^{\nu}_{1L}L^{-1}(1-\beta) \leq 0$.\\\\
Then,\\
\noindent $~~~~~~~~~~~~~\big(\tilde{A}_1\odot \tilde{A}_2\big)_{\alpha}^{L}=A_{1\alpha}^{L} \times A_{2\alpha}^{L}=[(a_1+r^{\mu}_{1L}R^{-1}(\alpha))(a_2-l^{\mu}_{2L}L^{-1}(\alpha)),~(a_1-l^{\mu}_{1L}L^{-1}(\alpha))(a_2+r^{\mu}_{2L}R^{-1}(\alpha))].$\\

\noindent Choosing $\alpha=\alpha_0$, we obtain
$$\big(\tilde{A}_1\odot \tilde{A}_2\big)_{\alpha_0}^{L}=[(a_1+r^{\mu}_{1L})(a_2-l^{\mu}_{2L}),(a_1-l^{\mu}_{1L})(a_2+r^{\mu}_{2L})]$$
which gives $~~~~~~~~~~~~~~~~~~~~~~~~~~~~~~~~~~~~r_L^{\mu}=(a_1-l^{\mu}_{1L})(a_2+r^{\mu}_{2L})-a_1a_2.$\\\\
Further, we have 
$$r_U^{\prime\mu}=(a_1+r_{1U}^{\prime\mu})(a_2+r_{2U}^{\prime\mu})-a_1a_2~~ \mbox{and} ~~r_L^{\nu}= \mbox{max}\{\left(a_1-l^{\nu}_{1L}\right)
\left(a_2-l^{\nu}_{2L}\right),\left(a_1+r^{\nu}_{1L}\right)\left(a_2+r^{\nu}_{2L} \right)\}-a_1a_2.$$
The proof that $r_U^{\prime\mu} \geq r_L^{\mu}$ and $r_L^{\nu} \geq r_L^{\mu}$ can be obtained on the lines of Sub-cases $2.1$ and $3.1$.\\

\noindent \underline{Sub-case 6.2.} $~~a_1-l^{\nu}_{1L}L^{-1}(1-\beta) \geq 0$.\\
This part can be proved on the lines of Sub-case $1.1$.\\

\noindent \textbf{Case 7.} $a_2+r_{2U}^{\prime\mu} <0$ and $ a_2+r_{2U}^{\prime\nu} \geq 0$.\\
Then, $~~~~~~~~~~~~~~~~~~~~~~~~~~~~a_2+r_{2U}^{\prime\mu}(R')^{-1}(\alpha) \leq~(\geq)~ 0~$ for $~ \alpha \leq~(\geq)~ R'\left(-\dfrac{a_2}{r_{2U}^{\prime\mu}}\right)$.\\

\noindent \underline{Sub-case 7.1.} $~~a_1-l^{\nu}_{1L}L^{-1}(1-\beta) \leq 0$.\\\\
Then,\\
\noindent $~~~~\big(\tilde{A}_1\odot \tilde{A}_2\big)_{\alpha}^{U}=A_{1\alpha}^{U} \times A_{2\alpha}^{U}=[(a_1+r^{\prime\mu}_{1U}(R')^{-1}(\alpha))(a_2-l^{\prime\mu}_{2U}(L')^{-1}(\alpha)),~(a_1-l^{\prime\mu}_{1U}(L')^{-1}(\alpha))(a_2+r^{\prime\mu}_{2U}(R')^{-1}(\alpha))].$\\

\noindent Taking $\alpha=\alpha_0$, we obtain
$$\big(\tilde{A}_1\odot \tilde{A}_2\big)_{\alpha_0}^{U}=[(a_1+r^{\prime\mu}_{1U})(a_2-l^{\prime\mu}_{2U}),(a_1-l^{\prime\mu}_{1U})(a_2+r^{\prime\mu}_{2U})].$$
It follows that $~~~~~~~~~~~~~~~~~~~~~~~~~~~~~~~r_U^{\prime\mu}=(a_1-l^{\prime\mu}_{1U})(a_2+r^{\prime\mu}_{2U})-a_1a_2.$\\\\
Also, we have \\\\
$~~~~~~~~~~~~r_L^{\mu}=\mbox{max} \{(a_1-l_{1L}^{\mu})(a_2+r_{2L}^{\mu}),~(a_1+r_{1L}^{\mu})(a_2+r_{2L}^{\mu}) \}-a_1a_2~~\mbox{and}~~r_U^{\prime\nu}=(a_1+r_{1U}^{\prime\nu})(a_2+r_{2U}^{\prime\nu})-a_1a_2.$\\\\
Further, the inequalities $r_U^{\prime\mu} \geq r_L^{\mu}$ and $r_U^{\prime\nu} \geq r_U^{\prime\mu}$ can be established on the lines of Sub-cases $2.1$ and $3.1$.\\

\noindent \underline{Sub-case 7.2.} $~~a_1-l^{\nu}_{1L}L^{-1}(1-\beta) \geq 0$.\\
This part can be proved on the similar lines as Sub-case $1.1$.\\

\noindent \textbf{Case 8.} $a_2+r_{2U}^{\prime\nu} <0$ and $ a_2+r_{2L}^{\nu} \geq 0$.\\\\
Then, $~~~~~~~~~~~~~~~~~~~~a_2+r_{2U}^{\prime\nu}(R')^{-1}(1-\beta) \leq~(\geq)~ 0~$ for $~ \beta \geq~(\leq)~ 1-R'\left(-\dfrac{a_2}{r_{2U}^{\prime\nu}}\right)$.\\

\noindent \underline{Sub-case 8.1.} $~~a_1-l^{\nu}_{1L}L^{-1}(1-\beta) \leq 0$.\\\\
Then, \\
\noindent $\big(\tilde{A}_1\odot \tilde{A}_2\big)_{\beta}^{U}=A_{1\beta}^{U} \times A_{2\beta}^{U}$\\
$~~~~~~~~~~~~~~~~~~~~~~=[(a_1+r^{\prime\nu}_{1U}(R')^{-1}(1-\beta))(a_2-l^{\prime\nu}_{2U}(L')^{-1}(1-\beta)),~(a_1-l^{\prime\nu}_{1U}(L')^{-1}(1-\beta))(a_2+r^{\prime\nu}_{2U}(R')^{-1}(1-\beta))].$\\

\noindent Choosing $\beta=\beta_0$, we get
$$\big(\tilde{A}_1\odot \tilde{A}_2\big)_{\beta_0}^{U}=[(a_1+r^{\prime\nu}_{1U})(a_2-l^{\prime\nu}_{2U}),(a_1-l^{\prime\nu}_{1U})(a_2+r^{\prime\nu}_{2U})].$$
This further yields $~~~~~~~~~~~~~~~~~~~~~r_U^{\prime\nu}=(a_1-l^{\prime\nu}_{1U})(a_2+r^{\prime\nu}_{2U})-a_1a_2.$\\\\
Moreover, we have\\\\
$~~~~~~~~~~~~~~~~~~~~~~~~~~~~~~r_L^{\nu}= \mbox{max}\left\{\left(a_1-l^{\nu}_{1L}\right)
\left(a_2-l^{\nu}_{2L}\right),\left(a_1+r^{\nu}_{1L}\right)\left(a_2+r^{\nu}_{2L} \right)\right\}-a_1a_2~~~$ and\\\\
$~~~~~~~~~~~~~~~~~~~~~~~~~~~~~~r_U^{\prime\mu}=\mbox{max} \left\{(a_1-l_{1U}^{\prime\mu})(a_2+r_{2U}^{\prime\mu}),~(a_1+r_{1U}^{\prime\mu})(a_2+r_{2U}^{\prime\mu}) \right\}-a_1a_2.$\\\\
Next, $r_L^{\nu} \geq r_U^{\prime\nu}$ and $r_U^{\prime\nu} \geq r_U^{\prime\mu}$ can be proved following the lines of Sub-cases $2.1$ and $3.1$.\\

\noindent \underline{Sub-case 8.2.} $~~a_1-l^{\nu}_{1L}L^{-1}(1-\beta) \geq 0$.\\
The proof of this part follows on the lines of Sub-case $1.1$.\\

\noindent \textbf{Case 9.} $ a_2+r_{2L}^{\nu} < 0$.\\\\
Then, $~~~~~~~~~~~~~~~~~~~a_2+r_{2L}^{\nu}R^{-1}(1-\beta) \leq~(\geq)~ 0~$ for $~ \beta \geq~(\leq)~ 1-R\left(-\dfrac{a_2}{r_{2L}^{\nu}}\right)$.\\
 
\noindent \underline{Sub-case 9.1.} $~~a_1-l^{\nu}_{1L}L^{-1}(1-\beta) \leq 0$.\\\\
Then, \\
$\big(\tilde{A}_1\odot \tilde{A}_2\big)_{\beta}^{L}=A_{1\beta}^{L} \times A_{2\beta}^{L}$\\
$~~~~~~~~~~~~~~~~~~~~~=[(a_1+r^{\nu}_{1L}R^{-1}(1-\beta))(a_2-l^{\nu}_{2L}L^{-1}(1-\beta)),~(a_1-l^{\nu}_{1L}L^{-1}(1-\beta))(a_2-l^{\nu}_{2L}L^{-1}(1-\beta))].$\\

\noindent Taking $\beta=\beta_0$, we obtain
$$\big(\tilde{A}_1\odot \tilde{A}_2\big)_{\beta_0}^{L}=[(a_1+r^{\nu}_{1L})(a_2-l^{\nu}_{2L}),(a_1-l^{\nu}_{1L})(a_2-l^{\nu}_{2L})].$$
It follows that $~~~~~~~~~~~~~~~~~~~~~~~~~~~~~~~~r_L^{\nu}=(a_1-l^{\nu}_{1L})(a_2-l^{\nu}_{2L})-a_1a_2.$\\
Further, \\\\
$~~~~~~~~~~~~~~~~~~~~~~~~r_U^{\prime\nu}=\mbox{max} \left\{(a_1+r_{1U}^{\prime\nu})(a_2+r_{2U}^{\prime\nu}),~(a_1-l_{1U}^{\prime\nu})(a_2+r_{2U}^{\prime\nu}) \right\}-a_1a_2~~$ and\\\\
$~~~~~~~~~~~~~~~~~~~~~~~~r_L^{\mu}=\mbox{max} \left\{(a_1-l_{1L}^{\mu})(a_2+r_{2L}^{\mu}),~(a_1+r_{1L}^{\mu})(a_2+r_{2L}^{\mu}) \right\}-a_1a_2.$\\\\
Finally, the proof of the inequalities $r_L^{\nu} \geq r_U^{\prime\nu}$ and $r_L^{\nu} \geq r_L^{\mu}$ can be obtained on the lines of Sub-cases $2.1$ and $3.1$.\\

\noindent \underline{Sub-case 9.2.} $~~a_1-l^{\nu}_{1L}L^{-1}(1-\beta) \geq 0$.\\
The proof of this part can be established on the lines of Sub-case $1.1$.\\

\noindent Finally, combining all the Cases $(1)-(9)$, we have\\
\noindent $\tilde{A}_1 \odot \tilde{A}_2=(a; l_{L}^{\mu},r_{L}^{\mu},l'^{\mu}_{U},r'^{\mu}_{U};l^{\nu}_{L},r^{\nu}_{L},l'^{\nu}_{U},r'^{\nu}_{U})_{LR}$ where\\
\noindent $a=a_1a_2,$\\
\noindent $l_L^{\mu}=a_1a_2-\mbox{min} \{(a_1-l_{1L}^{\mu})(a_2-l_{2L}^{\mu}),~(a_1+r_{1L}^{\mu})(a_2-l_{2L}^{\mu}) \},$\\
\noindent $r_L^{\mu}=\mbox{max} \{(a_1-l_{1L}^{\mu})(a_2+r_{2L}^{\mu}),~(a_1+r_{1L}^{\mu})(a_2+r_{2L}^{\mu}) \}-a_1a_2,$\\
\noindent $l_U^{\prime\mu}=a_1a_2-\mbox{min} \{(a_1-l_{1U}^{\prime\mu})(a_2-l_{2U}^{\prime\mu}),~(a_1+r_{1U}^{\prime\mu})(a_2-l_{2U}^{\prime\mu}) \},$\\
\noindent $r_U^{\prime\mu}=\mbox{max} \{(a_1-l_{1U}^{\prime\mu})(a_2+r_{2U}^{\prime\mu}),~(a_1+r_{1U}^{\prime\mu})(a_2+r_{2U}^{\prime\mu}) \}-a_1a_2,$\\
\noindent $l_L^{\nu}=a_1a_2-\mbox{min} \{(a_1-l_{1L}^{\nu})(a_2+r_{2L}^{\nu}),~(a_1+r_{1L}^{\nu})(a_2-l_{2L}^{\nu}) \},$\\
\noindent $r_L^{\nu}=\mbox{max} \{(a_1-l_{1L}^{\nu})(a_2-l_{2L}^{\nu}),~(a_1+r_{1L}^{\nu})(a_2+r_{2L}^{\nu}) \}-a_1a_2,$\\ 
\noindent $l_U^{\prime\nu}=a_1a_2-\mbox{min} \{(a_1-l_{1U}^{\prime\nu})(a_2-l_{2U}^{\prime\nu}),~(a_1+r_{1U}^{\prime\nu})(a_2-l_{2U}^{\prime\nu}) \},$\\
\noindent $r_U^{\prime\nu}=\mbox{max} \{(a_1+r_{1U}^{\prime\nu})(a_2+r_{2U}^{\prime\nu}),~(a_1-l_{1U}^{\prime\nu})(a_2+r_{2U}^{\prime\nu}) \}-a_1a_2$\\ 

\noindent where $l'^{\mu}_{U} \geq l^{\mu}_{L} >0,~~r'^{\mu}_{U} \geq r^{\mu}_{L} >0,~~l^{\nu}_{L} \geq l'^{\nu}_{U} >0,~~r^{\nu}_{L} \geq r'^{\nu}_{U} >0,~~l^{\nu}_{L} \geq l^{\mu}_{L} >0,~~r^{\nu}_{L} \geq r^{\mu}_{L} >0,~~l'^{\nu}_{U} \geq l'^{\mu}_{U} >0,$\\
$~~~~~~~~~~~~r'^{\nu}_{U} \geq r'^{\mu}_{U} >0.$\\\\
Hence, the result.\\

\noindent{\bf{Proposition 3.1.4.}} \textit{Let $\tilde{A}_1=(a_1; l_{1L}^{\mu},r_{1L}^{\mu},l'^{\mu}_{1U},r'^{\mu}_{1U};l^{\nu}_{1L},r^{\nu}_{1L},l'^{\nu}_{1U},r'^{\nu}_{1U})_{LR}$ be an $LR$-type IVIFN such that $a_1-l_{1U}^{\prime\nu} < 0,$\\
$a_1-l_{1U}^{\prime\mu} \geq 0$ and $\tilde{A}_2=(a_2; l_{2L}^{\mu},r_{2L}^{\mu},l'^{\mu}_{2U},r'^{\mu}_{2U}; l^{\nu}_{2L},r^{\nu}_{2L},l'^{\nu}_{2U},r'^{\nu}_{2U})_{LR}$ be any $LR$-type IVIFN. Then\\
  $\tilde{A}_1 \odot \tilde{A}_2=(a; l_{L}^{\mu},r_{L}^{\mu},l'^{\mu}_{U},r'^{\mu}_{U};l^{\nu}_{L},r^{\nu}_{L},l'^{\nu}_{U},r'^{\nu}_{U})_{LR}$ where\\ 
\noindent $a=a_1a_2,$\\
\noindent $l_L^{\mu}=a_1a_2-\mbox{min} \{(a_1+r_{1L}^{\mu})(a_2-l_{2L}^{\mu}),~(a_1-l_{1L}^{\mu})(a_2-l_{2L}^{\mu}) \},$\\
\noindent $r_L^{\mu}=\mbox{max} \{(a_1+r_{1L}^{\mu})(a_2+r_{2L}^{\mu}),~(a_1-l_{1L}^{\mu})(a_2+r_{2L}^{\mu}) \}-a_1a_2,$\\
\noindent $l_U^{\prime\mu}=a_1a_2-\mbox{min} \{(a_1+r_{1U}^{\prime\mu})(a_2-l_{2U}^{\prime\mu}),~(a_1-l_{1U}^{\prime\mu})(a_2-l_{2U}^{\prime\mu}) \},$\\
\noindent $r_U^{\prime\mu}=\mbox{max} \{(a_1+r_{1U}^{\prime\mu})(a_2+r_{2U}^{\prime\mu}),~(a_1-l_{1U}^{\prime\mu})(a_2+r_{2U}^{\prime\mu}) \}-a_1a_2,$\\ 
\noindent $l_L^{\nu}=a_1a_2-\mbox{min} \{(a_1-l_{1L}^{\nu})(a_2+r_{2L}^{\nu}),~(a_1+r_{1L}^{\nu})(a_2-l_{2L}^{\nu}) \},$\\
\noindent $r_L^{\nu}=\mbox{max} \{(a_1+r_{1L}^{\nu})(a_2+r_{2L}^{\nu}),~(a_1-l_{1L}^{\nu})(a_2-l_{2L}^{\nu}) \}-a_1a_2,$\\ 
\noindent $l_U^{\prime\nu}=a_1a_2-\mbox{min} \{(a_1-l_{1U}^{\prime\nu})(a_2+r_{2U}^{\prime\nu}),~(a_1+r_{1U}^{\prime\nu})(a_2-l_{2U}^{\prime\nu}) \},$\\
\noindent $r_U^{\prime\nu}=\mbox{max} \{(a_1+r_{1U}^{\prime\nu})(a_2+r_{2U}^{\prime\nu}),~(a_1-l_{1U}^{\prime\nu})(a_2-l_{2U}^{\prime\nu}) \}-a_1a_2,$\\ 
\noindent where the conditions for $LR-$type representation of $\tilde{A}_1\odot \tilde{A}_2$ are satisfied.}\\

\noindent{\bf{Proof.}} Similar to the Proposition $3.1.3$.\\

\noindent{\bf{Proposition 3.1.5.}} \textit{Let $\tilde{A}_1=(a_1; l_{1L}^{\mu},r_{1L}^{\mu},l'^{\mu}_{1U},r'^{\mu}_{1U};l^{\nu}_{1L},r^{\nu}_{1L},l'^{\nu}_{1U},r'^{\nu}_{1U})_{LR}$ be an $LR$-type IVIFN such that $a_1-l_{1U}^{\prime\mu} < 0,$\\
$ a_1-l_{1L}^{\mu} \geq 0$ and $\tilde{A}_2=(a_2; l_{2L}^{\mu},r_{2L}^{\mu},l'^{\mu}_{2U},r'^{\mu}_{2U}; l^{\nu}_{2L},r^{\nu}_{2L},l'^{\nu}_{2U},r'^{\nu}_{2U})_{LR}$ be any $LR$-type IVIFN. Then\\
  $\tilde{A}_1 \odot \tilde{A}_2=(a; l_{L}^{\mu},r_{L}^{\mu},l'^{\mu}_{U},r'^{\mu}_{U};l^{\nu}_{L},r^{\nu}_{L},l'^{\nu}_{U},r'^{\nu}_{U})_{LR}$ where\\ 
\noindent $a=a_1a_2,$\\
\noindent $l_L^{\mu}=a_1a_2-\mbox{min} \{(a_1-l_{1L}^{\mu})(a_2-l_{2L}^{\mu}),~(a_1+r_{1L}^{\mu})(a_2-l_{2L}^{\mu}) \},$\\
\noindent $r_L^{\mu}=\mbox{max} \{(a_1+r_{1L}^{\mu})(a_2+r_{2L}^{\mu}),~(a_1-l_{1L}^{\mu})(a_2+r_{2L}^{\mu}) \}-a_1a_2,$\\ 
\noindent $l_U^{\prime\mu}=a_1a_2-\mbox{min} \{(a_1-l_{1U}^{\prime\mu})(a_2+r_{2U}^{\prime\mu}),~(a_1+r_{1U}^{\prime\mu})(a_2-l_{2U}^{\prime\mu}) \},$\\
\noindent $r_U^{\prime\mu}=\mbox{max} \{(a_1+r_{1U}^{\prime\mu})(a_2+r_{2U}^{\prime\mu}),~(a_1-l_{1U}^{\prime\mu})(a_2-l_{2U}^{\prime\mu}) \}-a_1a_2,$\\ 
\noindent $l_L^{\nu}=a_1a_2-\mbox{min} \{(a_1-l_{1L}^{\nu})(a_2+r_{2L}^{\nu}),~(a_1+r_{1L}^{\nu})(a_2-l_{2L}^{\nu}) \},$\\
\noindent $r_L^{\nu}=\mbox{max} \{(a_1+r_{1L}^{\nu})(a_2+r_{2L}^{\nu}),~(a_1-l_{1L}^{\nu})(a_2-l_{2L}^{\nu}) \}-a_1a_2,$\\ 
\noindent $l_U^{\prime\nu}=a_1a_2-\mbox{min} \{(a_1-l_{1U}^{\prime\nu})(a_2+r_{2U}^{\prime\nu}),~(a_1+r_{1U}^{\prime\nu})(a_2-l_{2U}^{\prime\nu}) \},$\\
\noindent $r_U^{\prime\nu}=\mbox{max} \{(a_1+r_{1U}^{\prime\nu})(a_2+r_{2U}^{\prime\nu}),~(a_1-l_{1U}^{\prime\nu})(a_2-l_{2U}^{\prime\nu}) \}-a_1a_2,$\\ 
\noindent where the conditions for $LR-$type representation of $\tilde{A}_1\odot \tilde{A}_2$ are satisfied.}\\

\noindent{\bf{Proof.}} Similar to the Proposition $3.1.3$. \\

\noindent{\bf{Proposition 3.1.6.}} \textit{Let $\tilde{A}_1=(a_1; l_{1L}^{\mu},r_{1L}^{\mu},l'^{\mu}_{1U},r'^{\mu}_{1U};l^{\nu}_{1L},r^{\nu}_{1L},l'^{\nu}_{1U},r'^{\nu}_{1U})_{LR}$ be an $LR$-type IVIFN such that $a_1-l_{1L}^{\mu} < 0,$\\
$a_1 \geq 0$ and $\tilde{A}_2=(a_2; l_{2L}^{\mu},r_{2L}^{\mu},l'^{\mu}_{2U},r'^{\mu}_{2U};l^{\nu}_{2L},r^{\nu}_{2L},l'^{\nu}_{2U},r'^{\nu}_{2U})_{LR}$ be any $LR$-type IVIFN. Then\\
  $\tilde{A}_1 \odot \tilde{A}_2=(a; l_{L}^{\mu},r_{L}^{\mu},l'^{\mu}_{U},r'^{\mu}_{U};l^{\nu}_{L},r^{\nu}_{L},l'^{\nu}_{U},r'^{\nu}_{U})_{LR}$ where\\  
\noindent $a=a_1a_2,$\\
\noindent $l_L^{\mu}=a_1a_2-\mbox{min} \{(a_1-l_{1L}^{\mu})(a_2+r_{2L}^{\mu}),~(a_1+r_{1L}^{\mu})(a_2-l_{2L}^{\mu}) \},$\\
\noindent $r_L^{\mu}=\mbox{max} \{(a_1+r_{1L}^{\mu})(a_2+r_{2L}^{\mu}),~(a_1-l_{1L}^{\mu})(a_2-l_{2L}^{\mu}) \}-a_1a_2,$\\ 
\noindent $l_U^{\prime\mu}=a_1a_2-\mbox{min} \{(a_1-l_{1U}^{\prime\mu})(a_2+r_{2U}^{\prime\mu}),~(a_1+r_{1U}^{\prime\mu})(a_2-l_{2U}^{\prime\mu}) \},$\\
\noindent $r_U^{\prime\mu}=\mbox{max} \{(a_1+r_{1U}^{\prime\mu})(a_2+r_{2U}^{\prime\mu}),~(a_1-l_{1U}^{\prime\mu})(a_2-l_{2U}^{\prime\mu}) \}-a_1a_2,$\\ 
\noindent $l_L^{\nu}=a_1a_2-\mbox{min} \{(a_1-l_{1L}^{\nu})(a_2+r_{2L}^{\nu}),~(a_1+r_{1L}^{\nu})(a_2-l_{2L}^{\nu}) \},$\\
\noindent $r_L^{\nu}=\mbox{max} \{(a_1+r_{1L}^{\nu})(a_2+r_{2L}^{\nu}),~(a_1-l_{1L}^{\nu})(a_2-l_{2L}^{\nu}) \}-a_1a_2,$\\ 
\noindent $l_U^{\prime\nu}=a_1a_2-\mbox{min} \{(a_1-l_{1U}^{\prime\nu})(a_2+r_{2U}^{\prime\nu}),~(a_1+r_{1U}^{\prime\nu})(a_2-l_{2U}^{\prime\nu}) \},$\\
\noindent $r_U^{\prime\nu}=\mbox{max} \{(a_1+r_{1U}^{\prime\nu})(a_2+r_{2U}^{\prime\nu}),~(a_1-l_{1U}^{\prime\nu})(a_2-l_{2U}^{\prime\nu}) \}-a_1a_2,$\\ 
\noindent where the conditions for $LR-$type representation of $\tilde{A}_1\odot \tilde{A}_2$ are satisfied.}\\

\noindent{\bf{Proof.}} Similar to the Proposition $3.1.3$.\\

\noindent{\bf{Proposition 3.1.7.}} \textit{Let $\tilde{A}_1=(a_1; l_{1L}^{\mu},r_{1L}^{\mu},l'^{\mu}_{1U},r'^{\mu}_{1U};l^{\nu}_{1L},r^{\nu}_{1L},l'^{\nu}_{1U},r'^{\nu}_{1U})_{LR}$ be an $LR$-type IVIFN such that $a_1 < 0,$\\
$a_1+r_{1L}^{\mu} \geq 0$ and $\tilde{A}_2=(a_2; l_{2L}^{\mu},r_{2L}^{\mu},l'^{\mu}_{2U},r'^{\mu}_{2U}; l^{\nu}_{2L},r^{\nu}_{2L},l'^{\nu}_{2U},r'^{\nu}_{2U})_{LR}$ be any $LR$-type IVIFN. Then\\
  $\tilde{A}_1 \odot \tilde{A}_2=(a; l_{L}^{\mu},r_{L}^{\mu},l'^{\mu}_{U},r'^{\mu}_{U};l^{\nu}_{L},r^{\nu}_{L},l'^{\nu}_{U},r'^{\nu}_{U})_{LR}$ where\\  
\noindent $a=a_1a_2,$\\
\noindent $l_L^{\mu}=a_1a_2-\mbox{min} \{(a_1-l_{1L}^{\mu})(a_2+r_{2L}^{\mu}),~(a_1+r_{1L}^{\mu})(a_2-l_{2L}^{\mu}) \},$\\
\noindent $r_L^{\mu}=\mbox{max} \{(a_1+r_{1L}^{\mu})(a_2+r_{2L}^{\mu}),~(a_1-l_{1L}^{\mu})(a_2-l_{2L}^{\mu}) \}-a_1a_2,$\\
\noindent $l_U^{\prime\mu}=a_1a_2-\mbox{min} \{(a_1-l_{1U}^{\prime\mu})(a_2+r_{2U}^{\prime\mu}),~(a_1+r_{1U}^{\prime\mu})(a_2-l_{2U}^{\prime\mu}) \},$\\
\noindent $r_U^{\prime\mu}=\mbox{max} \{(a_1+r_{1U}^{\prime\mu})(a_2+r_{2U}^{\prime\mu}),~(a_1-l_{1U}^{\prime\mu})(a_2-l_{2U}^{\prime\mu}) \}-a_1a_2,$\\ 
\noindent $l_L^{\nu}=a_1a_2-\mbox{min} \{(a_1-l_{1L}^{\nu})(a_2+r_{2L}^{\nu}),~(a_1+r_{1L}^{\nu})(a_2-l_{2L}^{\nu}) \},$\\
\noindent $r_L^{\nu}=\mbox{max} \{(a_1+r_{1L}^{\nu})(a_2+r_{2L}^{\nu}),~(a_1-l_{1L}^{\nu})(a_2-l_{2L}^{\nu}) \}-a_1a_2,$\\ 
\noindent $l_U^{\prime\nu}=a_1a_2-\mbox{min} \{(a_1-l_{1U}^{\prime\nu})(a_2+r_{2U}^{\prime\nu}),~(a_1+r_{1U}^{\prime\nu})(a_2-l_{2U}^{\prime\nu}) \},$\\
\noindent $r_U^{\prime\nu}=\mbox{max} \{(a_1+r_{1U}^{\prime\nu})(a_2+r_{2U}^{\prime\nu}),~(a_1-l_{1U}^{\prime\nu})(a_2-l_{2U}^{\prime\nu}) \}-a_1a_2,$\\ 
\noindent where the conditions for $LR-$type representation of $\tilde{A}_1\odot \tilde{A}_2$ are satisfied.}\\

\noindent{\bf{Proof.}} Similar to the Proposition $3.1.3$.\\

\noindent{\bf{Proposition 3.1.8.}} \textit{Let $\tilde{A}_1=(a_1; l_{1L}^{\mu},r_{1L}^{\mu},l'^{\mu}_{1U},r'^{\mu}_{1U};l^{\nu}_{1L},r^{\nu}_{1L},l'^{\nu}_{1U},r'^{\nu}_{1U})_{LR}$ be an $LR$-type IVIFN such that $a_1+r_{1L}^{\mu} < 0,$\\
$a_1+r_{1U}^{\prime\mu} \geq 0$ and $\tilde{A}_2=(a_2; l_{2L}^{\mu},r_{2L}^{\mu},l'^{\mu}_{2U},r'^{\mu}_{2U}; l^{\nu}_{2L},r^{\nu}_{2L},l'^{\nu}_{2U},r'^{\nu}_{2U})_{LR}$ be any $LR$-type IVIFN. Then\\
  $\tilde{A}_1 \odot \tilde{A}_2=(a; l_{L}^{\mu},r_{L}^{\mu},l'^{\mu}_{U},r'^{\mu}_{U};l^{\nu}_{L},r^{\nu}_{L},l'^{\nu}_{U},r'^{\nu}_{U})_{LR}$ where\\  
\noindent $a=a_1a_2,$\\
\noindent $l_L^{\mu}=a_1a_2-\mbox{min} \{(a_1+r_{1L}^{\mu})(a_2+r_{2L}^{\mu}),~(a_1-l_{1L}^{\mu})(a_2+r_{2L}^{\mu}) \},$\\
\noindent $r_L^{\mu}=\mbox{max} \{(a_1-l_{1L}^{\mu})(a_2-l_{2L}^{\mu}),~(a_1+r_{1L}^{\mu})(a_2-l_{2L}^{\mu}) \}-a_1a_2,$\\ 
\noindent $l_U^{\prime\mu}=a_1a_2-\mbox{min} \{(a_1-l_{1U}^{\prime\mu})(a_2+r_{2U}^{\prime\mu}),~(a_1+r_{1U}^{\prime\mu})(a_2-l_{2U}^{\prime\mu}) \},$\\
\noindent $r_U^{\prime\mu}=\mbox{max} \{(a_1+r_{1U}^{\prime\mu})(a_2+r_{2U}^{\prime\mu}),~(a_1-l_{1U}^{\prime\mu})(a_2-l_{2U}^{\prime\mu}) \}-a_1a_2,$\\ 
\noindent $l_L^{\nu}=a_1a_2-\mbox{min} \{(a_1-l_{1L}^{\nu})(a_2+r_{2L}^{\nu}),~(a_1+r_{1L}^{\nu})(a_2-l_{2L}^{\nu}) \},$\\
\noindent $r_L^{\nu}=\mbox{max} \{(a_1+r_{1L}^{\nu})(a_2+r_{2L}^{\nu}),~(a_1-l_{1L}^{\nu})(a_2-l_{2L}^{\nu}) \}-a_1a_2,$\\ 
\noindent $l_U^{\prime\nu}=a_1a_2-\mbox{min} \{(a_1-l_{1U}^{\prime\nu})(a_2+r_{2U}^{\prime\nu}),~(a_1+r_{1U}^{\prime\nu})(a_2-l_{2U}^{\prime\nu}) \},$\\
\noindent $r_U^{\prime\nu}=\mbox{max} \{(a_1+r_{1U}^{\prime\nu})(a_2+r_{2U}^{\prime\nu}),~(a_1-l_{1U}^{\prime\nu})(a_2-l_{2U}^{\prime\nu}) \}-a_1a_2,$\\ 
\noindent where the conditions for $LR-$type representation of $\tilde{A}_1\odot \tilde{A}_2$ are satisfied.}\\

\noindent{\bf{Proof.}} Similar to the Proposition $3.1.3$.\\

\noindent{\bf{Proposition 3.1.9.}} \textit{Let $\tilde{A}_1=(a_1; l_{1L}^{\mu},r_{1L}^{\mu},l'^{\mu}_{1U},r'^{\mu}_{1U};l^{\nu}_{1L},r^{\nu}_{1L},l'^{\nu}_{1U},r'^{\nu}_{1U})_{LR}$ be an $LR$-type IVIFN such that $a_1+r_{1U}^{\prime\mu} < 0,$\\
$a_1+r_{1U}^{\prime\nu} \geq 0$ and $\tilde{A}_2=(a_2; l_{2L}^{\mu},r_{2L}^{\mu},l'^{\mu}_{2U},r'^{\mu}_{2U}; l^{\nu}_{2L},r^{\nu}_{2L},l'^{\nu}_{2U},r'^{\nu}_{2U})_{LR}$ be any $LR$-type IVIFN. Then\\
  $\tilde{A}_1 \odot \tilde{A}_2=(a; l_{L}^{\mu},r_{L}^{\mu},l'^{\mu}_{U},r'^{\mu}_{U};l^{\nu}_{L},r^{\nu}_{L},l'^{\nu}_{U},r'^{\nu}_{U})_{LR}$ where\\  
\noindent $a=a_1a_2,$\\
\noindent $l_L^{\mu}=a_1a_2-\mbox{min} \{(a_1+r_{1L}^{\mu})(a_2+r_{2L}^{\mu}),~(a_1-l_{1L}^{\mu})(a_2+r_{2L}^{\mu}) \},$\\
\noindent $r_L^{\mu}=\mbox{max} \{(a_1-l_{1L}^{\mu})(a_2-l_{2L}^{\mu}),~(a_1+r_{1L}^{\mu})(a_2-l_{2L}^{\mu}) \}-a_1a_2,$\\ 
\noindent $l_U^{\prime\mu}=a_1a_2-\mbox{min} \{(a_1+r_{1U}^{\prime\mu})(a_2+r_{2U}^{\prime\mu}),~(a_1-l_{1U}^{\prime\mu})(a_2+r_{2U}^{\prime\mu}) \},$\\
\noindent $r_U^{\prime\mu}=\mbox{max} \{(a_1-l_{1U}^{\prime\mu})(a_2-l_{2U}^{\prime\mu}),~(a_1+r_{1U}^{\prime\mu})(a_2-l_{2U}^{\prime\mu}) \}-a_1a_2,$\\ 
\noindent $l_L^{\nu}=a_1a_2-\mbox{min} \{(a_1-l_{1L}^{\nu})(a_2+r_{2L}^{\nu}),~(a_1+r_{1L}^{\nu})(a_2-l_{2L}^{\nu}) \},$\\
\noindent $r_L^{\nu}=\mbox{max} \{(a_1+r_{1L}^{\nu})(a_2+r_{2L}^{\nu}),~(a_1-l_{1L}^{\nu})(a_2-l_{2L}^{\nu}) \}-a_1a_2,$\\ 
\noindent $l_U^{\prime\nu}=a_1a_2-\mbox{min} \{(a_1-l_{1U}^{\prime\nu})(a_2+r_{2U}^{\prime\nu}),~(a_1+r_{1U}^{\prime\nu})(a_2-l_{2U}^{\prime\nu}) \},$\\
\noindent $r_U^{\prime\nu}=\mbox{max} \{(a_1+r_{1U}^{\prime\nu})(a_2+r_{2U}^{\prime\nu}),~(a_1-l_{1U}^{\prime\nu})(a_2-l_{2U}^{\prime\nu}) \}-a_1a_2,$\\ 
\noindent where the conditions for $LR-$type representation of $\tilde{A}_1\odot \tilde{A}_2$ are satisfied.}\\

\noindent{\bf{Proof.}} Similar to the Proposition $3.1.3$. \\

\noindent{\bf{Proposition 3.1.10.}} \textit{Let $\tilde{A}_1=(a_1; l_{1L}^{\mu},r_{1L}^{\mu},l'^{\mu}_{1U},r'^{\mu}_{1U};l^{\nu}_{1L},r^{\nu}_{1L},l'^{\nu}_{1U},r'^{\nu}_{1U})_{LR}$ be an $LR$-type IVIFN such that $a_1+r_{1U}^{\prime\nu} < 0,$\\
$a_1+r_{1L}^{\nu} \geq 0$ and $\tilde{A}_2=(a_2; l_{2L}^{\mu},r_{2L}^{\mu},l'^{\mu}_{2U},r'^{\mu}_{2U}; l^{\nu}_{2L},r^{\nu}_{2L},l'^{\nu}_{2U},r'^{\nu}_{2U})_{LR}$ be any $LR$-type IVIFN. Then\\
  $\tilde{A}_1 \odot \tilde{A}_2=(a; l_{L}^{\mu},r_{L}^{\mu},l'^{\mu}_{U},r'^{\mu}_{U};l^{\nu}_{L},r^{\nu}_{L},l'^{\nu}_{U},r'^{\nu}_{U})_{LR}$ where\\  
\noindent $a=a_1a_2,$\\
\noindent $l_L^{\mu}=a_1a_2-\mbox{min} \{(a_1+r_{1L}^{\mu})(a_2+r_{2L}^{\mu}),~(a_1-l_{1L}^{\mu})(a_2+r_{2L}^{\mu}) \},$\\
\noindent $r_L^{\mu}=\mbox{max} \{(a_1-l_{1L}^{\mu})(a_2-l_{2L}^{\mu}),~(a_1+r_{1L}^{\mu})(a_2-l_{2L}^{\mu}) \}-a_1a_2,$\\ 
\noindent $l_U^{\prime\mu}=a_1a_2-\mbox{min} \{(a_1+r_{1U}^{\prime\mu})(a_2+r_{2U}^{\prime\mu}),~(a_1-l_{1U}^{\prime\mu})(a_2+r_{2U}^{\prime\mu}) \},$\\
\noindent $r_U^{\prime\mu}=\mbox{max} \{(a_1-l_{1U}^{\prime\mu})(a_2-l_{2U}^{\prime\mu}),~(a_1+r_{1U}^{\prime\mu})(a_2-l_{2U}^{\prime\mu}) \}-a_1a_2,$\\ 
\noindent $l_L^{\nu}=a_1a_2-\mbox{min} \{(a_1-l_{1L}^{\nu})(a_2+r_{2L}^{\nu}),~(a_1+r_{1L}^{\nu})(a_2-l_{2L}^{\nu}) \},$\\
\noindent $r_L^{\nu}=\mbox{max} \{(a_1+r_{1L}^{\nu})(a_2+r_{2L}^{\nu}),~(a_1-l_{1L}^{\nu})(a_2-l_{2L}^{\nu}) \}-a_1a_2,$\\ 
\noindent $l_U^{\prime\nu}=a_1a_2-\mbox{min} \{(a_1+r_{1U}^{\prime\nu})(a_2+r_{2U}^{\prime\nu}),~(a_1-l_{1U}^{\prime\nu})(a_2+r_{2U}^{\prime\nu}) \},$\\
\noindent $r_U^{\prime\nu}=\mbox{max} \{(a_1-l_{1U}^{\prime\nu})(a_2-l_{2U}^{\prime\nu}),~(a_1+r_{1U}^{\prime\nu})(a_2-l_{2U}^{\prime\nu}) \}-a_1a_2,$\\ 
\noindent where the conditions for $LR-$type representation of $\tilde{A}_1\odot \tilde{A}_2$ are satisfied.}\\

\noindent{\bf{Proof.}} Similar to the Proposition $3.1.3$.\\

\noindent{\bf{Proposition 3.1.11.}} \textit{Let $\tilde{A}_1=(a_1; l_{1L}^{\mu},r_{1L}^{\mu},l'^{\mu}_{1U},r'^{\mu}_{1U};l^{\nu}_{1L},r^{\nu}_{1L},l'^{\nu}_{1U},r'^{\nu}_{1U})_{LR}$ be an $LR$-type IVIFN such that $a_1+r_{1L}^{\nu} < 0$\\
 and $\tilde{A}_2=(a_2; l_{2L}^{\mu},r_{2L}^{\mu},l'^{\mu}_{2U},r'^{\mu}_{2U};l^{\nu}_{2L},r^{\nu}_{2L}, l'^{\nu}_{2U},r'^{\nu}_{2U})_{LR}$ be any $LR$-type IVIFN. Then\\
  $\tilde{A}_1 \odot \tilde{A}_2=(a; l_{L}^{\mu},r_{L}^{\mu},l'^{\mu}_{U},r'^{\mu}_{U};l^{\nu}_{L},r^{\nu}_{L},l'^{\nu}_{U},r'^{\nu}_{U})_{LR}$ where\\  
\noindent $a=a_1a_2,$\\
\noindent $l_L^{\mu}=a_1a_2-\mbox{min} \{(a_1+r_{1L}^{\mu})(a_2+r_{2L}^{\mu}),~(a_1-l_{1L}^{\mu})(a_2+r_{2L}^{\mu}) \},$\\
\noindent $r_L^{\mu}=\mbox{max} \{(a_1-l_{1L}^{\mu})(a_2-l_{2L}^{\mu}),~(a_1+r_{1L}^{\mu})(a_2-l_{2L}^{\mu}) \}-a_1a_2,$\\ 
\noindent $l_U^{\prime\mu}=a_1a_2-\mbox{min} \{(a_1+r_{1U}^{\prime\mu})(a_2+r_{2U}^{\prime\mu}),~(a_1-l_{1U}^{\prime\mu})(a_2+r_{2U}^{\prime\mu}) \},$\\
\noindent $r_U^{\prime\mu}=\mbox{max} \{(a_1-l_{1U}^{\prime\mu})(a_2-l_{2U}^{\prime\mu}),~(a_1+r_{1U}^{\prime\mu})(a_2-l_{2U}^{\prime\mu}) \}-a_1a_2,$\\ 
\noindent $l_L^{\nu}=a_1a_2-\mbox{min} \{(a_1+r_{1L}^{\nu})(a_2+r_{2L}^{\nu}),~(a_1-l_{1L}^{\nu})(a_2+r_{2L}^{\nu}) \},$\\
\noindent $r_L^{\nu}=\mbox{max} \{(a_1-l_{1L}^{\nu})(a_2-l_{2L}^{\nu}),~(a_1+r_{1L}^{\nu})(a_2-l_{2L}^{\nu}) \}-a_1a_2,$\\ 
\noindent $l_U^{\prime\nu}=a_1a_2-\mbox{min} \{(a_1+r_{1U}^{\prime\nu})(a_2+r_{2U}^{\prime\nu}),~(a_1-l_{1U}^{\prime\nu})(a_2+r_{2U}^{\prime\nu}) \},$\\
\noindent $r_U^{\prime\nu}=\mbox{max} \{(a_1-l_{1U}^{\prime\nu})(a_2-l_{2U}^{\prime\nu}),~(a_1+r_{1U}^{\prime\nu})(a_2-l_{2U}^{\prime\nu}) \}-a_1a_2,$\\ 
\noindent where the conditions for $LR-$type representation of $\tilde{A}_1\odot \tilde{A}_2$ are satisfied.}\\

\noindent{\bf{Proof.}} Similar to the Proposition $3.1.3$.\\

\noindent{\bf{Proposition 3.1.12.}} \textit{Let $\tilde{A}_1=(a_1; l_{1L}^{\mu},r_{1L}^{\mu},l'^{\mu}_{1U},r'^{\mu}_{1U};l^{\nu}_{1L},r^{\nu}_{1L},l'^{\nu}_{1U},r'^{\nu}_{1U})_{LR}$ be an $LR$-type IVIFN such that $a_1-l_{1L}^{\nu} \geq 0$\\
 and $\tilde{A}_2=(a_2; l_{2L}^{\mu},r_{2L}^{\mu},l'^{\mu}_{2U},r'^{\mu}_{2U}; l^{\nu}_{2L},r^{\nu}_{2L}, l'^{\nu}_{2U},r'^{\nu}_{2U})_{LR}$ be any $LR$-type IVIFN. Then\\
  $\tilde{A}_1 \odot \tilde{A}_2=(a; l_{L}^{\mu},r_{L}^{\mu},l'^{\mu}_{U},r'^{\mu}_{U};l^{\nu}_{L},r^{\nu}_{L},l'^{\nu}_{U},r'^{\nu}_{U})_{LR}$ where\\ 
\noindent $a=a_1a_2,$\\
\noindent $l_L^{\mu}=a_1a_2-\mbox{min} \{(a_1-l_{1L}^{\mu})(a_2-l_{2L}^{\mu}),~(a_1+r_{1L}^{\mu})(a_2-l_{2L}^{\mu}) \},$\\
\noindent $r_L^{\mu}=\mbox{max} \{(a_1-l_{1L}^{\mu})(a_2+r_{2L}^{\mu}),~(a_1+r_{1L}^{\mu})(a_2+r_{2L}^{\mu}) \}-a_1a_2,$\\ 
\noindent $l_U^{\prime\mu}=a_1a_2-\mbox{min} \{(a_1-l_{1U}^{\prime\mu})(a_2-l_{2U}^{\prime\mu}),~(a_1+r_{1U}^{\prime\mu})(a_2-l_{2U}^{\prime\mu}) \},$\\
\noindent $r_U^{\prime\mu}=\mbox{max} \{(a_1-l_{1U}^{\prime\mu})(a_2+r_{2U}^{\prime\mu}),~(a_1+r_{1U}^{\prime\mu})(a_2+r_{2U}^{\prime\mu}) \}-a_1a_2,$\\ 
\noindent $l_L^{\nu}=a_1a_2-\mbox{min} \{(a_1-l_{1L}^{\nu})(a_2-l_{2L}^{\nu}),~(a_1+r_{1L}^{\nu})(a_2-l_{2L}^{\nu}) \},$\\
\noindent $r_L^{\nu}=\mbox{max} \{(a_1-l_{1L}^{\nu})(a_2+r_{2L}^{\nu}),~(a_1+r_{1L}^{\nu})(a_2+r_{2L}^{\nu}) \}-a_1a_2,$\\ 
\noindent $l_U^{\prime\nu}=a_1a_2-\mbox{min} \{(a_1-l_{1U}^{\prime\nu})(a_2-l_{2U}^{\prime\nu}),~(a_1+r_{1U}^{\prime\nu})(a_2-l_{2U}^{\prime\nu}) \},$\\
\noindent $r_U^{\prime\nu}=\mbox{max} \{(a_1+r_{1U}^{\prime\nu})(a_2+r_{2U}^{\prime\nu}),~(a_1-l_{1U}^{\prime\nu})(a_2+r_{2U}^{\prime\nu}) \}-a_1a_2,$\\  
\noindent where the conditions for $LR-$type representation of $\tilde{A}_1\odot \tilde{A}_2$ are satisfied.}\\

\noindent{\bf{Proof.}} Similar to the Proposition $3.1.3$.

\subsection{Proposed lexicographic criteria for the ranking of $LR$-type IVIFNs}
In the literature, there exist several ranking criteria to define the ordering of FNs, IFNs and IVIFNs. Most of the researchers have used linear ranking function which maps the set of these numbers to real-line and then a comparison is carried out on the basis of the usual ordering of real numbers. However, as pointed out in Pérez-Cañedo and Concepción-Morales \cite{ref48} that it may possible that two of these numbers may look different for the decision-maker but can be mapped to the same real number. To overcome this limitation, lexicographic ranking criteria with total order properties seems to be more logical and relevant \cite{ref19, ref33, ref47}. In this section, we discuss the lexicographic criteria for the ordering of $LR$-type IVIFNs.\\
 
\noindent{\bf{Definition 3.2.1}} \cite{ref47}. For $x,y \in \mathbb{R}^n$, the strict lexicographic inequality $x \prec_{lex} y$ holds, \textit{iff} there is $1 \leq i \leq n$ so that $x_j=y_j$ holds for $j<i$ and $x_i<y_i$. The weak lexicographic inequality $x \preceq_{lex} y$ holds, \textit{iff} $x \prec_{lex} y$ or $x=y$.\\

\noindent{\bf{Definition 3.2.2}} Let $\tilde{A}_1=(a_1; l_{1L}^{\mu},r_{1L}^{\mu},l'^{\mu}_{1U},r'^{\mu}_{1U};l^{\nu}_{1L},r^{\nu}_{1L},l'^{\nu}_{1U}, r'^{\nu}_{1U})_{LR}$ and $\tilde{A}_2=(a_2; l_{2L}^{\mu},r_{2L}^{\mu},l'^{\mu}_{2U},r'^{\mu}_{2U};l^{\nu}_{2L},r^{\nu}_{2L},l'^{\nu}_{2U},r'^{\nu}_{2U})_{LR}$ be two $LR$-type IVIFNs, then $\tilde{A}_1=\tilde{A}_2$ \textit{iff}
$$a_1=a_2, l_{1L}^{\mu}=l_{2L}^{\mu},r_{1L}^{\mu}=r_{2L}^{\mu},l'^{\mu}_{1U}=l'^{\mu}_{2U},r'^{\mu}_{1U}=r'^{\mu}_{2U},l^{\nu}_{1L}=l^{\nu}_{2L},r^{\nu}_{1L}=r^{\nu}_{2L},l'^{\nu}_{1U}=l'^{\nu}_{2U},r'^{\nu}_{1U}=r'^{\nu}_{2U}.$$ 

\noindent{\bf{Definition 3.2.3}} Let $\preceq_{lex}$ be the lexicographic order relation on $\mathbb{R}^n$. For $\tilde{A}=(a; l_{L}^{\mu},r_{L}^{\mu},l'^{\mu}_{U},r'^{\mu}_{U};l^{\nu}_{L},r^{\nu}_{L},l'^{\nu}_{U},r'^{\nu}_{U})_{LR} \in IV(\mathbb{R})$, let $S(\tilde{A})$ and $A(\tilde{A})$ be the score and accuracy indices of $\tilde{A}$, respectively; let $M(\tilde{A}):= a$ (the mean value of $\tilde{A}$), $C(\tilde{A}):= a-l_{L}^{\mu}$, $D(\tilde{A}):= a-l'^{\mu}_{U}$, $G(\tilde{A}):= a-l'^{\nu}_{U}$ and $H(\tilde{A}):= a-l^{\nu}_{L}$. Then, for any $\tilde{A}_1, \tilde{A}_2 \in IV(\mathbb{R})$, the strict inequality $\tilde{A}_1 \prec \tilde{A}_2$  holds \textit{iff}\\

\noindent $\big(S(\tilde{A}_1),A(\tilde{A}_1),M(\tilde{A}_1),C(\tilde{A}_1),D(\tilde{A}_1),G(\tilde{A}_1),H(\tilde{A}_1)\big) \prec_{lex} \big(S(\tilde{A}_2),A(\tilde{A}_2),M(\tilde{A}_2),C(\tilde{A}_2),D(\tilde{A}_2),G(\tilde{A}_2),H(\tilde{A}_2)\big).$\\

\noindent The weak inequality $\tilde{A}_1 \preceq \tilde{A}_2$ holds \textit{iff} either\\

\noindent $\big(S(\tilde{A}_1),A(\tilde{A}_1),M(\tilde{A}_1),C(\tilde{A}_1),D(\tilde{A}_1),G(\tilde{A}_1),H(\tilde{A}_1)\big) \prec_{lex} \big(S(\tilde{A}_2),A(\tilde{A}_2),M(\tilde{A}_2),C(\tilde{A}_2),D(\tilde{A}_2),G(\tilde{A}_2),H(\tilde{A}_2)\big)$\\

\noindent $\mbox{or}~~ \big(S(\tilde{A}_1),A(\tilde{A}_1),M(\tilde{A}_1),C(\tilde{A}_1),D(\tilde{A}_1),G(\tilde{A}_1),H(\tilde{A}_1)\big) = \big(S(\tilde{A}_2),A(\tilde{A}_2),M(\tilde{A}_2),C(\tilde{A}_2),D(\tilde{A}_2),G(\tilde{A}_2),H(\tilde{A}_2)\big).$\\

\noindent{\bf{Remark 3.2.1}} The specific order in which the functions $S, A, M, C, D, G~\mbox{and}~H$ appear in the ranking of $LR$-type IVIFNs gives a description of the relative importance of these functions in deciding the ordering of the IVIFNs. However, one may consider a different permutation of these parameters depending on the priority of these functions.\\
 
\noindent{\bf{Remark~ 3.2.2}} Since the expressions of all the functions $S,A,M,C,D,G,H$ are linear, therefore, each $\Phi \in \{S,A,M,\\C,D,G,H \}$ follows the linearity property, that is,
$$\Phi(\lambda_1 \tilde{A}_1\oplus \lambda_2 \tilde{A}_2)=\lambda_1\Phi(\tilde{A}_1)+\lambda_2\Phi(\tilde{A}_2)~~~~\mbox{for all}~~ \tilde{A}_1,\tilde{A}_2 \in IV(\mathbb{R})~ \mbox{and}~ \lambda_1,\lambda_2 \in \mathbb{R}.$$
 
\noindent{\bf{Theorem 3.2.1.}} \textit{The order relation on $IV(\mathbb{R})$ given in Definition 3.2.3, has the total order properties and yields a complete ranking on the set of all $LR$-type IVIFNs of the same type.}\\
\noindent{\bf{Proof.}} The order relation given in Definition $3.2.3$ is a total order due to its following properties:
\begin{enumerate}[$1.$]
\item (reflexivity) $\tilde{A}_1 \preceq \tilde{A}_1~~\forall~\tilde{A}_1 \in IV(\mathbb{R})$,
\item (anti-symmetry) $\tilde{A}_1 \preceq \tilde{A}_2$ and $\tilde{A}_2 \preceq \tilde{A}_1 \\
~~~~~~~~~\implies \tilde{A}_1=\tilde{A}_2~~\forall~\tilde{A}_1,\tilde{A}_2 \in IV(\mathbb{R})$,
\item (transitivity) $\tilde{A}_1 \preceq \tilde{A}_2$ and $\tilde{A}_2 \preceq \tilde{A}_3 \\
~~~~~~~~~\implies \tilde{A}_1 \preceq \tilde{A}_3~~\forall~\tilde{A}_1,\tilde{A}_2,\tilde{A}_3 \in IV(\mathbb{R})$,
\item (comparability) $\tilde{A}_1 \preceq \tilde{A}_2$ or $\tilde{A}_2 \preceq \tilde{A}_1~~\forall~\tilde{A}_1,\tilde{A}_2 \in IV(\mathbb{R})$.
\end{enumerate}
Further, it is to be noted that the system of linear equations $S(\tilde{A}_1)=S(\tilde{A}_2),~A(\tilde{A}_1)=A(\tilde{A}_2),~M(\tilde{A}_1)=M(\tilde{A}_2),\\
C(\tilde{A}_1)=C(\tilde{A}_2),~D(\tilde{A}_1)=D(\tilde{A}_2),~G(\tilde{A}_1)=G(\tilde{A}_2)$ and $~H(\tilde{A}_1)=H(\tilde{A}_2)$ has a unique solution, that is, $\tilde{A}_1=\tilde{A}_2$ as the absolute value of the determinant of this linear system is $2bd$ where $b=\int_0^1 R^{-1}(x)dx >0$ and $d=\int_0^1 R^{-1}(1-x)dx >0$. Hence, the result. \\

\section{Model formulation and solution algorithm}
Consider a linear programming problem in which each parameter $\tilde{c}_j,~ \tilde{a}_{ij},~\tilde{b}_i$ and decision variables $\tilde{x}_j$ are taken to be in the form of $LR$-type IVIFNs. The model of such a fully $LR$-type interval-valued intuitionistic fuzzy linear programming problem along-with unrestricted decision variables can be mathematically formulated as:
\begin{center}
$\hspace{-6.5cm}\textbf{(P1)} ~\max~\tilde{Z}= \displaystyle\sum_{j=1}^{n} \tilde{c}_j \odot \tilde{x}_j$\\
$\hspace{-1.1cm}\mbox{s.t.}~\displaystyle\sum_{j=1}^{n} \tilde{a}_{ij} \odot \tilde{x}_j = \tilde{b}_i,~~~\mbox{for}~i \in I_1:=\{ 1,2, \dots, m_1 \},$\\
$~~~~~~~~~~~~\displaystyle\sum_{j=1}^{n} \tilde{a}_{ij} \odot \tilde{x}_j \preceq \tilde{b}_i,~~~~\mbox{for}~i \in I_2:=\{m_1+1,m_1+2, \dots, m_2 \},$\\
$~~~~~~~~~~~\displaystyle\sum_{j=1}^{n} \tilde{a}_{ij} \odot \tilde{x}_j \succeq \tilde{b}_i,~~~~\mbox{for}~i \in I_3:=\{m_2+1,m_2+2, \dots, m \},$\\
$~~~~\tilde{x}_j~\mbox{are unrestricted in sign},~\mbox{for}~j \in J:=\{ 1,2, \dots, n \}$\\
\end{center}
where the inequalities "$\preceq$" and "$\succeq$" in the problem (P1) is in accordance with the lexicographic ranking given in Definition 3.2.3.\\

\noindent{\bf{Definition 4.1}} Let $X$ denote the set of all feasible solutions of (P1). A vector $\hat{\tilde{x}}=(\hat{\tilde{x}}_1,\hat{\tilde{x}}_2, \dots,\hat{\tilde{x}}_n) \in X$ is said to be an optimal solution of (P1) if
\begin{center}
$\displaystyle\sum_{j=1}^{n} \tilde{c}_j \odot \tilde{x}_j \preceq \displaystyle\sum_{j=1}^{n} \tilde{c}_j \odot \hat{\tilde{x}}_j,~~$ for all $\tilde{x}=(\tilde{x}_1,\tilde{x}_2, \dots, \tilde{x}_n) \in X$. 
\end{center}

 Next, we propose a method to find the unique optimal solution of (P1). We use lexicographic ranking criteria for ordering of $LR$-type IVIFNs. The steps of the proposed approach are as follows:
 \begin{itemize}
 \item[\textbf{Step 1.}] Let $~\tilde{a}_{ij}=(a_{ij};\alpha^{\mu}_{ijL},\beta^{\mu}_{ijL},\alpha'^{\mu}_{ijU},\beta'^{\mu}_{ijU};\alpha^{\nu}_{ijL},\beta^{\nu}_{ijL}, \alpha'^{\nu}_{ijU},\beta'^{\nu}_{ijU})_{LR}$, $~~\tilde{c}_{j}=(c_{j};\sigma^{\mu}_{jL},\rho^{\mu}_{jL},\sigma'^{\mu}_{jU}, \rho'^{\mu}_{jU};\sigma^{\nu}_{jL},\rho^{\nu}_{jL},\sigma'^{\nu}_{jU},\rho'^{\nu}_{jU})_{LR},$\\\\
$\tilde{b}_{i}=(b_{i};\gamma^{\mu}_{iL},\delta^{\mu}_{iL},\gamma'^{\mu}_{iU},\delta'^{\mu}_{iU};\gamma^{\nu}_{iL},\delta^{\nu}_{iL},\gamma'^{\nu}_{iU},\delta'^{\nu}_{iU})_{LR}$,~ and $~\tilde{x}_j=(x_j;\xi_{jL}^{\mu},\eta_{jL}^{\mu},\xi_{jU}^{\prime\mu},\eta_{jU}^{\prime\mu};\xi_{jL}^{\nu},\eta_{jL}^{\nu},\xi_{jU}^{\prime\nu},\eta_{jU}^{\prime\nu})_{LR}.$\\

\noindent Then, (P1) can be recast as:\\\\
\noindent$\max~ \displaystyle\sum_{j=1}^{n} (c_{j};\sigma^{\mu}_{jL},\rho^{\mu}_{jL},\sigma'^{\mu}_{jU}, \rho'^{\mu}_{jU};\sigma^{\nu}_{jL},\rho^{\nu}_{jL},\sigma'^{\nu}_{jU},\rho'^{\nu}_{jU})_{LR} \odot (x_j;\xi_{jL}^{\mu},\eta_{jL}^{\mu},\xi_{jU}^{\prime\mu},\eta_{jU}^{\prime\mu};\xi_{jL}^{\nu},\eta_{jL}^{\nu},\xi_{jU}^{\prime\nu},\eta_{jU}^{\prime\nu})_{LR}$\\\\
$~~~\mbox{s.t.}~~\displaystyle\sum_{j=1}^{n} (a_{ij};\alpha^{\mu}_{ijL},\beta^{\mu}_{ijL},\alpha'^{\mu}_{ijU},\beta'^{\mu}_{ijU};\alpha^{\nu}_{ijL},\beta^{\nu}_{ijL},\alpha'^{\nu}_{ijU},\beta'^{\nu}_{ijU})_{LR} \odot (x_j;\xi_{jL}^{\mu},\eta_{jL}^{\mu},\xi_{jU}^{\prime\mu},\eta_{jU}^{\prime\mu};\xi_{jL}^{\nu},\eta_{jL}^{\nu},\xi_{jU}^{\prime\nu},\eta_{jU}^{\prime\nu})_{LR} $\\\\
$~~~~~~~~~~=(b_{i};\gamma^{\mu}_{iL},\delta^{\mu}_{iL},\gamma'^{\mu}_{iU},\delta'^{\mu}_{iU};\gamma^{\nu}_{iL},\delta^{\nu}_{iL},\gamma'^{\nu}_{iU},\delta'^{\nu}_{iU})_{LR},~~~~~~~\mbox{for}~i \in I_1,$\\\\
$~~~~~~~~~~\displaystyle\sum_{j=1}^{n} (a_{ij};\alpha^{\mu}_{ijL},\beta^{\mu}_{ijL},\alpha'^{\mu}_{ijU},\beta'^{\mu}_{ijU};\alpha^{\nu}_{ijL},\beta^{\nu}_{ijL},\alpha'^{\nu}_{ijU},\beta'^{\nu}_{ijU})_{LR} \odot (x_j;\xi_{jL}^{\mu},\eta_{jL}^{\mu},\xi_{jU}^{\prime\mu},\eta_{jU}^{\prime\mu};\xi_{jL}^{\nu},\eta_{jL}^{\nu},\xi_{jU}^{\prime\nu},\eta_{jU}^{\prime\nu})_{LR} $\\\\
$~~~~~~~~~~~\preceq (b_{i};\gamma^{\mu}_{iL},\delta^{\mu}_{iL},\gamma'^{\mu}_{iU},\delta'^{\mu}_{iU};\gamma^{\nu}_{iL},\delta^{\nu}_{iL},\gamma'^{\nu}_{iU},\delta'^{\nu}_{iU})_{LR},~~~~~~~\mbox{for}~i \in I_2,$\\\\
$~~~~~~~~~~\displaystyle\sum_{j=1}^{n} (a_{ij};\alpha^{\mu}_{ijL},\beta^{\mu}_{ijL},\alpha'^{\mu}_{ijU},\beta'^{\mu}_{ijU};\alpha^{\nu}_{ijL},\beta^{\nu}_{ijL},\alpha'^{\nu}_{ijU},\beta'^{\nu}_{ijU})_{LR} \odot (x_j;\xi_{jL}^{\mu},\eta_{jL}^{\mu},\xi_{jU}^{\prime\mu},\eta_{jU}^{\prime\mu};\xi_{jL}^{\nu},\eta_{jL}^{\nu},\xi_{jU}^{\prime\nu},\eta_{jU}^{\prime\nu})_{LR}$\\\\
$~~~~~~~~~~~\succeq (b_{i};\gamma^{\mu}_{iL},\delta^{\mu}_{iL},\gamma'^{\mu}_{iU},\delta'^{\mu}_{iU};\gamma^{\nu}_{iL},\delta^{\nu}_{iL},\gamma'^{\nu}_{iU},\delta'^{\nu}_{iU})_{LR},~~~~~~~\mbox{for}~i \in I_3,$\\\\
$~~~~~~~~~~~(x_j;\xi_{jL}^{\mu},\eta_{jL}^{\mu},\xi_{jU}^{\prime\mu},\eta_{jU}^{\prime\mu};\xi_{jL}^{\nu},\eta_{jL}^{\nu},\xi_{jU}^{\prime\nu},\eta_{jU}^{\prime\nu})_{LR}~~~\mbox{are unrestricted in sign},~\mbox{for}~j \in J$.

\item[\textbf{Step 2.}] After applying the multiplication operation (Section $3.1$), let\\\\
\noindent $~~~~~~~~~~~~~~~~~~~~\tilde{c}_{j}\odot \tilde{x}_j=(p_{j};\tau^{\mu}_{jL},\omega^{\mu}_{jL},\tau'^{\mu}_{jU}, \omega'^{\mu}_{jU};\tau^{\nu}_{jL},\omega^{\nu}_{jL},\tau'^{\nu}_{jU},\omega'^{\nu}_{jU})_{LR}~~$ and \\\\
\noindent $~~~~~~~~~~~~~~~~~~~~\tilde{a}_{ij} \odot \tilde{x}_j=(m_{ij};s^{\mu}_{ijL},\lambda^{\mu}_{ijL},s'^{\mu}_{ijU},\lambda'^{\mu}_{ijU};s^{\nu}_{ijL},\lambda^{\nu}_{ijL},s'^{\nu}_{ijU},\lambda'^{\nu}_{ijU})_{LR}$.\\

\noindent Then, the problem (P1) in step 1 becomes:
\begin{center}
\textbf{(P2)} $\max~ \displaystyle\sum_{j=1}^{n}(p_{j};\tau^{\mu}_{jL},\omega^{\mu}_{jL},\tau'^{\mu}_{jU}, \omega'^{\mu}_{jU};\tau^{\nu}_{jL},\omega^{\nu}_{jL},\tau'^{\nu}_{jU},\omega'^{\nu}_{jU})_{LR}~~~~~~~~~~~~~~~~~~~$ \\
$\mbox{s.t.}~~\displaystyle\sum_{j=1}^{n}(m_{ij};s^{\mu}_{ijL},\lambda^{\mu}_{ijL},s'^{\mu}_{ijU},\lambda'^{\mu}_{ijU};s^{\nu}_{ijL},\lambda^{\nu}_{ijL},s'^{\nu}_{ijU},\lambda'^{\nu}_{ijU})_{LR} ~$\\
$~~~~~~~~~~~~~= (b_{i};\gamma^{\mu}_{iL},\delta^{\mu}_{iL},\gamma'^{\mu}_{iU},\delta'^{\mu}_{iU};\gamma^{\nu}_{iL},\delta^{\nu}_{iL},\gamma'^{\nu}_{iU},\delta'^{\nu}_{iU})_{LR},~~~\mbox{for}~i \in I_1,$\\
$~~~~~~\displaystyle\sum_{j=1}^{n}(m_{ij};s^{\mu}_{ijL},\lambda^{\mu}_{ijL},s'^{\mu}_{ijU},\lambda'^{\mu}_{ijU};s^{\nu}_{ijL},\lambda^{\nu}_{ijL},s'^{\nu}_{ijU},\lambda'^{\nu}_{ijU})_{LR}$\\
$~~~~~~~~~~~~~\preceq (b_{i};\gamma^{\mu}_{iL},\delta^{\mu}_{iL},\gamma'^{\mu}_{iU},\delta'^{\mu}_{iU};\gamma^{\nu}_{iL},\delta^{\nu}_{iL},\gamma'^{\nu}_{iU},\delta'^{\nu}_{iU})_{LR},~~~\mbox{for}~i \in I_2,$\\
$~~~~~~~\displaystyle\sum_{j=1}^{n}(m_{ij};s^{\mu}_{ijL},\lambda^{\mu}_{ijL},s'^{\mu}_{ijU},\lambda'^{\mu}_{ijU};s^{\nu}_{ijL},\lambda^{\nu}_{ijL},s'^{\nu}_{ijU},\lambda'^{\nu}_{ijU})_{LR}$\\
$~~~~~~~~~~~~~\succeq (b_{i};\gamma^{\mu}_{iL},\delta^{\mu}_{iL},\gamma'^{\mu}_{iU},\delta'^{\mu}_{iU};\gamma^{\nu}_{iL},\delta^{\nu}_{iL},\gamma'^{\nu}_{iU},\delta'^{\nu}_{iU})_{LR},~~~\mbox{for}~i \in I_3,$\\
$~~~~~~~~~~~~~~~~~~~~~~~~~~~~~~~~~~~~~~~~~~~~~~~~~(x_j;\xi_{jL}^{\mu},\eta_{jL}^{\mu},\xi_{jU}^{\prime\mu},\eta_{jU}^{\prime\mu};\xi_{jL}^{\nu},\eta_{jL}^{\nu},\xi_{jU}^{\prime\nu},\eta_{jU}^{\prime\nu})_{LR}~~\mbox{are unrestricted in sign},~\mbox{for}~j \in J$.
\end{center}

\item[\textbf{Step 3.}] Define $\tilde{l}_i:=\displaystyle\sum_{j=1}^{n}(m_{ij};s^{\mu}_{ijL},\lambda^{\mu}_{ijL},s'^{\mu}_{ijU},\lambda'^{\mu}_{ijU};s^{\nu}_{ijL},\lambda^{\nu}_{ijL},s'^{\nu}_{ijU},\lambda'^{\nu}_{ijU})_{LR}$ ~~and\\
 $\tilde{r}_i:=(b_{i};\gamma^{\mu}_{iL},\delta^{\mu}_{iL},\gamma'^{\mu}_{iU},\delta'^{\mu}_{iU};\gamma^{\nu}_{iL},\delta^{\nu}_{iL},\gamma'^{\nu}_{iU},\delta'^{\nu}_{iU})_{LR}$.\\

\noindent Further, by taking the order relation $\preceq_{lex}$ as in Definition $3.2.3$, according to the addition operation from Section 3.1 and the equality between $LR$-type IVIFNs given in Definition $3.2.2$, we can write the constraint set of (P2) as:\\\\
$\displaystyle\sum_{j=1}^{n}m_{ij}=b_{i},~\displaystyle\sum_{j=1}^{n}s^{\mu}_{ijL}=\gamma^{\mu}_{iL},~\displaystyle\sum_{j=1}^{n}\lambda^{\mu}_{ijL}=\delta^{\mu}_{iL},~\displaystyle\sum_{j=1}^{n}s'^{\mu}_{ijU}=\gamma'^{\mu}_{iU},~\displaystyle\sum_{j=1}^{n}\lambda'^{\mu}_{ijU}=\delta'^{\mu}_{iU}, \displaystyle\sum_{j=1}^{n}s^{\nu}_{ijL}=\gamma^{\nu}_{iL},$\\\\
$\displaystyle\sum_{j=1}^{n}\lambda^{\nu}_{ijL}=\delta^{\nu}_{iL},~\displaystyle\sum_{j=1}^{n}s'^{\nu}_{ijU}=\gamma'^{\nu}_{iU},~\displaystyle\sum_{j=1}^{n}\lambda'^{\nu}_{ijU}=\delta'^{\nu}_{iU},~~\mbox{for}~i \in I_1,~~~~~~~~~~~~~~~~~~~~~~~~~~~~~~~~~~~~~~~~~~~~~~~~~~~~~~~~~~~~~~~~~~~~~~~~~~~~~~(3a)$\\\\
$\big(S(\tilde{l}_i),A(\tilde{l}_i),M(\tilde{l}_i),C(\tilde{l}_i),D(\tilde{l}_i),G(\tilde{l}_i),H(\tilde{l}_i)\big)\preceq_{lex} \big(S(\tilde{r}_i),A(\tilde{r}_i),M(\tilde{r}_i),C(\tilde{r}_i),D(\tilde{r}_i),G(\tilde{r}_i),H(\tilde{r}_i)\big),~~\mbox{for}~i \in I_2,$\\$~~~~~~~~~~~~~~~~~~~~~~~~~~~~~~~~~~~~~~~~~~~~~~~~~~~~~~~~~~~~~~~~~~~~~~~~~~~~~~~~~~~~~~~~~~~~~~~~~~~~~~~~~~~~~~~~~~~~~~~~~~~~~~~~~~~~~~~~~~~~~~~~~~~~~~~~~~~~~~~~~~~~~~~~~~~~~~~~~~~~~~~~~~~(3b)$\\\\
$\big(S(\tilde{l}_i),A(\tilde{l}_i),M(\tilde{l}_i),C(\tilde{l}_i),D(\tilde{l}_i),G(\tilde{l}_i),H(\tilde{l}_i)\big)\succeq_{lex} \big(S(\tilde{r}_i),A(\tilde{r}_i),M(\tilde{r}_i),C(\tilde{r}_i),D(\tilde{r}_i),G(\tilde{r}_i),H(\tilde{r}_i)\big),~~\mbox{for}~i \in I_3,$\\\
$~~~~~~~~~~~~~~~~~~~~~~~~~~~~~~~~~~~~~~~~~~~~~~~~~~~~~~~~~~~~~~~~~~~~~~~~~~~~~~~~~~~~~~~~~~~~~~~~~~~~~~~~~~~~~~~~~~~~~~~~~~~~~~~~~~~~~~~~~~~~~~~~~~~~~~~~~~~~~~~~~~~~~~~~~~~~~~~~~~~~~~~~~~~(3c)$\\\\
$~\xi_{jU}^{\prime\mu} \geq \xi_{jL}^{\mu},~\eta_{jU}^{\prime\mu} \geq \eta_{jL}^{\mu},~\xi_{jL}^{\nu} \geq \xi_{jU}^{\prime\nu},~\eta_{jL}^{\nu} \geq \eta_{jU}^{\prime\nu},~\xi_{jL}^{\nu} \geq \xi_{jL}^{\mu},~\eta_{jL}^{\nu} \geq \eta_{jL}^{\mu},~\xi_{jU}^{\prime\nu} \geq \xi_{jU}^{\prime\mu},~\eta_{jU}^{\prime\nu} \geq \eta_{jU}^{\prime\mu},$\\
$~\xi_{jL}^{\mu} \geq 0,~\eta_{jL}^{\mu} \geq 0,~\xi_{jU}^{\prime\mu} \geq 0,~\eta_{jU}^{\prime\mu} \geq 0,~\xi_{jL}^{\nu} \geq 0,~\eta_{jL}^{\nu} \geq 0,~\xi_{jU}^{\prime\nu} \geq 0,~\eta_{jU}^{\prime\nu} \geq 0,~~\mbox{for}~j \in J.~~~~~~~~~~~~~~~~~~~~~~~~~~~~~~~~~~~~~~~~~~~~~~~~~~~~(3d)$\\\\
\noindent The objective function of the problem (P2) can be recast as:
$$\mbox{lex max}~~\big(S(\tilde{Z}),A(\tilde{Z}),M(\tilde{Z}),C(\tilde{Z}),D(\tilde{Z}), G(\tilde{Z}),H(\tilde{Z})\big)$$

\noindent where $\tilde{Z}= \displaystyle\sum_{j=1}^{n} \tilde{c}_j \odot \tilde{x}_j=\displaystyle\sum_{j=1}^{n}(p_{j};\tau^{\mu}_{jL},\omega^{\mu}_{jL},\tau'^{\mu}_{jU}, \omega'^{\mu}_{jU};\tau^{\nu}_{jL},\omega^{\nu}_{jL},\tau'^{\nu}_{jU},\omega'^{\nu}_{jU})_{LR}.$\\\\
Hence, the problem (P1) is finally transformed into the following optimization problem (P3):
\begin{center}
\textbf{(P3)}~~lex max$~\big(S(\tilde{Z}),A(\tilde{Z}),M(\tilde{Z}),C(\tilde{Z}),D(\tilde{Z}), G(\tilde{Z}),H(\tilde{Z})\big)$\\
$\mbox{subject to constraints $(3a)$ \textendash $(3d)$}~~~$.
\end{center}

\item[\textbf{Step 4.}] By introducing new binary variables $z_i^S, z_i^A, z_i^M, z_i^C,z_i^D, z_i^G~\mbox{and}~z_i^H$ for $i \in I_2 \cup I_3$ following Pérez-Cañedo and Concepción-Morales \cite{ref48} along-with a new constraint $(4p)$, for positive real values of $k$ and $K$ sufficiently small and large respectively, the lexicographic constraints $(3b)$ and $(3c)$ can be converted into the following set of constraints $(4a)$ \textendash $~(4p)$:\\\\
$kz_i^S \leq S(\tilde{r}_i)-S(\tilde{l}_i) \leq Kz_i^S,~~~~\mbox{for}~i \in I_2,\hspace{9.1cm}(4a)$\\\\
$-Kz_i^S+kz_i^A \leq A(\tilde{r}_i)-A(\tilde{l}_i) \leq Kz_i^A,~~\mbox{for}~i \in I_2, \hspace{8cm} (4b)$\\\\
$-K(z_i^S+z_i^A)+kz_i^M \leq M(\tilde{r}_i)-M(\tilde{l}_i) \leq Kz_i^M,~~\mbox{for}~i \in I_2, \hspace{7.05cm} (4c)$\\\\
$-K(z_i^S+z_i^A+z_i^M)+kz_i^C \leq C(\tilde{r}_i)-C(\tilde{l}_i) \leq Kz_i^C,~~\mbox{for}~i \in I_2, \hspace{6.7cm} (4d)$\\\\
$-K(z_i^S+z_i^A+z_i^M+z_i^C)+kz_i^D \leq D(\tilde{r}_i)-D(\tilde{l}_i) \leq Kz_i^D,~~\mbox{for}~i \in I_2, \hspace{6.05cm} (4e)$\\\\
$-K(z_i^S+z_i^A+z_i^M+z_i^C+z_i^D)+kz_i^G \leq G(\tilde{r}_i)-G(\tilde{l}_i) \leq Kz_i^G,~~\mbox{for}~i \in I_2, \hspace{5.4cm} (4f)$\\\\
$-K(z_i^S+z_i^A+z_i^M+z_i^C+z_i^D+z_i^G)+kz_i^H \leq H(\tilde{r}_i)-H(\tilde{l}_i) \leq Kz_i^H,~~\mbox{for}~i \in I_2, \hspace{4.55cm} (4g)$\\\\
$kz_i^S \leq S(\tilde{l}_i)-S(\tilde{r}_i) \leq Kz_i^S,~~~\mbox{for}~i \in I_3, \hspace{9.3cm} (4h)$\\\\
$-Kz_i^S+kz_i^A \leq A(\tilde{l}_i)-A(\tilde{r}_i) \leq Kz_i^A,~~\mbox{for}~i \in I_3, \hspace{8cm} (4i)$\\\\
$-K(z_i^S+z_i^A)+kz_i^M \leq M(\tilde{l}_i)-M(\tilde{r}_i) \leq Kz_i^M,~~\mbox{for}~i \in I_3, \hspace{7.05cm} (4j)$\\\\
$-K(z_i^S+z_i^A+z_i^M)+kz_i^C \leq C(\tilde{l}_i)-C(\tilde{r}_i) \leq Kz_i^C,~~\mbox{for}~i \in I_3, \hspace{6.3cm} (4k)$\\\\
$-K(z_i^S+z_i^A+z_i^M+z_i^C)+kz_i^D \leq D(\tilde{l}_i)-D(\tilde{r}_i)\leq Kz_i^D,~~\mbox{for}~i \in I_3, \hspace{5.5cm} (4l)$\\\\
$-K(z_i^S+z_i^A+z_i^M+z_i^C+z_i^D)+kz_i^G \leq G(\tilde{l}_i)-G(\tilde{r}_i) \leq Kz_i^G,~~\mbox{for}~i \in I_3, \hspace{4.7cm} (4m)$\\\\
$-K(z_i^S+z_i^A+z_i^M+z_i^C+z_i^D+z_i^G)+kz_i^H \leq H(\tilde{l}_i)-H(\tilde{r}_i) \leq Kz_i^H,~~\mbox{for}~i \in I_3, \hspace{3.45cm} (4n)$\\\\
$z_i^S, z_i^A, z_i^M, z_i^C, z_i^D, z_i^G, z_i^H \in \{0,1 \},~~~~~~~~\mbox{for}~~i \in I_2 \cup I_3, \hspace{7.15cm} (4o)$\\\\
$z_i^S \leq z_i^A \leq z_i^M \leq z_i^C \leq z_i^D \leq z_i^G \leq z_i^H,~~\mbox{for}~~i \in I_2 \cup I_3. \hspace{7.2cm} (4p)$

\item[\textbf{Step 5.}] Convert problem (P3) into the following mixed $0$ \textendash $~1$ lexicographic non-linear programming problem (P4):
\begin{center}
\noindent \hspace{-1cm}\textbf{(P4)}~~lex max$~\big(S(\tilde{Z}),A(\tilde{Z}),M(\tilde{Z}),C(\tilde{Z}),D(\tilde{Z}), G(\tilde{Z}),H(\tilde{Z})\big)$\\
\medspace
~~~~subject to constraints $(3a), (3d), (4a)$ \textendash $~(4p)$.
\end{center}

\item[\textbf{Step 6.}] Solve the problem (P4) by optimizing orderly one objective at a time subject to constraints $(3a), (3d), (4a)$ \textendash $~(4p)$; including all previously optimized objectives in the constraint set. Thus, by using the lexicographic method of multi-objective optimization \cite{ref23} and a suitable optimization solver, we can find the optimal solution $x_j,\xi_{jL}^{\mu},\eta_{jL}^{\mu},\xi_{jU}^{\prime\mu},\eta_{jU}^{\prime\mu},\xi_{jL}^{\nu},\eta_{jL}^{\nu},\xi_{jU}^{\prime\nu},\eta_{jU}^{\prime\nu}$ for $j \in J$. Hence, $\tilde{x}_j=(x_j;\xi_{jL}^{\mu},\eta_{jL}^{\mu},\xi_{jU}^{\prime\mu},\eta_{jU}^{\prime\mu};\xi_{jL}^{\nu}, \eta_{jL}^{\nu},\xi_{jU}^{\prime\nu},\eta_{jU}^{\prime\nu})_{LR}$ for $j \in J$ can be obtained.

\item[\textbf{Step 7.}] Finally, by substituting the values of $\tilde{x}_j$'s into $\displaystyle\sum_{j=1}^{n} \tilde{c}_j \odot \tilde{x}_j$, we obtain the optimal interval-valued intuitionistic fuzzy value of the problem (P1).
\end{itemize}

Now, we present Theorems $4.1$ and $4.2$ to prove the equivalence of the models (P1), (P3) and (P3), (P4), respectively. The similar results for fuzzy LPP were discussed by Pérez-Cañedo and Concepción-Morales \cite{ref47}. However, we have shown the results for $LR$-type IVIFLPPs.\\

\noindent \textbf{Theorem 4.1.} \textit{If $\hat{\tilde{x}}=(\hat{\tilde{x}}_1, \hat{\tilde{x}}_2, \hat{\tilde{x}}_3, \dots , \hat{\tilde{x}}_n)$ is an optimal solution of problem (P3), then it is also an optimal solution of (P1).}\\

\noindent \textbf{Proof.} We will prove the result by the method of contradiction. Let $\hat{\tilde{x}}=(\hat{\tilde{x}}_1, \hat{\tilde{x}}_2, \hat{\tilde{x}}_3, \dots , \hat{\tilde{x}}_n)$ be an optimal solution of problem (P3) but not an optimal solution of (P1). Then, there exists a feasible solution $\tilde{x}^{\ast}=(\tilde{x}_1^{\ast}, \tilde{x}_2^{\ast}, \dots, \tilde{x}_n^{\ast})$ of the problem (P1), such that
$$\displaystyle\sum_{j=1}^{n} \tilde{c}_j \odot \hat{\tilde{x}}_j \prec \displaystyle\sum_{j=1}^{n} \tilde{c}_j \odot \tilde{x}_j^{\ast}.$$
 In view of Definition $3.2.3$, $\tilde{x}^{\ast}$ is a feasible solution of the problem (P3) for which\\
 $\bigg(S\bigg(\displaystyle\sum_{j=1}^{n} \tilde{c}_j \odot \hat{\tilde{x}}_j\bigg),A\bigg(\displaystyle\sum_{j=1}^{n} \tilde{c}_j \odot \hat{\tilde{x}}_j\bigg),M\bigg(\displaystyle\sum_{j=1}^{n} \tilde{c}_j \odot \hat{\tilde{x}}_j\bigg),C\bigg(\displaystyle\sum_{j=1}^{n} \tilde{c}_j \odot \hat{\tilde{x}}_j\bigg),D\bigg(\displaystyle\sum_{j=1}^{n} \tilde{c}_j \odot \hat{\tilde{x}}_j\bigg),G\bigg(\displaystyle\sum_{j=1}^{n} \tilde{c}_j \odot \hat{\tilde{x}}_j\bigg),$\\\\
 $H\bigg(\displaystyle\sum_{j=1}^{n} \tilde{c}_j \odot \hat{\tilde{x}}_j\bigg)\bigg) \prec_{lex} \bigg(S\bigg(\displaystyle\sum_{j=1}^{n} \tilde{c}_j \odot \tilde{x}_j^{\ast}\bigg),A\bigg(\displaystyle\sum_{j=1}^{n} \tilde{c}_j \odot \tilde{x}_j^{\ast}\bigg),M\bigg(\displaystyle\sum_{j=1}^{n} \tilde{c}_j \odot 
\tilde{x}_j^{\ast}\bigg),C\bigg(\displaystyle\sum_{j=1}^{n} \tilde{c}_j \odot \tilde{x}_j^{\ast}\bigg),D\bigg(\displaystyle\sum_{j=1}^{n} \tilde{c}_j \odot \tilde{x}_j^{\ast}\bigg),$\\\\
$G\bigg(\displaystyle\sum_{j=1}^{n} \tilde{c}_j \odot \tilde{x}_j^{\ast}\bigg),H\bigg(\displaystyle\sum_{j=1}^{n} \tilde{c}_j \odot \tilde{x}_j^{\ast}\bigg)\bigg),$\\

\noindent that is, there exists a feasible solution of (P3) with the higher objective function value. This contradicts the fact that $\hat{\tilde{x}}$ is an optimal solution of (P3). Hence, the result.\\

\noindent \textbf{Theorem 4.2.} \textit{The optimal solution of the problem (P4) is also optimal for (P3). The converse of the statement is also true.}\\
\noindent \textbf{Proof.} Let $\tilde{l}_i:=\displaystyle\sum_{j=1}^{n}(m_{ij};s^{\mu}_{ijL},\lambda^{\mu}_{ijL},s'^{\mu}_{ijU},\lambda'^{\mu}_{ijU};s^{\nu}_{ijL},\lambda^{\nu}_{ijL},s'^{\nu}_{ijU},\lambda'^{\nu}_{ijU})_{LR},~~\tilde{r}_i:=(b_{i};\gamma^{\mu}_{iL},\delta^{\mu}_{iL},\gamma'^{\mu}_{iU},\delta'^{\mu}_{iU};\gamma^{\nu}_{iL},\delta^{\nu}_{iL},\gamma'^{\nu}_{iU},\delta'^{\nu}_{iU})_{LR}$\\
 and $z_i:=(z_i^S, z_i^A, z_i^M, z_i^C, z_i^D, z_i^G,z_i^H)$.\\

 In order to prove this theorem, we need to establish the equivalence between constraint sets $(3)[(3b),(3c)]$ and $(4)$\\$[(4a)-(4p)]$ for positive real values of $k$ and $K$ sufficiently small and large, respectively. Here, we will prove the result for the $I_2$ set of constraints only. For the remaining set of constraints, the result can be proved on the same lines. \\
To show that any solution satisfying constraint set (4) also satisfies set (3), we consider the following cases:
\begin{enumerate}[1.]
\item If $z_i=(1, \star,\star,\star,\star,\star,\star)$, where $\star=0~\mbox{or}~1$, then by substituting into constraint set (4), we get $k \leq S(\tilde{r}_i)-S(\tilde{l}_i) \leq K,$ which implies $S(\tilde{l}_i) < S(\tilde{r}_i)$; hence, according to Definition $3.2.3$, set of constraints (3) are satisfied.
\item If $z_i=(0, 1,\star,\star,\star,\star,\star)$, where $\star=0~\mbox{or}~1$, then by the constraint set (4) we get $0 \leq S(\tilde{r}_i)-S(\tilde{l}_i) \leq 0,$ and $k \leq A(\tilde{r}_i)-A(\tilde{l}_i) \leq K,$ which implies $S(\tilde{l}_i) = S(\tilde{r}_i)$ and $A(\tilde{l}_i) < A(\tilde{r}_i)$ ; hence, from Definition $3.2.3$, the constraint set (3) are satisfied.
\end{enumerate}
 
The remaining cases $z_i=(0,0,1,\star,\star,\star,\star)$, $z_i=(0,0,0,1,\star,\star,\star)$, $z_i=(0,0,0,0,1,\star,\star)$, $z_i=(0,0,0,0,0,1,\star)$, $z_i=(0,0,0,0,0,0,1)$ and $z_i=(0,0,0,0,0,0,0)$ can be proved similarly.\\

\noindent Now, it is to be shown that any solution satisfying constraint set $(3)$ also satisfies constraints $(4)$. For this, let us consider the following cases:
\begin{enumerate}[1.]
\item If $S(\tilde{r}_i)=S(\tilde{l}_i)$,  $A(\tilde{r}_i)=A(\tilde{l}_i)$, $M(\tilde{r}_i)=M(\tilde{l}_i)$, $C(\tilde{r}_i)=C(\tilde{l}_i)$, $D(\tilde{r}_i)=D(\tilde{l}_i)$, $G(\tilde{r}_i)=G(\tilde{l}_i)$ and $H(\tilde{r}_i)=H(\tilde{l}_i)$, then by constraint set (4), we get $z_i=(0,0,0,0,0,0,0)$.
\item If $S(\tilde{r}_i)>S(\tilde{l}_i)$, i.e., $S(\tilde{r}_i)-S(\tilde{l}_i)=s_i >0$, then by constraint $(4a)$, we get $kz_i^S \leq s_i \leq Kz_i^S$, this inequality is satisfied for $z_i^S=1$ and for positive real values of $k$ and $K$, sufficiently small and large. Further, $(4o)$ and $(4p)$ give $(z_i^A, z_i^M, z_i^C, z_i^D, z_i^G,z_i^H)=(1,1,1,1,1,1)$, due to which rest of the constraints $(4b)-(4g)$ are automatically satisfied. Also, $S(\tilde{r}_i)>S(\tilde{l}_i)$ implies that the constraint $(3b)$ is satisfied.
\item If $S(\tilde{r}_i)=S(\tilde{l}_i)$ and $A(\tilde{r}_i)>A(\tilde{l}_i)$ i.e., $A(\tilde{r}_i)-A(\tilde{l}_i)=a_i >0$, then by constraints $(4a)$ and $(4b)$, we get $s_i=0$ and $kz_i^A \leq a_i \leq Kz_i^A$, this inequality is satisfied for $z_i^S=0$, $z_i^A=1$ and for positive real values of $k$ and $K$, sufficiently small and large. Further, $(4o)$ and $(4p)$ yield $(z_i^M, z_i^C, z_i^D, z_i^G,z_i^H)=(1,1,1,1,1)$, due to which rest of the constraints $(4c)-(4g)$ are obviously satisfied. Further, $S(\tilde{r}_i)=S(\tilde{l}_i)$ and $A(\tilde{r}_i)>A(\tilde{l}_i)$ implies $(3b)$ holds.
\end{enumerate}

The remaining cases can be similarly obtained. This establishes that the problems (P3) and (P4) are equivalent.\\
 Hence proved.\\
 
\section{Advantages of the proposed method}
Some of the main advantages of the proposed algorithm over the existing approaches are listed below:
\begin{enumerate}
\item The existing method \cite{ref12} can only be used to solve IVIFLPPs in which all the decision variables are taken to be non-negative crisp parameters. However, our proposed method can be employed successfully to handle IVIFLPPs having all the decision variables represented by unrestricted IVIFNs.
\item The existing study \cite{ref53} defined a new product operator and the basic arithmetic operations on unrestricted $LR$-type IFNs. But, there is no study on $LR$-type IVIFNs. Consequently, we have introduced the definition of $LR$-type IVIFNs and developed the arithmetic operations on unrestricted $LR$-type IVIFNs.
\item The existing models \cite{ref33, ref47, ref48, ref53} can be used only to deal with $LR$-type FLPPs and IFLPPs. But, in this article, we have solved the IVIFLPP in which all the parameters and decision variables are represented by $LR$-type IVIFNs, which is more general. Hence, the proposed method can be successfully reduced to solve IFLPPs as well as FLPPs.
\item In the present study, we have considered the model parameters and decision variables to be $LR$-type IVIFNs, as a result, the proposed method can also be utilized to solve LPPs where parameters and variables are represented by Triangular/Trapezoidal IVIFNs.
\item The proposed algorithm can be used to solve the models having some/all decision variables as unrestricted $LR$-type IVIFNs or non-negative $LR$-type IVIFNs.
\end{enumerate}
 
\section{Numerical illustration}
In this section, we present a numerical example to demonstrate the steps involved in the proposed algorithm.\\
Consider the following LPP having all the parameters as $LR$-type IVIFN and unrestricted crisp variables:
\begin{center}
$\hspace{-1.5cm}\textbf{(S1)}~~\max~\tilde{Z}=\tilde{5}\odot x_1 \oplus\tilde{8} \odot x_2$\\
$\mbox{s.t.}~~\tilde{12}\odot x_1 \oplus \tilde{4} \odot x_2 =\tilde{100},$\\
$~~~~~~~\tilde{6}\odot x_1 \oplus \tilde{10} \odot x_2 \preceq \tilde{150},$\\
$~~~~~~~~~~~~~~~~~~~~x_1~\mbox{and}~x_2$ are unrestricted in sign
\end{center}
where\\
$\tilde{5}=(5;2,2,3,3;5,5,5,4)_{LR},~~\tilde{8}=(8;1,1,2,2;4,4,2,3)_{LR},~~\tilde{12}=(12;2,3,4,4;6,8,4,4)_{LR},$\\
 $\tilde{4}=(4;1,1,2,2;4,4,2,2)_{LR},~~\tilde{6}=(6;3,4,4,4;6,6,4,4)_{LR},~~\tilde{10}=(10;3,4,4,5;6,8,5,5)_{LR},$\\
$\tilde{100}=(100;25,35,50,50;80,100,50,50)_{LR},~~\tilde{150}=(150;50,60,50,70;120,100,80,70)_{LR},$ \\
$L(x)=R(x)=L'(x)=R'(x)=\max \{0,1-x\}~\forall~x\in \mathbb{R}.$\\\\
\textbf{Solution:}
\begin{enumerate}
\item[\textbf{Step 1.}] Substituting the expressions of various parameters, the problem (S1) can be re-written as follows:\\

$\max~\tilde{Z}=(5;2,2,3,3;5,5,5,4)_{LR} \odot x_1 \oplus (8;1,1,2,2;4,4,2,3)_{LR} \odot x_2$\\
$~~~\mbox{s.t.}~~(12;2,3,4,4;6,8,4,4)_{LR}\odot x_1 \oplus (4;1,1,2,2;4,4,2,2)_{LR} \odot x_2 =(100;25,35,50,50;80,100,50,50)_{LR},$\\
$~~~~~~~~~~~(6;3,4,4,4;6,6,4,4)_{LR}\odot x_1 \oplus (10;3,4,4,5;6,8,5,5)_{LR} \odot x_2 \preceq (150;50,60,50,70;120,100,80,70)_{LR},$\\
$~~~~~~~~~~x_1~\mbox{and}~x_2$ are unrestricted in sign.\\

\item[\textbf{Step 2.}] Using the multiplication operation (Corollary 3.1.1) on $LR$-type IVIFNs along-with the fact that
$$\max\{a,b\}=\displaystyle\frac{1}{2}(a+b+|a-b|),$$
\noindent the problem (S1) is converted to the following equivalent problem:\\

\noindent \textbf{(S2)} $\max~\tilde{Z}=\big(5x_1;2|x_1|,2|x_1|,3|x_1|,3|x_1|;5|x_1|,5|x_1|,\displaystyle\frac{1}{2}(x_1+9|x_1|),\displaystyle\frac{1}{2}(-x_1+9|x_1|)\big)_{LR} \oplus$\\
$~~~~~~~~~~~~~~~~~ \big(8x_2;|x_2|,|x_2|,2|x_2|,2|x_2|;4|x_2|,4|x_2|, \displaystyle\frac{1}{2}(-x_2+5|x_2|),\displaystyle\frac{1}{2}(x_2+5|x_2|)\big)_{LR}$\\
$~~~~~~~~~~~\mbox{s.t.}~~\big(12x_1;\displaystyle\frac{1}{2}(-x_1+5|x_1|),\displaystyle\frac{1}{2}(x_1+5|x_1|),4|x_1|,4|x_1|;-x_1+7|x_1|,x_1+7|x_1|,4|x_1|,4|x_1|\big)_{LR}$\\
$~~~~~~~~~~~~~~~~~~~~~~~~~~\oplus\big(4x_2;|x_2|,|x_2|,2|x_2|,2|x_2|;4|x_2|,4|x_2|,2|x_2|,2|x_2| \big)_{LR}=(100;25,35,50,50;80,100,50,50)_{LR},$\\
$~~~~~~~~~~~~~~~~~~\big(6x_1;\displaystyle\frac{1}{2}(-x_1+7|x_1|),\displaystyle\frac{1}{2}(x_1+7|x_1|),4|x_1|,4|x_1|;6|x_1|,6|x_1|,4|x_1|,4|x_1|\big)_{LR} \oplus \big(10x_2;\displaystyle\frac{1}{2}(-x_2+$\\
$~~~~~~~~~~~~~~~~~~~7|x_2|),\displaystyle\frac{1}{2}(x_2+7|x_2|),\displaystyle\frac{1}{2}(-x_2+9|x_2|),\displaystyle\frac{1}{2}(x_2+9|x_2|);-x_2+7|x_2|,x_2+7|x_2|,5|x_2|,5|x_2| \big)_{LR}$\\ $~~~~~~~~~~~~~~~~~~\preceq (150;50,60,50,70;120,100,80,70)_{LR},$\\
$~~~~~~~~~~~~~~~~~~~x_1~\mbox{and}~x_2$ are unrestricted in sign.\\

\item[\textbf{Step 3.}] By employing the lexicographic ordering as in Definition $3.2.3$, using the addition operation, equality between $LR$-type IVIFNs and the corresponding definitions of $S, A, M, C, D, G ~\mbox{and}~ H$, we get the following non-linear programming problem:\\
\noindent \textbf{(S3)} lex max $\bigg(\displaystyle\frac{1}{8}(x_1-x_2),\displaystyle\frac{1}{8}(79x_1+129x_2), 5x_1+8x_2, 5x_1+8x_2-2|x_1|-|x_2|,5x_1+8x_2-3|x_1|-2|x_2|,$\\
$~~~~~~~~~~~~~~~~~~~~~~~\displaystyle\frac{9}{2}x_1+\displaystyle\frac{17}{2}x_2-\displaystyle\frac{9}{2}|x_1|-\displaystyle\frac{5}{2}|x_2|,5x_1+8x_2-5|x_1|-4|x_2|\bigg) $\\
$~~~~~~~~~~~~~~~~\mbox{s.t.}~~~3x_1+x_2=25,~~-x_1+5|x_1|+2|x_2|=50,~~x_1+5|x_1|+2|x_2|=70,$\\
$~~~~~~~~~~~~~~~~~~~~~~~~2|x_1|+|x_2|=25,~~-x_1+7|x_1|+4|x_2|=80,~~x_1+7|x_1|+4|x_2|=100,$\\
$~~~~~~~~~~~~~~~~~~~~~~~~\bigg(\displaystyle\frac{x_1}{8},\displaystyle\frac{1}{8}\big(97x_1+164x_2\big),6x_1+10x_2,\displaystyle\frac{13}{2}x_1+\displaystyle\frac{21}{2}x_2-\displaystyle\frac{7}{2}|x_1|-\displaystyle\frac{7}{2}|x_2|, 6x_1+\displaystyle\frac{21}{2}x_2-4|x_1|-$\\
$~~~~~~~~~~~~~~~~~~~~~~~~\displaystyle\frac{9}{2}|x_2|,6x_1+10x_2-4|x_1|-5|x_2|,6x_1+11x_2-6|x_1|-7|x_2|\bigg)\\
~~~~~~~~~~~~~~~~~~~~~~~~\preceq_{lex} (7.5,300,150,100,100,70,30),$\\
$~~~~~~~~~~~~~~~~~~~~~~~~~x_1~\mbox{and}~x_2$ are unrestricted in sign.\\

\item[\textbf{Step 4.}] To convert the lexicographic constraint $\preceq_{lex}$ into its equivalent form, we introduce the binary variables $z_1^S,\\ z_1^A, z_1^M, z_1^C, z_1^D, z_1^G~\mbox{and}~z_1^H$ and for $k=10^{-4}$ and $K=1000$, we have the following set of constraints:\\
\noindent $kz_1^S \leq 7.5-\displaystyle\frac{1}{8}(x_1) \leq Kz_1^S, \hspace{10.45cm} (5a)$\\
\noindent $-Kz_1^S+kz_1^A \leq 300-\displaystyle\frac{1}{8}\big(97x_1+164x_2\big) \leq Kz_1^A,\hspace{7.35cm} (5b)$\\
\noindent $-K(z_1^S+z_1^A)+kz_1^M \leq 150-6x_1-10x_2 \leq Kz_1^M,\hspace{7.08cm} (5c)$\\
\noindent $-K(z_1^S+z_1^A+z_1^M)+kz_1^C \leq 100- \displaystyle\frac{13}{2}x_1-\displaystyle\frac{21}{2}x_2+\displaystyle\frac{7}{2}|x_1|+\displaystyle\frac{7}{2}|x_2|\leq Kz_1^C,\hspace{3.55cm} (5d)$\\
\noindent $-K(z_1^S+z_1^A+z_1^M+z_1^C)+kz_1^D \leq 100-6x_1-\displaystyle\frac{21}{2}x_2+4|x_1|+\displaystyle\frac{9}{2}|x_2| \leq Kz_1^D,\hspace{3.15cm} (5e)$\\
\noindent $-K(z_1^S+z_1^A+z_1^M+z_1^C+z_1^D)+kz_1^G \leq 70-6x_1-10x_2+4|x_1|+5|x_2| \leq Kz_1^G,\hspace{2.63cm} (5f)$\\
\noindent $-K(z_1^S+z_1^A+z_1^M+z_1^C+z_1^D+z_1^G)+kz_1^H \leq 30-6x_1-11x_2+6|x_1|+7|x_2| \leq Kz_1^H,\hspace{1.79cm} (5g)$\\
\noindent $z_1^S,z_1^A,z_1^M,z_1^C,z_1^D,z_1^G,z_1^H \in \{0,1\}, \hspace{9.23cm} (5h)$\\
$z_1^S \leq z_1^A \leq z_1^M \leq z_1^C \leq z_1^D \leq z_1^G \leq z_1^H.\hspace{8.9cm} (5i)$\\

\item[\textbf{Step 5.}] Finally, the equivalent mixed $0$ \textendash $~1$ lexicographic non-linear optimization problem obtained is as follows:\\
\textbf{(S4)} lex max $\bigg(\displaystyle\frac{1}{8}(x_1-x_2),\displaystyle\frac{1}{8}(79x_1+129x_2), 5x_1+8x_2,5x_1+8x_2-2|x_1|-|x_2|,5x_1+8x_2-3|x_1|-$\\
$~~~~~~~~~~~~~~~~~~~~~~~2|x_2|,\displaystyle\frac{9}{2}x_1+
\displaystyle\frac{17}{2}x_2-\displaystyle\frac{9}{2}|x_1|-\displaystyle\frac{5}{2}|x_2|,5x_1+8x_2-5|x_1|-4|x_2|\bigg)$\\
$~~~~~~~~~~~~~~~\mbox{s.t.}~~~3x_1+x_2=25,~~-x_1+5|x_1|+2|x_2|=50,~~x_1+5|x_1|+2|x_2|=70,$\\
$~~~~~~~~~~~~~~~~~~~~~~~~2|x_1|+|x_2|=25,~~-x_1+7|x_1|+4|x_2|=80,~~x_1+7|x_1|+4|x_2|=100,$\\
$~~~~~~~~~~~~~~~~~~~~~~~~\mbox{constraints}~(5a)-(5i),$\\
$~~~~~~~~~~~~~~~~~~~~~~~~x_1~\mbox{and}~x_2$ are unrestricted in sign. \\

\item[\textbf{Step 6.}] To solve the problem (S4), we start by optimizing the first component of the objective function along-with imposing all the constraints of problem (S4). As a result, we solve the single-objective non-linear programming problem (S4-1) given by:
\begin{center}
\textbf{(S4-1)} $\max~ S=\displaystyle\frac{1}{8}(x_1-x_2)~~~~~~~~~~~~~~~$\\
$~~~~~~~~~~~~~~\mbox{s.t.}~~\mbox{all the constraints of (S4)}.$
\end{center}

Here, we have solved all the crisp optimization models by using a software "LINGO$-17.0$" on a MacBook Air system with 1.8 GHz Dual-Core Intel Core i5 processor and 8 GB RAM. The optimal solution of (S4-1) gives the optimal objective value as $S=1.875$. Next, for optimizing the second component of objective function of (S4), we solve the problem (S4-2).
\begin{center}
\textbf{(S4-2)} $\max~ A=\displaystyle\frac{1}{8}(79x_1+129x_2)~~~~~$\\
$~~~~~~~~~~~~~\mbox{s.t.}~~~\mbox{all the constraints of (S4)},$\\
$~~~~~~~~~~~~\displaystyle\frac{1}{8}(x_1-x_2) \geq 1.875$.
\end{center}
 The optimal objective value of (S4-2) is $A=18.125$ and for optimizing the third objective, the problem (S4-3) is solved.
 \begin{center}
\textbf{(S4-3)} $\max~ M=5x_1+8x_2~~~~~~~~~~~~~~~$\\
$~~~~~~~~~~~~~\mbox{s.t.}~~~\mbox{all the constraints of (S4)},$\\
$~~~~~~~~~~~~\displaystyle\frac{1}{8}(x_1-x_2) \geq 1.875$,\\
$~~~~~~~~~~~~~~~~~~~~~~~~\displaystyle\frac{1}{8}(79x_1+129x_2) \geq 18.125.$
\end{center}
 The optimal solution of (S4-3) gives $M=10$. Further, the fourth objective is optimized by solving the problem (S4-4).
\begin{center}
\textbf{(S4-4)} $\max~ C=5x_1+8x_2-2|x_1|-|x_2|$\\
$~~~~~~\mbox{s.t.}~~~\mbox{all the constraints of (S4)},$\\
$~~~~~\displaystyle\frac{1}{8}(x_1-x_2) \geq 1.875,$\\
$~~~~~~~~~~~~~~~~~\displaystyle\frac{1}{8}(79x_1+129x_2) \geq 18.125,$\\
$5x_1+8x_2 \geq 10.$
\end{center}
 The optimal objective value of (S4-4) is $C=-15$ and for optimizing the fifth objective, we solve the problem (S4-5).
 \begin{center}
~~\textbf{(S4-5)} $\max~ D=5x_1+8x_2-3|x_1|-2|x_2|$\\
$~~~~~~\mbox{s.t.}~~~\mbox{all the constraints of (S4)},$\\
$~~~~\displaystyle\frac{1}{8}(x_1-x_2) \geq 1.875,$\\
$~~~~~~~~~~~~~~~~\displaystyle\frac{1}{8}(79x_1+129x_2) \geq 18.125,$\\
$5x_1+8x_2 \geq 10,~$\\
$~~~~~~~~~~~~~~~~~~~~~~~~~5x_1+8x_2-2|x_1|-|x_2| \geq -15.$
\end{center}
 The optimal solution of (S4-5) gives the optimal objective value as $D=-30$ and to optimize the sixth objective, the problem (S4-6) is solved.
 \begin{center}
~~~\textbf{(S4-6)} $\max~ G=\displaystyle\frac{9}{2}x_1+\displaystyle\frac{17}{2}x_2-\displaystyle\frac{9}{2}|x_1|-\displaystyle\frac{5}{2}|x_2|$\\
$\mbox{s.t.}~~~\mbox{all the constraints of (S4)},$\\
$\displaystyle\frac{1}{8}(x_1-x_2) \geq 1.875,~$\\
$~~~~~~~~~~~\displaystyle\frac{1}{8}(79x_1+129x_2) \geq 18.125,$\\
$5x_1+8x_2 \geq 10,~~~~~~$\\
$~~~~~~~~~~~~~~~~~~~~5x_1+8x_2-2|x_1|-|x_2| \geq -15,$\\
$~~~~~~~~~~~~~~~~~~~~~~5x_1+8x_2-3|x_1|-2|x_2| \geq -30.$
\end{center}
 The optimal solution of (S4-6) gives $G=-55$. Finally, the seventh objective is optimized by solving the problem (S4-7).
 \begin{center}
\textbf{(S4-7)} $\max~ H=5x_1+8x_2-5|x_1|-4|x_2|~~~~~~~~~$\\
$\mbox{s.t.}~~~\mbox{all the constraints of (S4)},~~~~~~$\\
$\displaystyle\frac{1}{8}(x_1-x_2) \geq 1.875,~~~~~~~~$\\
$~~~~\displaystyle\frac{1}{8}(79x_1+129x_2) \geq 18.125,$\\
$5x_1+8x_2 \geq 10,~~~~~~~~~~~~~$\\
$~~~~~~~~~~~~~5x_1+8x_2-2|x_1|-|x_2| \geq -15,$\\
$~~~~~~~~~~~~~~~5x_1+8x_2-3|x_1|-2|x_2| \geq -30,$\\
$~~~~~~~~~~~~~~~~~~~~\displaystyle\frac{9}{2}x_1+\displaystyle\frac{17}{2}x_2-\displaystyle\frac{9}{2}|x_1|-\displaystyle\frac{5}{2}|x_2| \geq -55.$
\end{center}

 The optimal objective value and optimal solution of problem (S4-7) are respectively equal to $H=-60,$ and $x_1=10, ~x_2=-5$.\\

\item[\textbf{Step 7.}] By putting the optimal solution values as obtained in step 6, we get the unique optimal IVIF value
 $$\tilde{Z}=(10;25,25,40,40;70,70,65,50)_{LR}.$$
 \end{enumerate}

\section{An application in production planning}

A bicycle manufacturing company produces two models of bicycles, namely, road bikes and mountain bikes. Steel alloy and rubber constitute the primary raw materials required in the production process of the main body of the bicycles. It is observed from the past experiences that no complete unit of a road bicycle can be manufactured if $1.5$ hours of skilled labourers are employed per unit. On the other hand, the company can't bear to spend more than $2.5$ hours of labour time per unit of the road bike as it will reduce the efficiency of the process. From the past data, it is estimated that complete manufacturing of each unit of a road bike requires about $1.75$ to $2.25$ hours of skilled labourers. It was also judged that the production curve of road bikes peaks near to $2$ hours of labour time per unit. Further, it is known from the judgement of the decision-maker that values of all the parameters follow linear variations and involve uncertainty as well as some inherent hesitation. In this context, all the information related to the skilled labour hours is summarized in Table \ref{table2} while the data referring to the requirement of raw material is tabulated in Table \ref{table3}.

The manager of the company has approximated that about $100$ hours of labour time is available for the production process. Further, due to the various uncontrollable factors, there occurs a fluctuation in the supply and consumption of Steel alloy and rubber in the manufacturing firm. Thus, the manager has estimated that around $300$ units of Steel alloy and nearly $120$ units of rubber will be available for this production run.  However, the manufacturing firm can purchase the additional required units of Steel alloy at a price of about $\$10$ or can sell leftover units at the same price. The company estimated the selling price per unit of road bike at a price nearly $\$80$ and that of mountain bike at about $\$120$. Finally, the firm manager wants to find the number of optimal units of Steel alloy need to be purchased or sold and units of the road bikes and mountain bikes be produced consequently so as to maximize the total profit.\\


\begin{table*}[width=1\linewidth,cols=5,pos=h]
\caption{Data related to the requirement of skilled labour hours per bicycle} \label{table2}
\begin{tabular*}{\tblwidth}{@{} CCCCC@{} }
\toprule
Type of bicycle & Time in which & Estimated time & Peak Production  & Time unacceptable \\
& no unit prepared & range to prepare & Time  &  by the company\\
& completely & each unit & & for each unit$^1$\\
\midrule
Road bike & \centering$1.5$ & $1.75$ to $2.25$ & $2$ & $2.5$\\
\midrule
Mountain bike & $3$ &  $3.75$ to $4.5$ & $4$ & $4.7$\\
\bottomrule
\end{tabular*}
\vspace{1ex}

 {\raggedright $^1$ \begin{footnotesize}
 {This represents the labour hours which, if employed for a unit of the bicycle then it will reduce the efficiency of the production run.}
 \end{footnotesize} \par}
\end{table*}

\begin{table*}[width=1\linewidth,cols=5,pos=h]
\caption{Raw material requirement (in Kilograms) per unit of bicycle} \label{table3}
\begin{tabular*}{\tblwidth}{@{} CCCCCC@{} }
\toprule
Type of bicycle & Type of &  Amount of material that & Estimated range & Amount of material & Units of material \\
& material & can't produce one & of material & giving maximum  & can't be utilized\\
& & complete frame of bicycle & per unit of bike & production rate & per bicycle$^2$ \\
\midrule
Road bike & Steel alloy & $5$ & $6$ to $8$ & $7$ & $8.5$\\
& Rubber & $1.5$ & $1.75$ to $2.3$ & $2$ & $2.7$\\
\midrule
Mountain bike & Steel alloy & $6.5$ &  $7.2$ to $8.9$ & $8$ & $10$\\
& Rubber & $3$ & $3.5$ to $4.7$ & $4$ & $5$\\
\bottomrule
\end{tabular*}
\vspace{1ex}

 {\raggedright $^2$ \begin{footnotesize}
 {It describes the units of the raw material which if used per frame of the bicycle then it will reduce the profit of the company.}
 \end{footnotesize} \par}
\end{table*}

\textbf{Problem formulation:} Since all the parameters of the problem are known to be in the estimated/uncertain form, therefore, it is more relevant to represent these estimated numbers by $LR$-type IVIFNs. Following, the data given in Tables \ref{table2} and \ref{table3}, the given estimated parameters can be presented as follows:
\begin{enumerate}[$1.$]
\item Labour time (hrs.):
\begin{enumerate}[$(a).$]
\item For Road bike: $~~~~~~~~~~~~\tilde{2}=(2;0.1,0.1,0.25,0.25;0.5,0.5,0.3,0.3)_{LR}.$
\item For Mountain bike: $~~~~~\tilde{4}=(4;0.2,0.2,0.25,0.5;1,0.7,0.5,0.6)_{LR}$.
\end{enumerate}
\item Raw material (Kgs.):
\begin{enumerate}[$(a).$]
\item For Road bike:
\begin{enumerate}[$(i)$.]
\item Steel alloy: $~~~~~\tilde{7}=(7;0.5,0.5,1,1;2,1.5,1.5,1)_{LR}.$
\item Rubber: $~~~~~~~~~~\tilde{2}=(2;0.1,0.2,0.25,0.3;0.5,0.7,0.3,0.5)_{LR}.$
\end{enumerate}
\item For Mountain bike:
\begin{enumerate}[$(i)$.]
\item Steel alloy: $~~~~~\tilde{8}=(8;0.5,0.5,0.8,0.9;1.5,2,1.2,1.5)_{LR}.$
\item Rubber: $~~~~~~~~~~\tilde{4}=(4;0.3,0.3,0.5,0.7;1,1,0.7,0.8)_{LR}.$
\end{enumerate}
\end{enumerate}
\item Availability:
\begin{enumerate}[$(a).$]
\item Labour time (hrs.): $~~~~~\tilde{100}=(100;8,8,10,10;20,22,15,15)_{LR}.$
\item Steel alloy (Kgs.): $~~~~~~\tilde{300}=(300;10,12,15,15;30,30,20,25)_{LR}.$
\item Rubber (Kgs.): $~~~~~~~~~~~\tilde{120}=(120;10,8,15,15;30,30,18,20)_{LR}.$
\end{enumerate}
\item Estimated Cost (\$):
\begin{enumerate}[$(a).$]
\item For Steel alloy: $~~~~~~~~~~~~\tilde{10}=(10;1,1.5,2,2;4,5,3,3.5)_{LR}.$
\item For Road bike: $~~~~~~~~~~~~~\tilde{80}=(80;5,5,7,7;10,10,8,9)_{LR}.$
\item For Mountain bike: $~~~~~\tilde{120}=(120;8,7,10,10;15,15,12,11)_{LR}.$
\end{enumerate}
\end{enumerate}

Based on this data and as per the statement of the problem, it can be formulated as follows:
\begin{center}
~~~~\textbf{(M1)}~$\max~\tilde{Z}=\tilde{80} \odot \tilde{x}_1 \oplus \tilde{120} \odot \tilde{x}_2 \oplus \tilde{10} \odot \tilde{y}_1$\\
$~~~~~~\mbox{s.t.}~\tilde{7} \odot \tilde{x}_1 \oplus \tilde{8} \odot \tilde{x}_2 \oplus \tilde{y}_1 =\tilde{300},$\\
$~~~\tilde{2} \odot \tilde{x}_1 \oplus \tilde{4} \odot \tilde{x}_2 \preceq \tilde{120},$\\
$~~~\tilde{2} \odot \tilde{x}_1 \oplus \tilde{4} \odot \tilde{x}_2 \preceq \tilde{100},$\\
$~~~~~~~~\tilde{x}_1, \tilde{x}_2 \succeq 0$, $~\tilde{y}_1$ unrestricted
\end{center}
\noindent where $~\tilde{x}_1=(x_1;\xi_{1L}^{\mu},\eta_{1L}^{\mu},\xi_{1U}^{\prime\mu},\eta_{1U}^{\prime\mu};\xi_{1L}^{\nu},\eta_{1L}^{\nu},\xi_{1U}^{\prime\nu},\eta_{1U}^{\prime\nu})_{LR}= \mbox{Number of units of road bike to be produced},$\\\\
$\tilde{x}_2=(x_2;\xi_{2L}^{\mu},\eta_{2L}^{\mu},\xi_{2U}^{\prime\mu},\eta_{2U}^{\prime\mu};\xi_{2L}^{\nu},\eta_{2L}^{\nu},\xi_{2U}^{\prime\nu},\eta_{2U}^{\prime\nu})_{LR}= \mbox{Number of units of mountain bike to be produced}$,\\\\
$\tilde{y}_1=(y_1;\pi_{1L}^{\mu},\theta_{1L}^{\mu},\pi_{1U}^{\prime\mu},\theta_{1U}^{\prime\mu};\pi_{1L}^{\nu},\theta_{1L}^{\nu},\pi_{1U}^{\prime\nu},\theta_{1U}^{\prime\nu})_{LR}= \mbox{Number of additional units of Steel alloy to be purchased or sold}$.\\

\noindent \textbf{Solution:} Substituting the values of the parameters, the problem (M1) can be further expressed as:\\

\noindent \textbf{(M2)}~max~$\tilde{Z}=(80;5,5,7,7;10,10,8,9)_{LR} \odot \tilde{x}_1 \oplus (120;8,7,10,10;15,15,12,11)_{LR} \odot \tilde{x}_2 \oplus$\\
$~~~~~~~~~~~~~~~~~~~~~~~~~~~(10;1,1.5,2,2;4,5,3,3.5)_{LR} \odot \tilde{y}_1$\\
$~~~~~~~~~~~~~\mbox{s.t.}~(7;0.5,0.5,1,1;2,1.5,1.5,1)_{LR} \odot \tilde{x}_1 \oplus (8;0.5,0.5,0.8,0.9;1.5,2,1.2,1.5)_{LR} \odot \tilde{x}_2 \oplus \tilde{y}_1$\\
$~~~~~~~~~~~~~~~~~~ =(300;10,12,15,15;30,30,20,25)_{LR},$\\
$~~~~~~~~~~~~~~~~~~~(2;0.1,0.2,0.25,0.3;0.5,0.7,0.3,0.5)_{LR} \odot \tilde{x}_1 \oplus (4;0.3,0.3,0.5,0.7;1,1,0.7,0.8)_{LR} \odot \tilde{x}_2 $\\
$~~~~~~~~~~~~~~~~~~~\preceq (120;10,8,15,15;30,30,18,20)_{LR},$\\
$~~~~~~~~~~~~~~~~~~~(2;0.1,0.1,0.25,0.25;0.5,0.5,0.3,0.3)_{LR} \odot \tilde{x}_1 \oplus (4;0.2,0.2,0.25,0.5;1,0.7,0.5,0.6)_{LR} \odot \tilde{x}_2$\\
$~~~~~~~~~~~~~~~~~~~ \preceq (100;8,8,10,10;20,22,15,15)_{LR},$\\
$~~~~~~~~~~~~~~~~~~~~\tilde{x}_1, \tilde{x}_2 \succeq 0$, $~\tilde{y}_1$ unrestricted.\\

\noindent Now, since the input data follows the linear trend, therefore, $L(x)=R(x)=L'(x)=R'(x)=\mbox{max} \{0,1-x \}~ \forall~x \in \mathbb{R}$.

\noindent Further, using "LINGO-17.0" software for solving the equivalent crisp model obtained after applying Steps $1-7$ of the proposed algorithm, we get the following optimal solution of problem (M1):\\
$\tilde{x}_1=(1.71;0,0,0,0;1.7,0,0.64,0)_{LR},~~\tilde{x}_2=(3.05;1.02,1.13,1.5,1.18;2.08,2.13,1.5,1.97)_{LR}~$ and\\
$\tilde{y}_1=(263.69;0,0,0,0;0,0,0,0)_{LR}$\\
with the optimal IVIF value of the objective function equals 
$$\tilde{Z}=(3138.94;410.44,569.18,735.66,723.92;1454.84,1669.39,1050.39,1229.97)_{LR}
$$

\subsection{Results and discussion}

\begin{figure*} 
\includegraphics[scale=0.11]{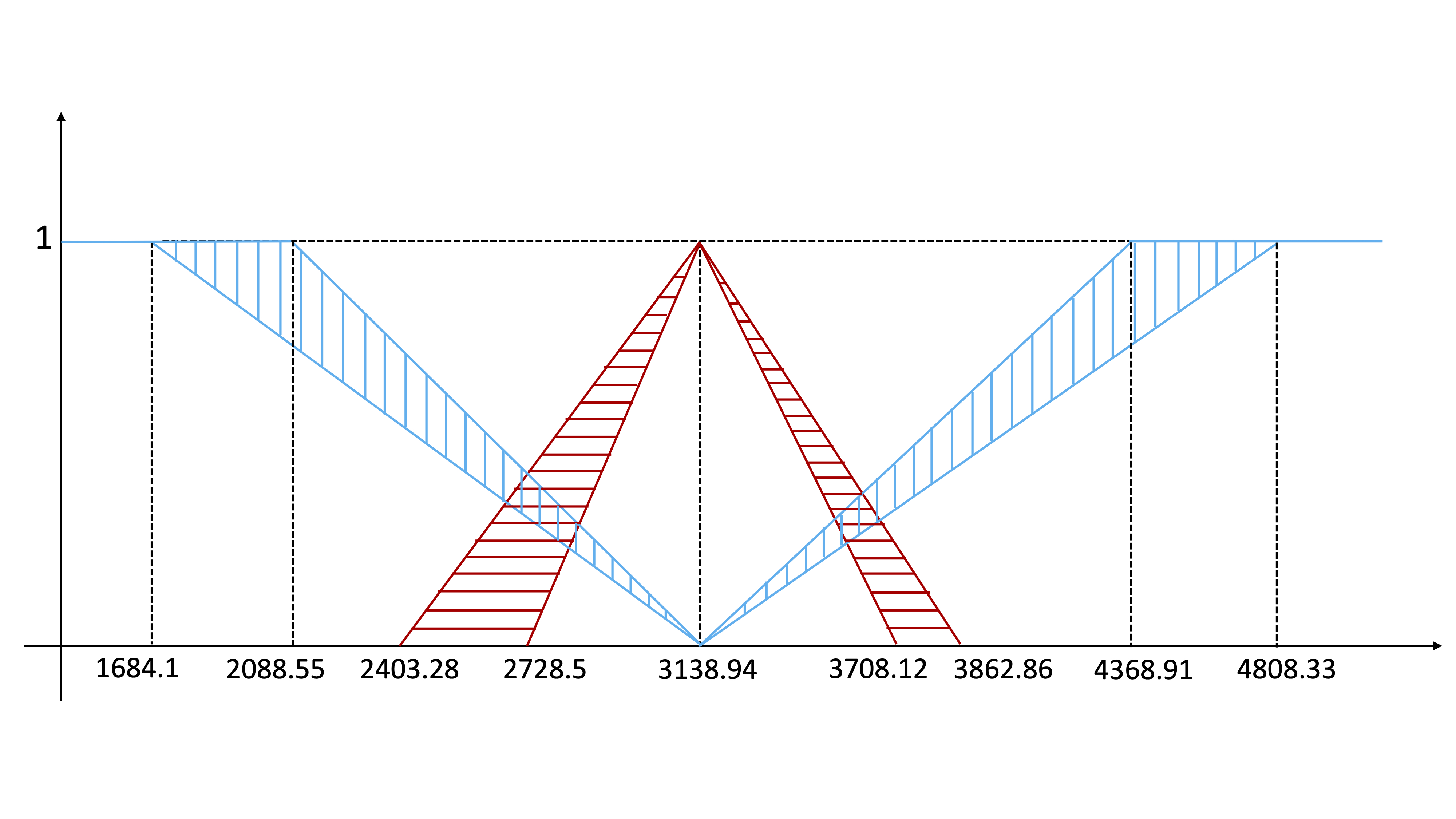}
\vspace{-1cm}
\caption{Representation of the optimal profit function as an $LR$-type TIVIFN} \label{fig4}
\end{figure*}

The optimal solution of problem (M1) calls for selling nearly $263$ Kgs. of leftover Steel alloy, manufacturing about $2$ units of the road bike and $3$ units of mountain bike to gain the maximum profit in the given scenario. The graphical representation of the objective function value $\tilde{Z}$ as an $LR$-type IVIFN is shown in Fig. \ref{fig4}. The interpretation of the profit function can be viewed as follows:\\
The company's acceptance increases if the profit value increases from nearly $\$2403.28$ to $\$3138.94$ while the degree of attainability of profit decreases if profit further increases from $\$3138.94$ to $\$3862.86$. The manager is fully satisfied at a profit of $\$3138.94$. However, when the profit increases from nearly $\$1684.1$ to $\$3138.94$, the degree of non \textendash attainability (or rejection) decreases continuously while the non \textendash attainability degree increases if the profit grows from $\$3138.94$ to $\$4808.33$.  Further, the company's profit can't go below $\$1684.1$ and a profit above $\$4808.33$ is also not achievable.

\subsection{Managerial insights}

Our modelling of linear optimization problems in $LR$-type IVIF environment integrates two significant variations in the data to solve the realistic problems. Firstly, this modelling allows different types of variation in the input data using the suitable choice of $L$ and $R$ functions. Secondly, the model parameters are taken as IVIFNs with interval degrees which handle the uncertain data in a most appropriate manner. Hence, it is crucial for a policy-maker to understand and judge how the optimal strategy varies using different $L,R$ functions and to evaluate the optimal solution which suits best to the concerned organization.

The formulation, solution and analysis of the production planning problem (M1) successfully incorporate the uncertain and vague data in the model and further provides the flexible and optimal production strategy to the company manager.

\subsection{Comparison with other cases}

The problem (M1) is also solved by transforming the inequality constraints of the model into equality by introducing the non-negative $LR$-type IVIF slack variables. Additionally, the model is alternatively solved using the usual order relation $\preceq$ in place of $\preceq_{lex}$. The values obtained are mentioned in Table \ref{table4}.

Further, by comparing the optimal values of $\tilde{Z}$, it is observed that the solution obtained by the proposed methodology yields better results. However, it may also be noticed that the optimal solutions are close to each other. Moreover, using the Definition $3.2.3$ for comparing the objective function values, we get
$$S(\tilde{Z}_{slack})=-75.34 < S(\tilde{Z}_{order~ relation~ \preceq})=-58.96 < S(\tilde{Z}_{proposed})=-30.89$$

\noindent implying that $\tilde{Z}_{slack} \prec \tilde{Z}_{order~ relation~ \preceq} \prec \tilde{Z}_{proposed}$. Therefore, it can be inferred that the proposed algorithm yields better results than the two alternative approaches. 

\begin{table*}[width=1.05\linewidth,cols=5,pos=h]
\caption{Optimal solution obtained using other cases and proposed algorithm} \label{table4}
\begin{tabular*}{\tblwidth}{@{} CLLL@{} }
\toprule
Solution & Proposed method & Adding non-negative $LR$-type  & Using order relation $\preceq $ \\
& & IVIF slack variables & instead of $\preceq_{lex}$\\
\midrule
$\tilde{x}_1$ & $(1.71;0,0,0,0;$ & $(0;0,0,0,0;$& $(0;0,0,0,0.1;$\\
& $~~1.7,0,0.64,0)_{LR}$& $~~0,0.3,0,0.875)_{LR}$ & $~~0,0.794,0,0.687)_{LR}$ \\
$\tilde{x}_2$ & $(3.05;1.02,1.13,1.5,1.18;$ & $(3.682;0,0,0,0;$& $(5.75;0,0,0,0;$\\
 & $~~2.08,2.13,1.5,1.97)_{LR}$& $~~0,0,0,0)_{LR}$ & $~~1.5,0,0.294,0)_{LR}$ \\
$\tilde{y}_1$ & $(263.69;0,0,0,0;$ & $(246.67;7,8.67,9.67,9;$&$(240;6.25,8.25,9,8.25;$ \\
 & $~~0,0,0,0)_{LR}$& $~~16.67,20,12,15)_{LR}$ & $~~9,18.25,9,15.25)_{LR}$ \\
$\tilde{Z}$ & $(3138.94;410.44,569.18,735.66,723.92;$& $(2908.54;339.13,495.48,607.52,638.16;$&$(3090;342.25,495.12,609.5,645.2;$ \\
 & $~~1454.84,1669.39,1050.39,1229.97)_{LR}$& $~~1141.93,1615.58,868.19,1184.22)_{LR}$ & $~~1257.75,1631.46,883.75,1170.27)_{LR}$ \\
\bottomrule
\end{tabular*}
\end{table*}

\section{Conclusions and future research scope}

The study proposes the definition of $LR$-type interval-valued intuitionistic fuzzy numbers and defines the basic operation on unrestricted $LR$-type IVIFNs with the help of $\alpha$-cut and $\beta$-cut. We have also derived the expressions of the various ranking indices for these numbers. Further, a methodology has been proposed to solve a class of unrestricted fully $LR$-type  IVIF linear programming problems. The current study also generalizes the results and theory of fuzzy, intuitionistic fuzzy and $LR$-type fuzzy / intuitionistic fuzzy numbers. Finally, it is to be pointed out that all the existing models of LPPs \cite{ref33, ref42, ref47, ref48, ref53} can also be solved using the proposed algorithm. However, many real-world problems which fail to have crisp parameters can only be solved efficiently using the proposed technique by representing the uncertain data using intervals. Moreover, for ranking of the $LR$-type IVIFNs, a lexicographic criteria $(S, A, M, C, D, G, H)$ has been used. However, the choice of this ranking criterion is not fixed. There are a total of $7!$ permutations possible. It totally depends on the decision-maker's attitude towards the preferences of various parameters involved in the lexicographic ranking. Furthermore, the practical applicability of the proposed model and solution algorithm is demonstrated by solving the production planning problem of a manufacturing company. Additionally, the optimal solution values are also compared with the two alternative approaches that can be used to solve the problem. The comparative results conclude that the proposed technique perform better than the alternate methods.\\

Some of the future directions of the present work are as follows: 
\begin{enumerate}
\item It can be observed that although the method is able to handle a very generalized class of linear optimization problems under uncertain conditions but the algorithm involves a large number of computation steps. Since, to obtain the final optimal value, one needs to solve seven mixed $0-1$ integer non-linear programming problems. Therefore, it will be interesting to devise a more computationally efficient approach to handle such problems.
\item In the future, different optimization problems such as supply chain problems, portfolio optimization problems, transportation problems, etc. can be solved under the $LR$-type IVIF environment.
\item An important futuristic research aspect would be to formulate and to devise a solution algorithm for multi-objective LPPs under IVIF scenario.
\item In future research, the methodology can be extended to solve the quadratic programming problems, non-linear / fractional problems having parameters / variables as $LR$-type IVIFNs. 
\end{enumerate}


\printcredits  

\section*{Declaration of competing interest}
We wish to confirm that there are no known conflicts of interest associated with this article and there has been no significant financial support for this work that could have influenced its outcome.

\section*{Funding}
This research did not receive any specific grant from funding agencies in the public, commercial, or not-for-profit sectors.

\section*{Data availability}
No data was used for the research described in the article.

\section*{Acknowledgement} 
The first author is thankful to the Ministry of Human Resource Development, India, for financial support, to carry out this research work.







\end{document}